\documentclass[11pt]{article}
\pdfoutput=1

\usepackage[T1]{fontenc}
\usepackage{amsmath}
\usepackage{amssymb}
\usepackage{epsfig}
\usepackage{graphicx}
\usepackage{booktabs}
\usepackage{theorem}

\newtheorem{theorem}{Theorem}[section]

\newtheorem{proposition}[theorem]{Proposition}

\theorembodyfont{\upshape}

\newtheorem{assumption}[theorem]{Assumption}
\newtheorem{example}[theorem]{Example}
\newtheorem{remark}[theorem]{Remark}

\newcommand{\note}[1]
{$^{(!)}$\marginpar[{\hfill\tiny{\sf{#1}}}]{\tiny{\sf{(!) #1}}}}

\newcommand{\proofbegin}{{\sc Proof. }}

\newcommand{\proofend}{\hfill{\fboxsep0mm%
   \fbox{\rule{0mm}{1.2ex}\hspace{1.2ex}}}\par\medskip}
\newcommand{\rnum}{{\mathbf R}}

\newcommand{\nnum}{{\mathbf N}}

\newcommand{\znum}{{\mathbf Z}}
\newcommand{\mean}{{\mathbf E}}

\newcommand{\var}{{\rm var}}
\newcommand{\cov}{{\rm cov}}

\newcommand{\tilX}{\tilde{X}}

\newcommand{\pconv}{\stackrel{P}{\rightarrow}}
\newcommand{\dconv}{\stackrel{{\cal D}}{\rightarrow}}

\newcommand{\argmin}{\mathop{\rm argmin}}

\newcounter{alphcount}
{\begin{list}{(\alph{alphcount})\hfill}%
             {\usecounter{alphcount}\labelwidth2em%
              \leftmargin2em\labelsep0em\topsep0em%
              \itemsep0em plus 0.5ex\parsep0em}}{\end{list}}
\newcounter{romancount}
\newenvironment{romanlist}%
{\begin{list}{(\roman{romancount})\hfill}%
             {\usecounter{romancount}\labelwidth2em%
              \leftmargin2em\labelsep0em\topsep0em%
              \itemsep0em plus 0.5ex\parsep0em}}{\end{list}}
\newcounter{storecount}

 \def\doublespace{\baselineskip=\normalbaselineskip \multiply\baselineskip by 6
 \divide\baselineskip by 5}

\begin{document}
\doublespace
\begin{titlepage}
\title{\bf Locally Stationary Processes}
\bigskip
\bigskip

\author{Rainer Dahlhaus\\
    Institut f\"{u}r Angewandte Mathematik\\
    Universit\"{a}t Heidelberg\\
    Im Neuenheimer Feld 294\\
    69120 Heidelberg\\
    Germany}
\date{}
\maketitle

\vspace*{1.5cm}

\begin{abstract}
\noindent The article contains an overview over locally stationary processes. At the beginning time varying autoregressive processes are discussed in detail - both as as a deep example and an important class of locally stationary processes. In the next section a general framework for time series with time varying finite dimensional parameters is discussed with special emphasis on nonlinear locally stationary processes. Then the paper focusses on linear processes where a more general theory is possible. First a general definition for linear processes is given and time varying spectral densities are discussed in detail. Then the Gaussian likelihood theory is presented for locally stationary processes. In the next section the relevance of empirical spectral processes for locally stationary time series is discussed. Empirical spectral processes play a major role in
proving theoretical results and provide a deeper understanding of
many techniques. The article concludes with an overview of other results for locally stationary processes.

\end{abstract}
\setcounter{page}{0}

\vfill

\noindent {\small {\bf Acknowledgement:} I am grateful to Suhasini Subba Rao for helpful comments on an earlier version which lead to significant improvements.}
\end{titlepage}

\section{Introduction} \label{intro}

Stationarity has played a major role in time series analysis for several decades. For stationary processes there exist a large variety of models and powerful methods, such as bootstrap methods or methods based on the spectral density. Furthermore, there are important mathematical tools such as the ergodic theorem or several central limit theorems. As an example we mention the likelihood theory for Gaussian processes which is well developed.

During recent years the focus has turned to nonstationary time series. Here the situation is
more difficult: First, there exists no natural generalization from stationary to nonstationary time series and second, it is often not clear how to set down a meaningful asymptotics for
nonstationary processes. An exception are nonstationary models which are generated by a time invariant generation mechanism -- for examples integrated or cointegrated models. These models have attracted a lot of attention during recent years. For general nonstationary processes ordinary asymptotic considerations are often contradictory to the idea of nonstationarity since future observations of a nonstationary process may not contain any information at all on
the probabilistic structure of the process at present. For this reason the theory of locally stationary processes is based on infill asymptotics originating from nonparametric statistics.

As a consequence valuable asymptotic concepts such as consistency,
asymptotic normality, efficiency, LAN-expansions, neglecting higher order
terms in Taylor expansions, etc. can be used in the theoretical
treatment of statistical procedures for such processes. This leads to several meaningful results also for the original non-rescaled case such as the comparison of different estimates, the approximations for
the distribution of estimates and bandwidth selection (for a detailed example see Remark~\ref{ImplicationsNonRescaled}).

The type of processes which can be described with this infill asymptotics are processes which locally at each time point are close to a stationary process but whose characteristics (covariances, parameters, etc.) are gradually changing in an unspecific way as time evolves. The simplest example for such a process may be an AR(p)-process whose parameters are varying in time. The infill asymptotic approach means that time is rescaled to the unit interval. For time varying AR-processes this is explained in detail in the next section. Another example are GARCH-processes which have recently been investigated by several authors -- see Section~\ref{NonlinearProcAndDerivative}.

The idea of having locally approximately a stationary process was also the starting point of Priestley's theory of processes with evolutionary spectra
(Priestley (1965) -- see also Priestley (1988), Granger and Hatanaka (1964), Tj{\o}stheim (1976) and M\'{e}lard and Herteleer-de-Schutter (1989) among others).
Priestley considered processes having a time varying spectral
representation
\begin{equation*}
X_{t} = \int^{\pi}_{-\pi} \mbox{exp} (i \lambda t) \, \tilde{A}_{t} (\lambda) \,d \xi
(\lambda), \quad t \in \znum
\label{1.1}
\end{equation*}
with an orthogonal increment process $\xi (\lambda)$ and a time varying
transfer function $\tilde{A}_{t}(\lambda)$. (Priestley mainly looked at continuous
time processes, but the theory is the same).
Also within this approach asymptotic considerations (e.g. for
judging the efficiency of a local covariance estimator) are not possible
or meaningless from an applied view. Using the above mentioned infill asymptotics means in this case basically to replace $\tilde{A}_{t} (\lambda)$ with some function $A (t/T,\lambda)$ -- see (\ref{tvSpectralRepresentation}).

Beyond the above cited references on processes with evolutionary spectra there has also been work on processes with time varying parameters which does not use the infill asymptotics discussed in this paper (cf. Subba Rao (1970); Hallin (1986) among others). Furthermore, there have been several papers on inference for processes with time varying parameters -- mainly within the engineering literature (cf. Grenier (1983), Kayhan et.al. (1994) among others).

The paper is organized as follows: In Section~\ref{DeepExample}
we start with time varying autoregressive processes as a deep example and an important class of locally stationary processes. There we mark
many principles and problems addressed at later stages with higher generality. In Section~\ref{NonlinearProcAndDerivative} we present a more general framework for time series with time varying finite dimensional parameters and show how nonparametric inference can be done and theoretically handled. We also introduce derivative processes which play a major role in the derivations. The results cover in particular nonlinear processes such as GARCH-processes with time varying parameters.

If one restrict to linear processes or even more to Gaussian processes then a much more general theory is possible which is developed in the subsequent sections. In Section~\ref{SectionSpectralDensity} we give a general definition for linear processes and discuss time varying spectral densities in detail. Section~\ref{SectionLikelihood} then contains the Gaussian likelihood theory for locally stationary processes. In Section~\ref{SectionEmpiricalSpectralProcesses} we discuss the relevance of  empirical spectral processes for locally
stationary time series. Empirical spectral processes play a major role in
proving theoretical results and provide a deeper understanding of
many techniques.

\section{Time varying autoregressive processes -- a deep example} \label{DeepExample}

We now discuss time varying autoregressive processes in detail. In particular we mark
many principles and problems addressed at later stages with higher generality. Consider the time varying AR(1) process
\begin{equation} \label{AR1-nonrescaled}
X_{t} + \alpha_t X_{t-1} = \sigma_t \, \varepsilon_{t} \qquad \mbox{ with }
\varepsilon_{t} \mbox{ iid } {\cal N}(0,1).
\end{equation}
We now apply infill asymptotics that is we rescale the parameter curves $\alpha_t$
and $\sigma_t$ to the unit interval. This means that we replace them by $\alpha (\frac {t} {T})$ and  $\sigma
(\frac {t} {T})$ with curves $\alpha (\cdot): [0,1] \rightarrow (-1,1)$ and $\sigma(\cdot):[0,1]
\rightarrow (0,\infty)$ leading in the general AR(p)-case to the definition given in
(\ref{ARpRescaled}) below. Formally this results in replacing $X_{t}$ by a triangular array of
observations $\big( X_{t,T}\, ; \; t=1,\ldots,T; \; T \in \mathbb{N} \big)$ where $T$ is the sample
size.

We now indicate again the reason for this rescaling. Suppose we fit the
parameteric model $\alpha_{\theta,t} := b + ct+dt^{2}$ to the nonrescaled model (\ref{AR1-nonrescaled}) which we assume to be observed for $t = 1, ..., T$. It is easy to construct different
estimators for the parameters (e.g. the least squares estimator, the  maximum likelihood estimator
or a moment estimator) but it is nearly impossible to derive the finite sample properties of these
estimators. On the other hand classical non-rescaled asymptotic considerations for comparing these
estimators make no sense since with $t \rightarrow \infty$ also $\alpha_{\theta,t} \rightarrow \infty$
while e.g. $|\alpha_t|$ may be less than one within the observed segment -- i.e. the resulting asymptotic
results are without any relevance for the observed stretch of data. By rescaling $\alpha_t$ and
$\sigma_t$ to the unit interval as described above we overcome these problems. As $T$ tends to
infinity more and more observations of each local structure become available and we obtain a
reasonable framework for a meaningful asymptotic analysis of statistical procedures allowing to
retain such powerful tools as consistency, asymptotic normality, efficiency, LAN-expansions, etc.
for nonstationary processes. For example the results on asymptotic normality of an estimator
obtained in this framework may be used to approximate the distribution of the estimator in the
finite sample situation. It is important to note that classical asymptotics for stationary
processes arises as a special case of this infill asymptotics in case where all parameter curves are
constant.

Unfortunately infill asymptotics does not describe the physical
behavior of the process as $T \rightarrow \infty$. This may be unusual for time series analysis but it
has been common in other branches of statistics for many years. We remark that all statistical
methods and procedures stay the same or can easily be translated from the rescaled processes to the
original non-rescaled processes. A more complicated example on how the results of the rescaled case transfer to
the non-rescaled case is given in Remark~\ref{ImplicationsNonRescaled}.

In the following we therefore consider time varying autoregressive \big(tvAR(p)\big) processes
defined by
\begin{equation} \label{ARpRescaled}
X_{t,T} + \sum_{j=1}^{p}\alpha_{j}(\frac{t}{T})\; X_{t-j,T}=\sigma(\frac{t}{T})\;
\varepsilon_{t} ,\mbox{ } t\in\znum
\end{equation}
where the $\varepsilon_{t}$ are independent random variables with mean zero
and variance 1. We assume $\sigma(u)=\sigma(0)$, $\alpha_j(u)=\alpha_j(0)$ for $u<0$ and
$\sigma(u)=\sigma(1)$, $\alpha_j(u)=\alpha_j(1)$ for $u>1$. In addition we usually assume some smoothness conditions on $\sigma(\cdot)$ and the $\alpha_{j}(\cdot)$. In addition one may include a time varying mean by replacing $X_{t-j,T}$ in (\ref{ARpRescaled}) by $X_{t-j,T}-\mu(\frac{t-j}{T})$ -- see Section~\ref{SectionAdditionalTopics}.6.

In some neighborhood of a fixed time point $u_0=t_0/n$ the process $X_{t,T}$ can be approximated by
the stationary process $\tilde{X}_t (u_0)$ defined by
\begin{equation} \label{ARpApproximation}
\tilde{X}_{t} (u_0) + \sum_{j=1}^{p}\alpha_{j}(u_0)\; \tilde{X}_{t-j} (u_0)=\sigma(u_0)\,
\varepsilon_{t} ,\mbox{ } t\in\znum.
\end{equation}
It can be shown (see Section~\ref{NonlinearProcAndDerivative}) that we have under suitable regularity conditions
\begin{eqnarray} \label{AR-LocalApproximation}
\big|X_{t,T} - \tilde{X}_{t}(u_0)\big| = O_{p}\Big(\Big|\frac{t}{T}-u_0\Big| + \frac{1}{T} \Big)
\end{eqnarray}
which justifies the notation ``locally stationary process''. $X_{t,T}$ has an unique time varying spectral density which is locally the same
as the spectral density of $\tilde{X}_t (u)$, namely
\begin{equation} \label{tvARSpectralDensity}
f(u, \lambda) :=
{\sigma^2 (u) \over 2\pi} \ \Big|1 + \sum_{j=1}^p\alpha_j(u)\exp (-ij\lambda)\Big|^{-2}
\end{equation}
(see Example~\ref{ExampleLocStat}). Furthermore it has locally in some sense the same autocovariance
\begin{equation*} \label{ARAutoCovariance}
c(u,j) := \int _{-\pi}^\pi e^{ij\lambda} f(u, \lambda) d\lambda,
\hspace {.4cm} j\in {\bf Z}
\end{equation*}
since $\cov (X_{[uT],T}, X_{[uT] + k,T}) = c(u,k) + O(T^{-1})$ uniformly in $u$ and $k$ (cf.(\ref{tvCovarianceProperty})). This justifies to term $c(u, k)$ the local covariance function of $X_{t,T}$ at time $u=
t/T$.

\begin{figure}
\centering
\includegraphics[width=340pt,keepaspectratio]{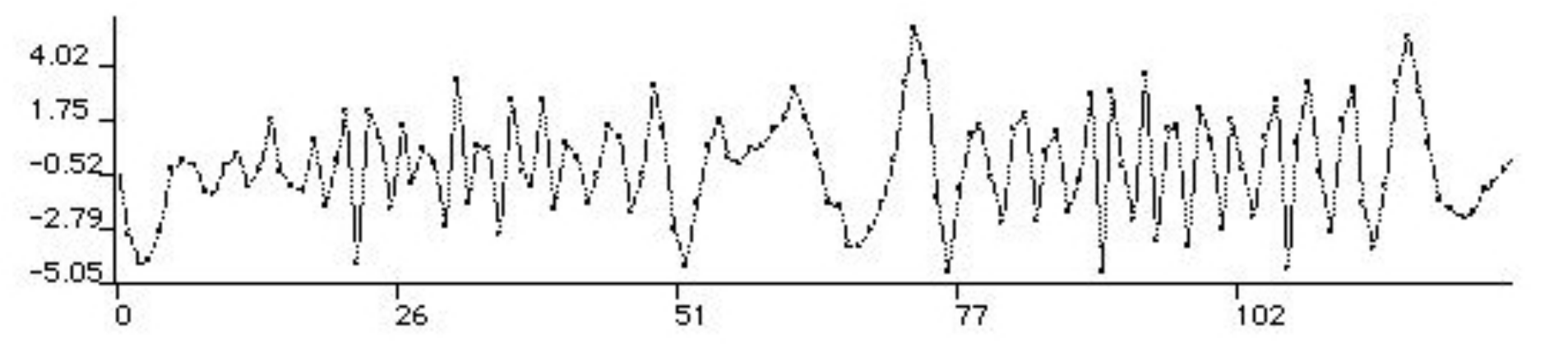}
\caption{\footnotesize T=128 realizations of a time varying AR(2)-model}
\label{fig1}
\end{figure}

\begin{figure}
\centering
\includegraphics[width=400pt,keepaspectratio]{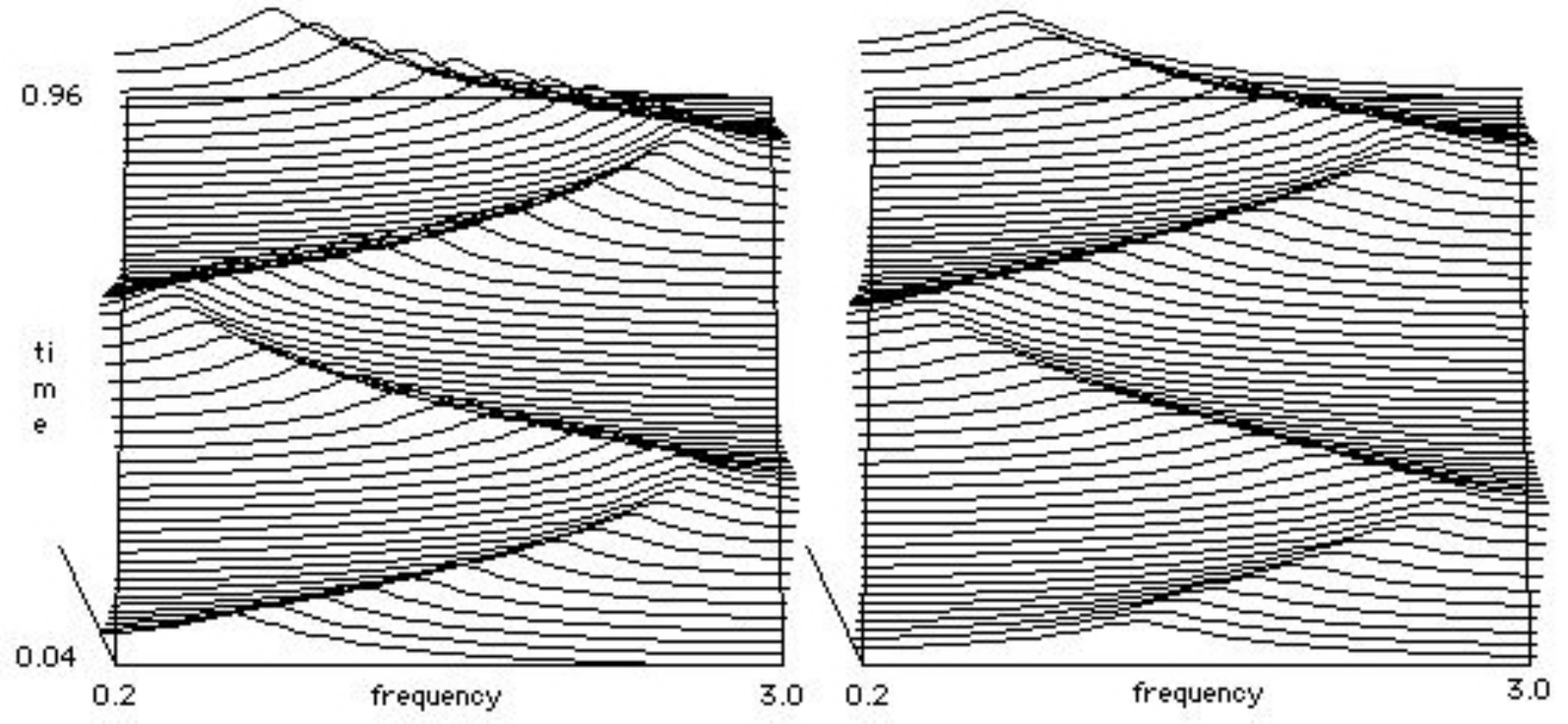}
\caption{\footnotesize True and estimated time varying spectrum of a tvAR(2)-process}
\label{fig2}
\end{figure}

As an example Figure~\ref{fig1} shows $T=128$ observations of a tvAR(2)-process with mean $0$ and parameters
$\sigma(u)\equiv 1$, $\alpha_1(u)\equiv-1.8\cos(1.5-\cos4{\pi} u)$, $\alpha_2(u)=0.81$ and Gaussian innovations $\varepsilon_t$. The parameters are chosen in a way such that for fixed $u$ the complex roots of the
characteristic polynomial are $\frac{1}{0.9}\exp[\pm i(1.5-\cos 4{\pi} u)]$, that is they are close to the unit circle and their phase varies cyclically with $u$.
As could be expected from these roots the observations show a periodic behavior with time varying
period-length. The left picture of Figure~\ref{fig2} shows the true time
varying spectrum of the process. One clearly sees that the location of the peak is also time varying (it is located at frequency $1.5\!-\!\cos 4{\pi} u$).

\medskip
\underline{1. Local estimation by stationary methods on segments}

 An ad-hoc method which works in
nearly all cases for locally stationary processes is to do inference via stationary methods on segments. The idea is that the process $X_{t,T}$ is almost
stationary  on a reasonably small segment $\{t : |t/T-u_0|\leq b/2 \}$. The parameter of interest (or
the correlation, spectral density, etc) is estimated by some classical method and the resulting
estimate is assigned to the midpoint $u_0$ of the segment. By shifting the segment this finally leads
to an estimate of the unknown parameter curve (time varying correlation, time varying spectral
density, etc). An important modification of this method is obtained when more weight is
put on data in the center of the interval than at the edges. This can often be achieved by
using a data taper on the segment or by using a kernel type estimate.

Since we use observations from the process $X_{t,T}$ (instead of $\tilde{X}_{t}(u_0)$) the procedure
causes a bias which depends on the degree of non-stationarity of the process on the segment. It is
possible to evaluate this bias and to use the resulting expression for an optimal choice of the
segment length. To demonstrate this we now discuss the estimation of the  AR coefficient functions by
classical Yule-Walker estimates on segments. Since the approximating process $\tilde{X}_{t}(u_0)$ is
stationary we obtain from (\ref{ARpApproximation})  that the Yule-Walker equations hold locally at
time $u_0$, that is we have with $\boldsymbol{\alpha}(u_0):= \big( \alpha_1(u_0) ,...,\alpha_p(u_0)\big)'$
\begin{equation} \label{TheoreticalYuleWalker}
\boldsymbol{\alpha}(u_0) = -  R(u_0)^{-1} \, r(u_0) \qquad \mbox{and} \qquad \sigma^2 (u_0)= c(u_0, 0) + \boldsymbol{\alpha}(u_0)' \,r (u_0)
\end{equation}
where $r(u_0):= \big(
c(u_0,1),...,c(u_0,p)\big)'$ and
  $R(u_0):=\{ c(u_0,i-j)\}_{i,j=1,...,p}$.

To estimate $\boldsymbol{\alpha}(u_0)$ we use the classical Yule Walker estimator  on the segment $[u_0T] -
N/2 + 1, \ldots , [u_0T] + N/2$ (ordinary time) or on $[u_0-b_T/2, u_0+ b_T/2]$ (rescaled time with
bandwidth $b_T := N/T$), that is
\begin{equation} \label{LocalYuleWalkerEst}
\hat{\boldsymbol{\alpha}}_T (u_0) = - \hat{R}_T(u_0)^{-1} \,
\hat{r}_T (u_0) \qquad \mbox{and} \qquad \hat{\sigma}^2_{T}(u_0) =\hat{c}_T(u_0, 0)+\hat{\boldsymbol{\alpha}}_T (u_0)' \,
\hat{r}_T (u_0)
\end{equation}
where $\hat{r}_T (u_0):= (\hat{c}_T(u_0,1),...,\hat{c}_T(u_0,p))'$ and $\hat{R}_T (u_0):=\{\hat{c}_T(u_0,i-j)\}_{i,j=1,...,p}$ with some covariance estimator $\hat{c}_T(u_0,j)$.

Before we discuss the properties of this estimator we first discuss different
covariance estimates and their properties.

\medskip
\underline{2. Local covariance estimation}

The covariance estimate with data taper on
the segment $[u_0T] \!- \!N/2 \!+ \!1, \ldots , [u_0T]\! + \!N/2$ is
\begin{equation} \label{TaperedCovarianceEstimate}
  {\hat c}_{T}(u_0,k):=
\frac{1}{H_N}\sum _{s,t=1 \atop s-t=j}^{N}h(\frac{s}{N})\, h(\frac{t}{N})\,
X_{[u_0T]-\frac{N}{2}+s ,T}\, X_{[u_0T]-\frac{N}{2}+t,T}.
\end{equation}
where $h:[0,1] \rightarrow {\bf R}$ is a data taper with $h(x)=h(1-x)$,
$H_N:=\sum_{j=0}^{N-1}h^2(\frac{j}{N} )\sim N \int_0^1 h^2(x)\,dx$ is the normalizing factor. The data taper usually is largest at $x=1/2$ and decays slowly to $0$ at the edges. For
$h(x) = \chi_{(0,1]} (x) $ we obtain the classical non-tapered covariance estimate.

An asymptotically equivalent (and from a certain viewpoint more intuitive estimator) is the kernel
density estimator
\begin{equation} \label{KernelCovarianceEstimate}
\tilde{c}_{T} (u_0,k) : = \frac{1}{b_{T}T} \sum_{t} K \left(
\frac{u_0-(t+k/2)/T}{b_{T}} \right) X_{t,T} X_{t+k,T}
\end{equation}
where $K: \rnum \rightarrow [0, \infty)$ is a kernel with $K(x) = K(-x)$, $\int K(x) dx = 1$,
$K(x) =0$ for $x \not\in [-1/2, 1/2]$ and $b_{T}$ is the bandwidth. Also equivalent is
\begin{equation} \label{KernelCovarianceEstimate2}
\tilde{\tilde{c}}_{T} (u_0,i,j) : = \frac{1}{b_{T}T} \sum_{t} K \left(
\frac{u_0-t/T}{b_{T}} \right) X_{t-i,T} X_{t-j,T}
\end{equation}
with $i-j=k$ which appears in least square regression -- cf. Example~\ref{ExampleNonlinearModels}(i). If $K(x) = h(x)^{2}$ all three estimators are equivalent in the sense that they lead to the same asymptotic bias, variance and mean squared error. For reasons of clarity a few remarks are in order:

1) The classical stationary method on a segment is in this case the estimator
without data taper which is the same as the kernel estimator with a rectangular kernel.\\
2) A first step towards a better estimate (as it is proved below) is to put higher weights
in the middle and lower weights at the edges of the observation domain in order to cope in a better
way with the nonstationarity of $X_{t,T}$ on the segment. In this context this may be either
achieved by using a kernel estimate or a data-taper which is asymptotically equivalent. This is
straightforward for local covariance estimates and local Yule-Walker estimates and can usually also be applied
to other estimation problems.\\
3) Data-tapers have also been used for stationary time series (in particular in spectral
estimation, but also with Yule Walker estimates and covariance estimation where they give positive
definite autocovariances with a lower bias). Thus the reason for using data-tapers for segment
estimates is twofold: reducing the bias due to nonstationarity on the segment and reducing the
(classical) bias of the procedure as a stationary method.

We now determine the mean squared error of the above estimators. Furthermore, we determine the
optimal segment length $N$ and show that weighted estimates are  better than ordinary
estimates.

\begin{theorem}
Suppose $X_{t,T}$ is locally stationary with mean $0$. Under suitable regularity conditions (in particular second order smoothness of $c (\cdot,k)$) we
have for $\hat{c}_{T}(u_0,k)$, $\tilde{c}_{T}(u_0,k)$ and $\tilde{\tilde{c}}_{T}(u_0,i,j)$ with $K(x) = h(x)^{2}$ and $b_T = N/T$
\begin{flalign*} \label{}
(i) \qquad \;  &\mean \, \hat{c}_{T}(u_0,k) = c (u_0,k) + \frac{1}{2} \, b_{T}^{2} \int x^{2} K(x) \,
dx \left[ \frac{\partial^{2}}{\partial^{2}u} \, c(u_0,k) \right] + o
(b^{2}_{T}) + O \!\left( \frac{1}{b_{T}T} \right) & \,
\end{flalign*}
and
\begin{flalign*} \label{}
(ii) \;\;&{\rm var} \big(\hat{c}_{T}(u_0,k)\big) = \frac{1}{b_{T}T} \int^{1/2}_{-1/2}
K(x)^{2} dx \sum^{\infty}_{\ell =
- \infty} c(u_0,\ell) \, \big[c(u_0,\ell) + c (u_0, \ell + 2k)\big] + o \!\left( \frac{1}{b_{T}T} \right). & \,
\end{flalign*}
\end{theorem}
\proofbegin (i) see Dahlhaus (1996c), (ii) is omitted (the form of the asymptotic variance is the same as in the stationary case).

\medskip

Note that the above bias of order $b_{T}^{2}$ is solely due to
nonstationarity which is measured by $\frac{\partial^{2}}{\partial u^{2}} c(u_0,k)$. If the process
is stationary this second derivative is zero and the bias disappears. The bandwidth $b_{T}$ may now be chosen to minimize the mean squared error.

\begin{remark} [Minimizing the mean squared error] \label{MinimalMSE} $\,$\\
Let $\mu (u_0) :=  \frac{\partial^{2}}{\partial^{2}u_0} \, c(u_0,k)$,  $\tau (u_0) := \sum^{\infty}_{\ell =
- \infty} c(u_0,\ell) \, \big[c(u_0,\ell) + c (u_0, \ell + 2k)\big]$, $d_K :=\int x^{2} K(x) \, dx$ and $v_K:=\int K(x)^{2}\, dx$. Then we have for the mean squared error
\begin{equation} \label{AsymptoticMSE}
\mean \, \big|\hat{c}_{T} (u_0,k) - c (u_0,k) \big|^{2} = \frac {b^{4}} {4} \,d_K^{2}\,
\mu(u_0)^{2} + \frac{1}{bT}\, v_K \, \tau(u_0) +o\big(b^{4}+\frac{1}{bT}\big).
\end{equation}
It can be shown (cf. Priestley, 1981, Chapter 7.5) that this MSE gets minimal for
\begin{equation} \label{OptimalKernel}
K(x)=K_{opt}(x)= 6 x(1-x),   \hspace {.4cm} 0\leq x \leq 1
\end{equation}
and
\begin{equation} \label{OptimalBandwidth}
  b=b_{opt}(u_0)= C(K_{opt})^{1/5}
  \bigg[\frac{\tau (u_0)}{\mu(u_0)^{2}}\bigg]^{1/5}T^{-1/5}
\end{equation}
where $C(K) = v_K / d_K^{2}$. In this case we have with $c(K) = v_K d_K^{1/2}$
\begin{equation} \label{OptimalMSE}
T^{4/5}\, \mean \, \big|\hat{c}_{T} (u_0,k) - c (u_0,k) \big|^{2}
   = \frac{5}{4} \, c(K_{opt})^{4/5} \,\mu(u_0)^{2/5} \, \tau(u_0)^{4/5}+ o(1).
\end{equation}

$\mu (u_0) =  \frac{\partial^{2}}{\partial^{2}u_0} \, c(u_0,k)$ measures the ``degree of
nonstationarity'' while $\tau (u_0)$ measures the variability of the estimate at time $u_0$. The
segment length $N_{opt}=b_{opt}T$ gets larger if $\mu (u_0)$ gets smaller, i.e. if the process is
closer to stationarity (in this case: if the $k$-th order covariance is more constant/more linear
in time). At the same time the mean squared error decreases. The results are similar to kernel
estimation in nonparametric regression. A yet unsolved problem  is how to adaptively determine the
bandwidth from the observed process.
\end{remark}

\underline{3. Segment selection and asymptotic mean squared-error for local Yule-Walker estimates}

For the local Yule-Walker estimates from (\ref{LocalYuleWalkerEst}) with the covariances ${\hat c}_{T}(u_0,k)$
as defined in (\ref{TaperedCovarianceEstimate}) Dahlhaus and Giraitis (1998) have proved (see also Example~\ref{LocalARestimatesDahGir})
\begin{equation*} \label{}
\mean \, \hat{\boldsymbol{\alpha}}_T (u_0) = \boldsymbol{\alpha} (u_0) - \frac {b^{2}} {2} \, d_K \, {\boldsymbol \mu} (u_0) + o (b^{2})
\end{equation*}
with
\begin{equation*} \label{}
{\boldsymbol \mu}(u_{0})  =  R(u_0)^{-1} \, \Big[ \Big(\frac {\partial^{2}} {\partial u^{2}} R(u)\Big) {\boldsymbol \alpha}(u_0) + \Big(\frac {\partial^{2}} {\partial u^{2}} r(u)\Big) \Big]_{u=u_0}
\end{equation*}
and
\begin{equation*} \label{}
\var\big( \hat{\boldsymbol{\alpha}}_T (u_0) \big)  = \frac {1} {bT}\, v_K  \,\sigma^{2}(u_0) \,R(u_0)^{-1} + o \Big(\frac {1} {bT}\Big).
\end{equation*}
Thus, we obtain for $\mean \big\|\hat{\boldsymbol{\alpha}}_T (u_0) - \boldsymbol{\alpha} (u_0) \big\|^{2}$ the same expression as in (\ref{AsymptoticMSE}) with $\tau(u_0)=\sigma^{2}(u_0) \, {\rm tr} \{R(u_0)^{-1}\}$ and $\mu (u_0)^{2}$ replaced by $\|{\boldsymbol \mu}(u_{0})\|^{2}$. With these changes the optimal bandwidth is given by (\ref{OptimalBandwidth}) and the optimal mean squared error by  (\ref{OptimalMSE}).

\begin{remark} [Implications for non-rescaled processes] \label{ImplicationsNonRescaled}
 Suppose
that we observe data from a (non-rescaled) tvAR(p)-process
\begin{equation} \label{ARpNonRescaled}
X_{t} + \sum_{j=1}^{p} \alpha_{tj}\; X_{t-j}=\sigma_t\;
\varepsilon_{t} ,\mbox{ } t\in\znum.
\end{equation}
In order to estimate $\boldsymbol{\alpha}_t$ at some time $t_0$ we may use the  segment
Yule-Walker estimator as given in (\ref{LocalYuleWalkerEst}).
The theoretically optimal segment length is given by (\ref{OptimalBandwidth}) as
\begin{equation} \label{NoptTVAR}
  N_{opt}(u_0)= C(K_{opt})^{1/5}
  \bigg[\frac{\tau (u_0)}{\|{\boldsymbol \mu}(u_{0})\|^{2}}\bigg]^{1/5}T^{4/5}
\end{equation}
which at first sight depends on $T$ and the rescaling.\\[2pt]
Suppose that we have parameter functions
$\tilde{a}_{j}(\cdot)$ and some $T>t_{0}$ with $\tilde{a}_{j}(\frac
{t_{0}} {T})=\alpha_j(t_{0})$ (i.e. the original function has been
rescaled to the unit interval) and we denote by $\tilde{R}$, $\tilde{r}$ and $\tilde{\boldsymbol{\alpha}}$ the corresponding parameters in the rescaled world (i.e. $\tilde{R}(u_0)=R(t_0)$ etc.). Then
\begin{equation*} \label{}
\tau(u_0)=\tilde{\sigma}^{2}(u_0) \, {\rm tr} \{\tilde{R}(u_0)^{-1}\} = \sigma^{2}(t_0) \, {\rm tr} \{R(t_0)^{-1}\}
\end{equation*}
and (with the second order difference as an approximation of the second derivative)
\begin{align*} \label{}
{\boldsymbol \mu}(u_{0}) & = \tilde{R}(u_0)^{-1} \, \bigg[ \Big(\frac {\partial^{2}} {\partial u^{2}} \tilde{R}(u)\Big) \tilde{\boldsymbol \alpha}(u_0) + \Big(\frac {\partial^{2}} {\partial u^{2}} \tilde{r}(u)\Big) \bigg]_{u=u_0}\\[6pt]
& \approx R(t_0)^{-1} \, \bigg[ \frac {R(t_0) - 2 R(t_0 -1) + R(t_0 -2)} {1/T^{2}} \, \boldsymbol a(t_0) +  \frac {r(t_0) - 2 r(t_0 -1) + r(t_0 -2)} {1/T^{2}} \bigg]
\end{align*}
Plugging this into (\ref{NoptTVAR}) reveals that $T$ drops out completely and the optimal segment length can completely be determined in terms of the original non-rescaled process. This is a nice example on how the asymptotic considerations in the rescaled world can be transferred with benefit to the original non-rescaled world.
\hfill $\Box$
\end{remark}

These considerations justify the asymptotic
approach of this paper: While it is not possible to set down a
meaningful asymptotic theory for the non-rescaled model (\ref{AR1-nonrescaled}) an approach using the rescaled
model (\ref{ARpRescaled}) leads to meaningful results also for the  model (\ref{AR1-nonrescaled}).
Another example for this relevance is the construction of confidence intervals for the local Yule-Walker estimates
from the central limit theorem in Dahlhaus and Giraitis (1998), Theorem 3.2.
\medskip

\underline{4. Parametric Whittle-type estimates -- a first approach}

We now assume that the $p+1$-dimensional parameter curve ${\boldsymbol \theta} (\cdot)=
\big(\alpha_1 (\cdot),\ldots,\alpha_p (\cdot),\sigma^{2}(\cdot)\big)'$ is parameterized by a finite
dimensional parameter $\eta \in {\mathbf R}^q$ that is ${\boldsymbol \theta}(\cdot)={\boldsymbol
\theta}_\eta (\cdot)$. An example studied below is where the AR-coefficients are modeled by
polynomials. Another example is where the AR-coefficients are modeled by a parametric transition curve as in Section~\ref{DeepExample}.6(iv). In particular when the length of the time series is short this may be a proper choice.
We now show how the stationary Whittle likelihood can be generalized to the locally stationary case
(another generalization is given in (\ref{GeneralizedWhittleLikelihood})).

If we were looking for a \underline{nonparametric} estimate for the parameter curve ${\boldsymbol \theta} (\cdot)$ we
could apply the stationary Whittle estimate on a segment leading to
\begin{equation} \label{NonparametricWhittleEstimates}
\hat{\boldsymbol{\theta}}_{T}^{W}(u_0) := \argmin_{\theta \in \Theta}
\mathcal{L}_{T}^{W}(u_0, \boldsymbol{\theta})
\end{equation}
with the Whittle likelihood
\begin{equation} \label{NonparametricWhittleLikelihood}
\mathcal{L}_{T}^{W}(u_0, \boldsymbol{\theta}) := \frac{1}{4\pi} \int^{\pi}_{-\pi} \biggl\{
\log 4\pi^{2} f_{\boldsymbol{\theta}} (\lambda) + \frac{I_{T}(u_0,\lambda)}{f_{\boldsymbol{\theta}}
(\lambda)} \biggr\} \, d \lambda
\end{equation}
with the tapered periodogram on a segment about $u_0$, that is
\begin{equation} \label{PeriodogramOnSegment}
I_{T} (u_0,\lambda) := \frac{1}{2\pi H_{N}} \bigg| \sum_{s = 1}^{N} h \!
\left( \frac{s}{N} \right) X_{[u_0T] - N/2 + s,T} \, \exp \big(\!-i \lambda s\big)
\bigg|^{2}.
\end{equation}
Here $h(\cdot)$ is a data taper as in (\ref{TaperedCovarianceEstimate}). For $h(x) = \chi_{(0,1]}
(x) $ we obtain the non-tapered periodogram. The properties of this nonparametric estimate are
discussed later - in particular in Example~\ref{LocalWhittleEstimatesBiasVariance} and at the end of Example~\ref{ExampleLocalWhittle}. In case of a tvAR(p)-process $\hat{\boldsymbol{\theta}}_{T}(u_0)$ is exactly the local Yule-Walker estimate defined in (\ref{LocalYuleWalkerEst}) with the covariance-estimate given in (\ref{TaperedCovarianceEstimate}).

Suppose now that we want to fit globally the parametric model ${\boldsymbol \theta}(\cdot)={\boldsymbol \theta}_\eta (\cdot)$ to the data, that is we have the time varying spectrum $f_\eta (u,\lambda) := f_{\boldsymbol{\theta}_\eta(u)}(\lambda)$. Since $\mathcal{L}_{T}^{W}(u, \boldsymbol{\theta})$ is an approximation of the
Gaussian log-likelihood on the segment $\{[uT] - N/2 + 1,\ldots,[uT] + N/2\}$ a reasonable approach is to use
\begin{equation} \label{BlockWhittleEstimate}
\hat{\eta}_{T}^{BW} := \argmin_{\eta \in \Theta_\eta}
\mathcal{L}_{T}^{BW}(\eta)
\end{equation}
with the \underline{block Whittle likelihood}
\begin{equation} \label{BlockWhittleLikelihood}
{\mathcal L}_T^{BW} (\eta) :=\frac{1}{4 {\pi}}\frac{1}{M} \sum_{j=1}^{M}
\int_{- \pi}^{\pi} \left\{ \log 4\pi^{2} f_{\eta}(u_j,\lambda)
+\frac{I_T(u_j,\lambda)}{f_{\eta}(u_j,\lambda)} \right\}
d\lambda.
\end{equation}
Here $u_j:=t_j/T$ with $t_j:=S(j-1)+N/2~~(j=1,\ldots,M)$ i.e. we calculate the likelihood on overlapping segments which we shift each time by $S$. Furthermore $T=S(M-1)~+~N$. A better justification of the form of the likelihood is provided by the asymptotic Kullback-Leibler information divergence derived in Theorem~\ref{LimitLikelihood}.

As discussed above the reason for using data-tapers is twofold: they reduce the bias due to nonstationarity on the segment and they reduce the leakage (already known from the stationary case). It is remarkable that the taper in this case does not lead to an increase of the asymptotic variance if the segments are overlapping (cf. Dahlhaus (1997), Theorem~3.3).

The properties of the above estimate are discussed in Dahlhaus (1997) including consistency, asymptotic normality, model selection and the behavior if the model is misspecified. The estimate is asymptotically efficient if $S/N \rightarrow 0$.

\begin{figure*}
\centering \hspace*{0.3cm}
\includegraphics[width=330pt,keepaspectratio]{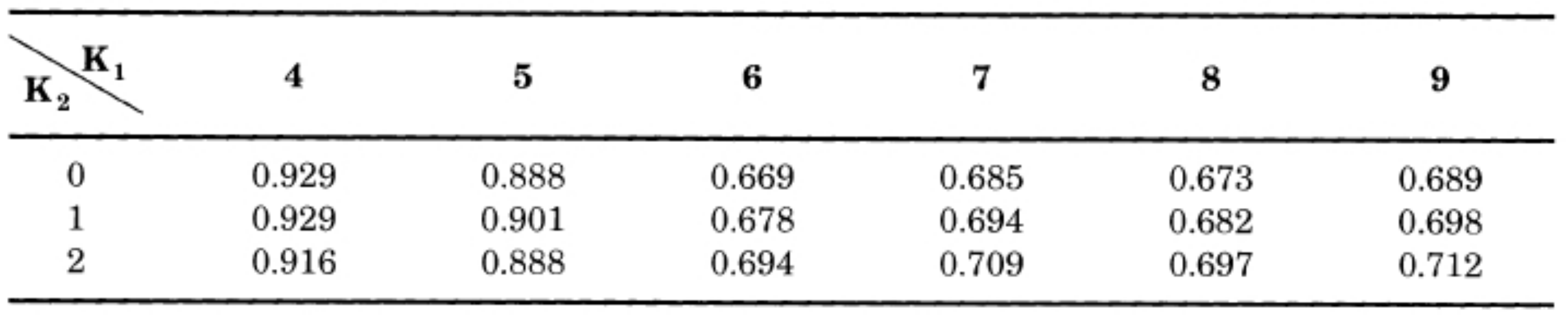}\\ 
Table 1: {\footnotesize Values for AIC for $p=2$ and different polynomial orders}
\label{}
\end{figure*}

As an example we now fit a tvAR(p)-model to the data from Figure~\ref{fig1} and estimate the parameters by minimizing ${\mathcal L}_T^{BW}(\eta)$. The AR-coefficients are modeled
as polynomials with different orders. Thus, we fit the
model
\begin{equation*} \label{}
\alpha_j(u)=\sum_{k=0}^{K_j} b_{jk} \, u^k \quad (j=1,\ldots,p) \quad \mbox{and} \quad \sigma (u) \equiv c
\end{equation*}
to the data. The model orders $p,K_1,\ldots,K_p$ are chosen
by minimizing the AIC-criterion

\begin{equation*} \label{}
AIC(p,K_1,\ldots,K_p)=\log\hat{\sigma}^2(p,K_1,\ldots,K_p)+2\,\big(p+1+\sum^p_{j=1}
  K_j\big)\,/\,T.
\end{equation*}

\begin{figure}
\centering
\includegraphics[width=350pt,keepaspectratio]{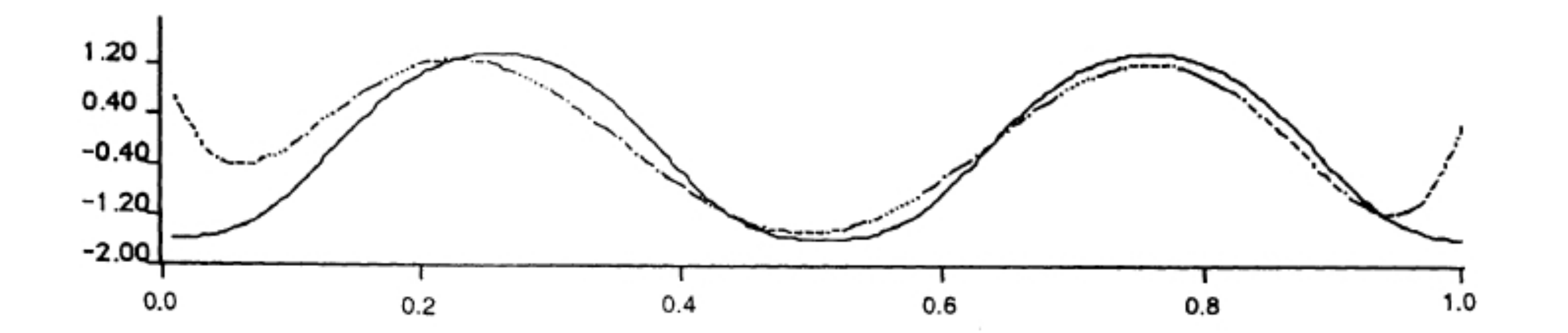}
\caption{\footnotesize True and estimated parameter curve $\alpha_{1}(\cdot)$}
\label{fig3}
\end{figure}

\begin{figure}
\centering
\includegraphics[width=305pt,keepaspectratio]{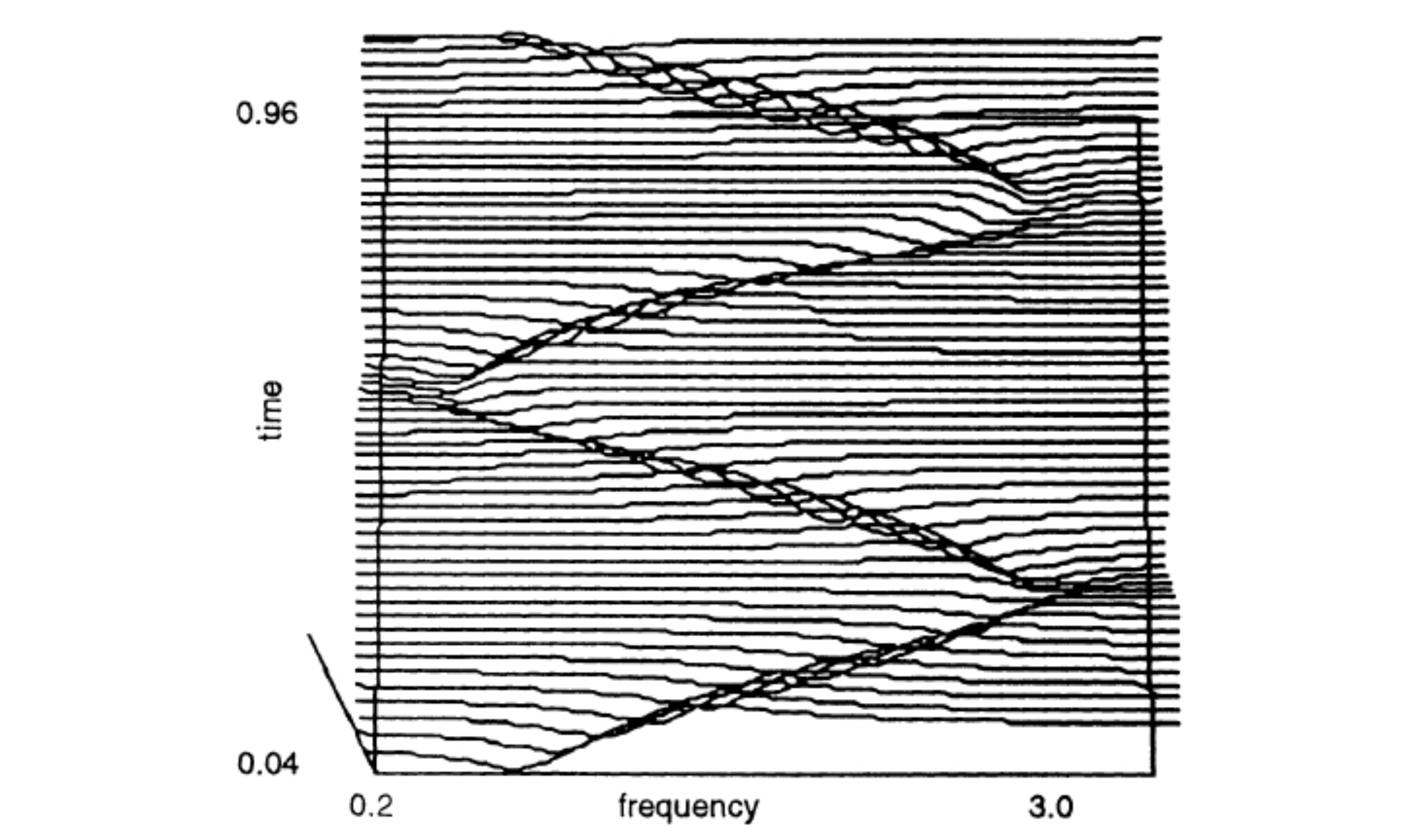} 
\caption{\footnotesize Difference of estimated and true spectrum}
\label{fig4}
\end{figure}

Table 1 shows these values for $p=2$ and different $K_1$ and $K_2$. The values for other
$p$ turned out to be larger. Thus, a model with $p=2, K_1=6, K_2=0$
is fitted. The function $\alpha_1(u)$ and its estimate are plotted in Figure~\ref{fig3}. For $\hat{a}_2(u)$ we obtain $0.71$
(a constant is fitted because of $K_2=0$) while the true
$\alpha_2(u)$ is 0.81. Furthermore, $\hat{\sigma}^2=1.71$ while
$\sigma^2=1.0$. The corresponding (parametric) estimate of the spectrum is the right picture of Figure~\ref{fig2} and the difference to the true spectrum is plotted in Figure~\ref{fig4}.

Given the small sample size the quality of the fit is remarkable. Two negative
effects can be observed. First, the fit of $\alpha_1(u)$ becomes rather bad outside $u_1=0.063$ and
$u_M=0.938$. This is not surprising, due to the behavior of a polynomial
and the fact that the use of ${\mathcal L}_T^{BW}(\eta)$ as a distance only punishes bad
fits inside the interval $[u_1,u_M]$. This end effect
improves if one chooses $K_1=8$ instead of $K_1=6$. A better way
seems to modify $\mathcal {L}_T^{BW}(\eta)$ and to include periodograms of shorter lengths
at the edges. The second effect is
that the peak in the spectrum is underestimated. This bias is in part due to the non-stationarity of the process on intervals $(u_j-N/(2T),u_j+N/(2T))$ where $I_T(u_j,\lambda)$ is calculated.

We mention that the above estimates can be written in closed form and calculated without an optimization routine. More generally this holds for tvAR(p)-models if $\sigma^2$ is constant and $\alpha_j(u)=\sum_{k=1}^K b_{jk}\, f_k(u)$ with some functions $f_1(u),\ldots,f_K(u)$ (in the above case $f_k(u)=u^{k-1}$). For details see Dahlhaus (1997), Section 4.

A closer look at the above estimate reveals that it is somehow the outcome of a two step procedure where in the first step the periodogram is calculated on  segments (which implicitly includes some smoothing with bandwidth $b= N/T$) and afterwards the AR(p)-process with the above polynomials is fitted to the outcome (instead of a direct fit of the  AR(p)-model and the polynomials to the data). We now make this more precise.

With the above form of the spectrum $f_{\eta}(u,\lambda)$ (cf.(\ref{tvARSpectralDensity})) and Kolmogorov's
formula, (cf. Brockwell and Davis, 1991, Theorem 5.8.1) we obtain with $\hat{R}_T(u_j)$ and $\hat{r}_T (u_j)$ as defined
in (\ref{LocalYuleWalkerEst}) after some straightforward calculations
\begin{align*}
{\mathcal L}_T^{BW}(\eta)&=\frac{1}{2}{\frac{1}{M}\sum_{j=1}^M} \Big[
\log 4\pi^{2} \sigma_{\eta}^2 (u_j) + \frac{1}{\sigma_{\eta}^2(u_j)} \Big( \hat{c}_T (u_j,0) -
\hat{r}_T (u_j)^{'}
  \hat{R}_T(u_j)^{-1} \hat{r}_T (u_j) \Big)\Big]\\
&+ \frac{1}{2}{\frac{1}{M}\sum_{j=1}^M} \frac{1}{\sigma_{\eta}^2(u_j)} \Big[ \big( \hat{R}_T(u_j) \,
  \boldsymbol{\alpha}_{\eta}(u_j)
  + \hat{r}_T (u_j) \big)' \hat{R}_T(u_j)^{-1}
  \big( \hat{R}_T(u_j) \, \boldsymbol{\alpha}_{\eta}(u_j)
  + \hat{r}_T (u_j) \big) \Big].
\end{align*}

We now plug in the Yule-Walker estimate $\hat{\boldsymbol{\alpha}}_T (u) = - \hat{R}_T(u)^{-1} \,
\hat{r}_T (u)$ with asymptotic variance proportional to $\sigma^{2}(u) \,R(u)^{-1}$ and $\hat{\sigma}^2_{T}(u) =\hat{c}_T(u, 0)- \hat{r}_T (u)' \,\hat{R}_T(u)^{-1} \,\hat{r}_T (u)$ with asymptotic variance $2\, \sigma^{2}(u)$. Since
$\log x=(x-1)-\frac{1}{2}(x-1)^2+o((x-1)^2)$ we  obtain
\begin{align*}
{\mathcal L}_T^{BW}(\eta) &= \frac{1}{2}{\frac{1}{M}\sum_{j=1}^M} \frac {1} {2\, \sigma_{\eta}^4(u_j)} \Big[
\sigma_{\eta}^2(u_j) - \hat{\sigma}^2_{T}(u_j)\Big]^{2}\\
&+ \frac{1}{2}{\frac{1}{M}\sum_{j=1}^M} \Big[ \Big(
  \boldsymbol{\alpha}_{\eta}(u_j)
  - \hat{\boldsymbol{\alpha}}_T (u_j) \Big)' \sigma_{\eta}^2(u_j)^{-1}\,\hat{R}_T(u_j)
  \Big(  \boldsymbol{\alpha}_{\eta}(u_j)
 - \hat{\boldsymbol{\alpha}}_T (u_j) \Big) \Big]\\
& + \frac{1}{2}{\frac{1}{M}\sum_{j=1}^M} \log 4\pi^{2} \hat{\sigma}^2_{T}(u_j) + \frac {1} {2} + o\Big(\big(\frac {\sigma_{\eta}^2(u_j)-\hat{\sigma}^2_{T}(u)} {\sigma_{\eta}^2(u_j)}\big)^{2}\Big).
\end{align*}
If the model is correctly specified then we have for $\eta$ close to the minimum: $\sigma_{\eta}^2(u_j)^{-1}\,\hat{R}_T(u_j) \approx \sigma^{2}(u_j)^{-1} \,R(u_j)$ and $2\, \sigma_{\eta}^4(u_j) \approx 2\, \sigma^{2}(u_j)$ which means that $\hat {\eta}_T$ is approximately obtained by a
weighted least squares fit
of $\boldsymbol{\alpha}_{\eta}(u)$ and $\sigma^2_{\eta}(u)$ to the Yule-Walker
estimates on the segments. The method works in this case since the (parametric!) model fitted in the second step is somehow `smoother' than the first smoothing implicitly induced by using the periodogram on a segment. However, we would clearly run into problems if the fitted polynomials were of high order or if even $K_{j}=K_{j}(T) \rightarrow \infty$ as $T \rightarrow \infty$.

A good alternative seems to use the quasi-likelihood ${\mathcal L}_{T}^{GW} (\eta)$ from
(\ref{GeneralizedWhittleLikelihood}) or (in particular for AR(p)-models) the conditional likelihood estimate
from (\ref{FullConditionalLikelihood}) with $\ell_{t,T} (\cdot)$ as in (\ref{ConditionalLikelihoodARp}) for which the estimator can explicitly be calculated if $\sigma(\cdot) \equiv c$. For $\sigma_{0}(\cdot) \neq c$ iterative or approximative solutions are needed. The properties of this estimator have not been investigated yet. In any case the benefit of the likelihood ${\mathcal L}_T^{BW}(\eta)$ and even more of the improved likelihood ${\mathcal L}_{T}^{GW} (\eta)$ are their generality because they can be applied to arbitrary parametric models which can be identified from the second order spectrum.

Furthermore, algorithmic issues, such as in-order algorithms (e.g. generalizations of the Levinson-Durbin algorithm) need to be developed.

\medskip

\underline{5. Inference for nonparametric tvAR-models -- an overview}

In the last section we studied parametric estimates for tvAR(p)-models. This is an
important option if the length of the time series is short or if we have specific parametric models
in mind. In general however one would prefer nonparametric models. For
nonparametric statistics a large variety of different estimates are available (local polynomial
fits, estimation under shape restrictions, wavelet methods etc) and it turns out that it is not too
difficult to apply such methods to tvAR(p)-models and moreover also to other possibly nonlinear models (while the derivation of the corresponding theory may be very challenging).
 A key role is played by the conditional likelihood at time $t$ which in the tvAR(p)-case is
\begin{align}
\ell_{t,T} ({\boldsymbol \theta})
&:= - \log f_{{\boldsymbol \theta}} \big(X_{t,T} \big| X_{t-1,T},\ldots,X_{1,T}\big)\label{ConditionalLikelihoodAtTimet}\\
& \;= \frac{1}{2} \,
\textrm{log} \big( 2\pi\,\sigma^2 \big) + \frac{1}{2\,\sigma^2}
\;\Big( X_{t,T} + \sum_{j=1}^{p} \alpha_{j}\; X_{t-j,T} \Big)^{2} \label{ConditionalLikelihoodARp}
\end{align}
where ${\boldsymbol \theta} = \big(\alpha_1,\ldots,\alpha_p,\sigma^{2}\big)'$ and its approximation $\ell^{\,*}_{t,T} ({\boldsymbol \theta})$ defined in
(\ref{LocalWhittleLikelihoodAtTimet}). As a simple example consider the estimation of the curve $\alpha_1 (\cdot)$ of
a tvAR(1)-process by a local linear fit given by $\hat{\alpha}_1 (\cdot)=\hat{c}_0$ where
\begin{equation} \label{61}
(\hat{c}_0,\hat{c}_1) = \argmin_{c_0,c_1}
 \frac{1}{bT} \sum^{T}_{t=1} K \! \Big(
\frac{u_0 -t/T}{b} \Big) \, \bigg( X_{t,T} + \Big[c_0 + c_{1}
\Big(\frac{t}{T}-u_0 \Big)\Big] X_{t-1,T} \bigg)^{2}
\end{equation}
or more generally (with vectors $\boldsymbol{c}_0$ and $\boldsymbol{c}_1$) given by $\hat{\boldsymbol{\theta}}(u_0) =\hat{\boldsymbol{c}}_0$ with
\begin {equation} \label{LocalLinearEstimate}
(\hat{\boldsymbol{c}}_0,\hat{\boldsymbol{c}}_1) = \argmin_{\boldsymbol{c}_0,\boldsymbol{c}_1} \frac{1}{bT} \sum^{T}_{t=1} K \! \Big( \frac{u_0
-t/T}{b} \Big) \;\ell_{t,T} \Big( \boldsymbol{c}_0 + \boldsymbol{c}_{1} \big(\frac{t}{T}-u_0 \big) \Big).
\end{equation}

Besides this local linear estimate many other estimates can be constructed based on the conditional
likelihood $\ell_{t,T} ({\boldsymbol \theta})$ from above:
\begin{enumerate}
\item A \underline{kernel estimate} defined by
\begin{equation} \label{LocalKernelEstimate}
\hat{\boldsymbol{\theta}} (u_0) = \argmin_{\boldsymbol{\theta}} \frac{1}{bT} \sum^{T}_{t=1} K
\! \Big( \frac{u_0 -t/T}{b} \Big) \, \ell_{t,T}\! \big(\boldsymbol{\theta} \big).
\end{equation}
This estimate is studied in Section~\ref{NonlinearProcAndDerivative}. We are convinced that
it is equivalent to the local Yule-Walker estimate from (\ref{LocalYuleWalkerEst}) with $K(x) = h(x)^{2}$,
$b = N/T$ and that all results from \underline{3.} are exactly the same for this estimate.
\item A \underline{local polynomial fit} defined by $\hat{\boldsymbol{\theta}}(u_0) =\hat{\boldsymbol{c}}_0$ with
\begin {equation} \label{LocalPolynomialFit}
(\hat{\boldsymbol{c}}_0,\dots,\hat{\boldsymbol{c}}_d)' =
\argmin_{\boldsymbol{c}_0,\dots,\boldsymbol{c}_d} \frac{1}{bT} \sum^{T}_{t=1} K \! \Big( \frac{u_0
-t/T}{b} \Big) \;\ell_{t,T}  \bigg( \sum^{d}_{j=0} \boldsymbol{c}_{j} (\frac{t}{T} - u_0)^{j} \bigg).
\end{equation}
Local polynomial fits for tvAR(p)-models have been investigated by Kim (2001)
and Jentsch (2006).
\item An \underline{orthogonal series estimate} (e.g. a \underline{wavelet estimate}) defined by
\begin {equation} \label{OrthogonalSeriesEstimate}
\bar{\boldsymbol{\beta}} = \argmin_{\boldsymbol{\beta}} \frac{1}{T} \sum^{T}_{t=1} \ell_{t,T} \! \biggl( \sum^{J(T)}_{j=1}
\boldsymbol{\beta}_{j} \psi_{j} \big(\frac{t}{T}\big) \biggr)
\end {equation}
together with some shrinkage of $\bar{\boldsymbol{\beta}}$ to obtain $\hat{\boldsymbol{\beta}}$ and $\hat{\boldsymbol{\theta}} (u_0) = \sum^{J(T)}_{j=1}
\hat{\boldsymbol{\beta}}_{j} \psi_{j} (u_0)$. Usually $J(T) \rightarrow \infty$ as $T \rightarrow \infty$. Such an estimate has been investigated for a truncated wavelet expansion for tvAR(p)-models in Dahlhaus, Neumann and von
Sachs (1999).
\item A \underline{nonparametric maximum likelihood estimate} defined by
\begin {equation} \label{GlobalNonparametricEstimate}
\hat{\boldsymbol{\theta}}(\cdot) = \argmin_{\boldsymbol{\theta}(\cdot) \in \Theta}
\frac{1}{T} \sum^{T}_{t=1} \ell_{t,T}\! \Big(\boldsymbol{\theta}\big(\frac {t} {T}\big) \Big)
\end {equation}
where $\Theta$ is an adequate function space, for example a space of curves under shape
restrictions such as monotonicity constraints. In Dahlhaus and Polonik (2006) the estimation of a monotonic variance
function in a tvAR-model is studied, including explicit
algorithms involving isotonic regression.
\item A \underline{parametric fit} for the curves $\boldsymbol{\theta}(\cdot)=\boldsymbol{\theta}_\eta
    (\cdot)$ with $\eta \in {\mathbf R}^q$ defined by
\begin {equation} \label{FullConditionalLikelihood}
\hat{\eta} = \argmin_{\eta} \frac{1}{T} \sum^{T}_{t=1} \ell_{t,T}\!
 \Big(\boldsymbol{\theta}_\eta \big(\frac {t} {T}\big) \Big)
\end {equation}
The resulting estimate has not been investigated yet. It is presumably very close to the exact MLE studied in Theorem~\ref{AsymptoticPropertiesGWhittle}.
\end{enumerate}
\begin{remark} \label{}
(i) In the tvAR(p)-case the situation simplifies a lot if $\sigma^{2}(\cdot) \equiv c$. In that case the estimates for $\boldsymbol{\alpha} (\cdot)$ and $\sigma^{2}$ ``split'' and $\ell_{t,T} ({\boldsymbol \theta})$ can in all cases be replaced by $\big( X_{t,T} + \sum_{j=1}^{p} \alpha_{j}\; X_{t-j,T} \big)^{2}$ leading to least squares type estimates.\\[4pt]
(ii) All estimates from above can be transferred to other models by using the
conditional likelihood (\ref{ConditionalLikelihoodAtTimet}) for the specific model.
The kernel estimate will be investigated in Section~\ref{NonlinearProcAndDerivative}.\\[4pt]
(iii) As mentioned above an alternative choice is to replace $\ell_{t,T} ({\boldsymbol \theta})$ by the local generalized
Whittle likelihood $\ell^{\,*}_{t,T} ({\boldsymbol \theta})$ from
(\ref{LocalWhittleLikelihoodAtTimet}). With that likelihood several estimates from above have been investigated -- see the detailed discussion at the end of Section~\ref{SectionLikelihood}. In that case the $d$-dimensional parameter curve
$\boldsymbol{\theta}(\cdot)= \big(\theta_1(\cdot),\dots,\theta_d(\cdot)\big)'$ must be uniquely
identifiable from the time varying spectrum $f(u,\lambda) = f_{\boldsymbol{\theta}(u)}(\lambda)$.
\hfill $\Box$
\end{remark}

\medskip

\underline{6. Shape- and transition curves}

There exist several alternative models for tvAR-processes -- in particular models where specific characteristics of the time series are modeled by a curve. Below we give 4 examples where we restrict ourselves to tvAR(2)-models. Suppose we have a stationary  AR(2)-model with complex roots $\frac {1} {r} \exp( i \phi)$ and $\frac {1} {r} \exp( - i \phi)$, that is with parameters $a_1 = -2 r \cos(\phi)$, $a_2 = r^{2}$, and variance $\sigma^{2}$. The corresponding process shows a quasi-periodic behavior with period of length $\frac {2\pi} {\phi}$, that is with  frequency $\phi$. The more $r$ gets closer to $1$ the more the shape of the process gets closer to a sine-wave. The amplitude is proportional to $\sigma$ (if $\sigma$ (say in (\ref{ARpRescaled})) is replaced by $c \cdot \sigma$, then $X_t$ is replaced by $c \cdot X_t$).

In the specific tvAR(2)-case we can now consider the following shape- and transition-models for quasi-periodic processes:

\underline{(i) Model with a time varying amplitude curve:}
\begin{equation*} \label{}
a_1(\cdot), \, a_2 (\cdot) \; \mbox{constant}; \quad \sigma(\cdot) \; \mbox{time varying}.
\end{equation*}
Chandler and Polonik (2006) use this model with a unimodal $\sigma(\cdot)$ and a nonparametric maximum likelihood estimate for the discrimination of earthquakes and explosions. The properties of the estimator have been investigated in Dahlhaus and Polonik (2006).

\underline{(ii) Model with a time varying frequency curve:}
\begin{equation*} \label{}
a_1 (\cdot) = -2 \,r \cos\big(\phi (\cdot)\big), \; a_2 (\cdot) = r^{2}   \mbox{ with } r \; \mbox{constant and } \phi(\cdot) \mbox{ time varying}, \; \sigma (\cdot) \mbox{ constant}.
\end{equation*}
The model in Figure 1 is of this form with $r=0.9$ and $ \phi (u) = 1.5-\cos4{\pi} u$.

\underline{(iii) Model with a time varying period-distinctiveness:}
\begin{equation*} \label{}
a_1 (\cdot) = -2 \,r (\cdot) \cos(\phi ), \; a_2 (\cdot) = r(\cdot)^{2}   \mbox{ with } r (\cdot) \; \mbox{time varying and } \phi \mbox{ constant}, \; \sigma (\cdot) \mbox{ constant}.
\end{equation*}

\underline{(iv) Transition models:} Amado and Ter\"{a}svirta (2011)
have recently used the logistic transition function to model parameter
transitions in GARCH-models. The simplest transition function is
\begin{equation*} \label{}
G\Big(\frac {t} {T}; \gamma,c\Big) := \Big[ 1 + \exp\Big\{-\gamma\,\Big(\frac {t} {n}-c\Big)\Big\}\Big]^{-1}.
\end{equation*}
Since $G(0; \gamma,c) \approx 0$ and $G(1; \gamma,c) \approx 1$ the model
\begin{equation*} \label{}
a_1 (u) = a_{1}^{\scriptsize \mbox{start}} + G(u; \gamma,c) \, \big(a_{1}^{\mbox{\scriptsize end}}-a_{1}^{\scriptsize \mbox{start}}\big), \; a_{2} (u) = a_{2}^{\scriptsize \mbox{start}} + G(u; \gamma,c) \, \big(a_{2}^{\mbox{\scriptsize end}}-a_{2}^{\scriptsize \mbox{start}}\big)
\end{equation*}
is a parametric model for a smooth transition from the AR-model with parameters $(a_{1}^{\scriptsize \mbox{start}},a_{2}^{\scriptsize \mbox{start}})$ at $u=0$
to the model with parameters $(a_{1}^{\scriptsize \mbox{end}},a_{2}^{\scriptsize \mbox{end}}\big)$ at $u=1$. Here $c$ and $\gamma$ are the location and the `smoothness' of transition respectively. More general transition models (in particular with more states) may be found in Amado and Ter\"{a}svirta (2011). $G(\cdot; \gamma,c)$ may also be replaced by a (nonparametric) function $G(\cdot)$ with $G(0)=0$ and $G(1)=0$.

It is obvious that all methods from subsection 5 can be applied in cases (i)-(iv) to estimate the constant parameters and the shape- and transition-curves. We mention that the theoretical results for local Whittle estimates of Dahlhaus and Giraitis (1998) apply to these models (cf. Example~\ref{LocalWhittleEstimatesBiasVariance}), the uniform convergence result for the local generalized Whittle estimate in Theorem~\ref{uniform-convergence-Whittle}, the asymptotic results of Dahlhaus and Neumann (2001) where the parameter curves are estimated by a nonlinear wavelet method, the results of Dahlhaus and Polonik (2006) on nonparametric maximum likelihood estimates under shape constraints, and the results for parametric models in Theorem~\ref{AsymptoticPropertiesGWhittle} on the MLE and the generalized Whittle estimator, and in Dahlhaus (1997) on the block Whittle estimator.

\section{Local likelihoods, derivative processes and nonlinear models with time varying parameters} \label{NonlinearProcAndDerivative}

In this section we present a more general framework for time series with time varying finite dimensional parameters ${\boldsymbol \theta} (\cdot)$ and show how nonparametric inference can be done and theoretically handled. Typically such models result from the generalization of classical parametric models to the time varying case. If we restrict ourselves to linear processes or even more to Gaussian processes then a much more general theory is possible which is developed in the subsequent sections. Large parts of the present section are based on the ideas presented in Dahlhaus and Subba Rao (2006) where time varying ARCH-models have been investigated.

The key idea is to use at each time point $u_0 \in (0,1)$ the stationary approximation $\tilde{X}_{t}(u_0)$ to the original process $X_{t,T}$ and to calculate the bias resulting from the use of this approximation. This will end in Taylor-type expansions of $X_{t,T}$ in terms of so-called derivative processes. These expansions play a major role in the theoretical derivations.

Suppose for example that we estimate the multivariate parameter curve ${\boldsymbol \theta} (\cdot)$ by minimizing the (negative) local conditional log-likelihood, that is
\begin{equation*} \label{}
\hat{{\boldsymbol \theta}}_{T}^{C}(u_0) := \argmin_{{\boldsymbol \theta} \in \Theta}
\mathcal{L}_{T}^{C}(u_0, {\boldsymbol \theta})
\end{equation*}
with
\begin{equation} \label{localConditionalLikelihood}
\mathcal{L}_{T}^{C}(u_0, {\boldsymbol \theta}) :=
\frac{1}{T} \sum^{T}_{t=1}  \;\frac{1}{\,b}\, K \Big(\frac{u_0-t/T}{b}\Big) \; \ell_{t,T} ({\boldsymbol \theta})
\end{equation}
and
\begin{equation*} \label{}
\ell_{t,T} ({\boldsymbol \theta}) := - \log f_{{\boldsymbol \theta}} \big(X_{t,T} \big| X_{t-1,T},\ldots,X_{1,T}\big)
\end{equation*}

where $K$ is symmetric, has compact support $[-\frac {1} {2},\frac {1} {2}]$ and fulfills
$\int_{-1/2}^{1/2} K(x) \, dx = 1$. We assume that $b=b_{T} \rightarrow 0$ and $b\,T \rightarrow
\infty$ as $T \rightarrow \infty$. Two examples for this likelihood are given below.

We approximate $\mathcal{L}_{T}^{C}(u_0, {\boldsymbol \theta})$ with $\tilde{\mathcal{L}}_{T}^{C}(u_0, {\boldsymbol \theta})$ which is the same function but with $\ell_{t,T} ({\boldsymbol \theta})$ replaced by
\begin{equation*} \label{}
\tilde{\ell}_t (u_0,{\boldsymbol \theta}) := - \log f_{{\boldsymbol \theta}} \big( \tilde{X}_{t}(u_0) \big| \tilde{X}_{t-1}(u_0),\ldots,\tilde{X}_{1}(u_0)\big),
\end{equation*}
which means that $X_{t,T}$ is replaced by its stationary approximation $\tilde{X}_{t}(u_0)$. Usually this is the local conditional likelihood for the process $\tilde{X}_{t}(u_0)$.

\begin{example} \label{ExampleNonlinearModels}
(i) Consider the tvAR(p) process defined in (\ref{ARpRescaled}) together with its stationary
approximation at time $u_0$ given by (\ref{ARpApproximation}). Under suitable regularity conditions
it can be shown that $X_{t,T} = \tilde{X}_{t}(u_0) + O_{p}\big(\big|\frac{t}{T}-u_0\big| + \frac {1} {T}\big)$ \big(cf.(\ref{ExpansionDerivativeProcesses})\big).
In case where the $\varepsilon_t$ are Gaussian the conditional likelihood at time $t$  is given by
\begin{equation} \label{LocalLikelihoodARpCase}
\ell_{t,T} ({\boldsymbol \theta}) = \frac{1}{2} \,
\textrm{log} \big( 2\pi\,\sigma^2 \big) + \frac{1}{2\,\sigma^2}
\;\Big( X_{t,T} + \sum_{j=1}^{p} \alpha_{j}\; X_{t-j,T} \Big)^{2}
\end{equation}
where ${\boldsymbol \theta}= (\alpha_1,\ldots,\alpha_p,\sigma^{2})'$. It is easy to show that the resulting estimate is the same as in (\ref{LocalYuleWalkerEst})
but with $\hat{r}_T (u_0):= \big(\tilde{\tilde{c}}_{T}(u_0,0,1),...,\tilde{\tilde{c}}_{T}(u_0,0,p)\big)'$ and $\hat{R}_T (u_0):=\{\tilde{\tilde{c}}_{T}(u_0,i,j)\}_{i,j=1,...,p}$ with the local covariance estimator $\tilde{\tilde{c}}_{T}(u,i,j)$ as defined in (\ref{KernelCovarianceEstimate2}).

(ii) A tvARCH(p) model where $\{X_{t,T}\}$ is assumed
to satisfy the representation
\begin{eqnarray}
\label{eq:tvarch-p}
 && X_{t,T} = \sigma_{t,T}Z_{t} \nonumber\\
\textrm{ where } &&  \sigma_{t,T}^{2}  = \alpha_{0}(\frac{t}{T}) +
\sum_{j=1}^{p}\alpha_{j}(\frac{t}{T})X_{t-j,N}^{2}
\quad \textrm{for } t=1,\ldots,N
\end{eqnarray}
with $Z_{t}$ being independent, identically distributed random
variables with $\mean Z_{t} = 0$, $\mean Z_{t}^{2}=1$.

The corresponding
stationary approximation $\tilde{X}_{t}(u_0)$ at time $u_{0}$ is given by
\begin{eqnarray}
\label{eq:tvarch-rescale}
 && \tilX_{t}(u_{0}) = \sigma_{t}(u_{0})\, Z_{t} \nonumber\\
\textrm{ where } &&  \sigma_{t}(u_{0})^{2}  = \alpha_{0}(u_{0}) +
\sum_{j=1}^{p}\alpha_{j}(u_{0})\, \tilX_{t-j}(u_{0})^{2}\quad \textrm{for } t \in \znum.
\end{eqnarray}
It is shown in Dahlhaus and Subba Rao (2006) that $\{X_{t,T}^{2}\}$ as defined above has an almost surely
well-defined unique solution in the set of all causal solutions and $X_{t,T}^{2} = \tilX_{t}(u_{0})^{2} +
O_{p}(|\frac{t}{T} - u_{0}| + \frac{1}{N})$. In case where the  $Z_t$ are Gaussian the conditional likelihood is given by
\begin{equation} \label{LocalLikelihoodTVArch}
\ell_{t,T} ({\boldsymbol \theta}) = \frac{1}{2} \log
w_{t,T}({\boldsymbol {\boldsymbol \theta}}) +
\frac{X_{t,T}^{2}}{2 \,w_{t,T}({\boldsymbol {\boldsymbol \theta}})}
\textrm{ with }  w_{t,T}({\boldsymbol {\boldsymbol \theta}}) = \alpha_{0} +
\sum_{j=1}^{p} \alpha_{j}X_{t-j,T}^{2}
\end{equation}
where ${\boldsymbol \theta}= (\alpha_0,\ldots,\alpha_p)'$. Dahlhaus and Subba Rao (2006) prove consistency of the resulting estimate also in case where the true process is not Gaussian. As an alternative Fryzlewicz et.al. (2008) propose a kernel normalized-least-squares estimator which has a closed
form and thus has some advantages over the above kernel estimate for small samples.

(iii) Another example is a tvGARCH(p,q)-process -- see Example~\ref{ExampleTVGARCH}. \hfill $\Box$
\end{example}

We now discuss the derivation of the asymptotic bias, mean squared error, consistency and asymptotic normality of $\hat{{\boldsymbol \theta}}_{T}(u_0)$ for an ``arbitrary'' local minimum distance function $\mathcal{L}_{T}(u_0, {\boldsymbol \theta})$ (keeping in mind the above local conditional likelihood). The results are obtained by approximating $\mathcal{L}_{T}(u_0, {\boldsymbol \theta})$ with
$\tilde{\mathcal{L}}_{T}(u_0, {\boldsymbol \theta})$ which is the same function but with $X_{t,T}$ replaced by its
stationary approximation $\tilde{X}_{t}(u_0)$. Typically both, $\mathcal{L}_{T}(u_0, {\boldsymbol \theta})$
and $\tilde{\mathcal{L}}_{T}(u_0, {\boldsymbol \theta})$ will converge to the same limit-function which we denote
by $\mathcal{L} (u_0, {\boldsymbol \theta})$. Let
\begin{equation*} \label{}
{\boldsymbol \theta}_0 (u_0) := \argmin_{{\boldsymbol \theta} \in \Theta}
\mathcal{L} (u_0, {\boldsymbol \theta}).
\end{equation*}
If the model is correctly specified then typically ${\boldsymbol \theta}_0 (u_0)$ is the true curve. Furthermore, let
\begin{equation*} \label{}
\mathcal{B}_{T}(u_0, {\boldsymbol \theta}) := \mathcal{L}_{T}(u_0, {\boldsymbol \theta}) - \tilde{\mathcal{L}}_{T}(u_0, {\boldsymbol \theta}).
\end{equation*}

The following two results describe how the asymptotic properties of $\hat{{\boldsymbol \theta}}_{T}(u_0)$ can be derived. They should be regarded as a general roadmap and the challenge is to prove the conditions in a specific situation which may be quite difficult.

\begin{theorem} \label{NonlinearConsistency}
(i) Suppose that $\Theta$ is compact with $\boldsymbol{\theta}_0 (u_0) \in Int(\Theta)$, the function $\mathcal{L} (u_{0},{\boldsymbol \theta})$
is continuous in ${\boldsymbol \theta}$ and the minimum ${\boldsymbol \theta}_0 (u_0)$ is unique. If
\begin{equation} \label{Consist1}
\sup_{{\boldsymbol \theta} \in \Theta}
\big|\tilde{\mathcal{L}}_{T}(u_{0},{\boldsymbol \theta}) -
\mathcal{L} (u_{0},{\boldsymbol \theta})\big| \pconv 0,
\end{equation}
and
\begin{equation} \label{ConditionB1}
\sup_{{\boldsymbol \theta} \in \Theta}
\big|\mathcal{B}_{T}(u_0, {\boldsymbol \theta})\big| \pconv 0
\end{equation}
then
\begin{equation} \label{ConsistencyResult1}
\hat{{\boldsymbol \theta}}_{T}(u_0) \pconv {\boldsymbol \theta}_0 (u_0).
\end{equation}
(ii) Suppose in addition that $\mathcal{L} (u,{\boldsymbol \theta})$ and ${\boldsymbol \theta}_{0} (u)$ are uniformly continuous in
$u$ and ${\boldsymbol \theta}$ and the convergence in (\ref{Consist1}) and (\ref{ConditionB1}) is uniformly in $u_0 \in [0,1]$. Then
\begin{equation} \label{ConsistencyResult2}
\sup_{u_0 \in [0,1]} \big|\hat{{\boldsymbol \theta}}_{T}(u_0) - {\boldsymbol \theta}_0 (u_0)\big| \pconv 0.
\end{equation}
\end{theorem}
\proofbegin The proof of (i) is standard -- cf. the proof of Theorem~2 in Dahlhaus and Subba Rao (2006). The proof of (ii) is a straightforward generalization.

\medskip

Note that in (i) all conditions apart from (\ref{ConditionB1}) are conditions on the stationary process $\tilde{X}_{t}(u_0)$ with (fixed) parameter ${\boldsymbol \theta} (u_0)$ and the stationary likelihood / minimum-distance function $\tilde{\mathcal{L}}_{T}(u_0, {\boldsymbol \theta})$. These properties are usually known from existing results on stationary processes. It only remains to verify the condition (\ref{ConditionB1}) which can be done by using the expansion (\ref{ExpansionDerivativeProcesses})  in terms of derivative processes (see the discussion below). (ii) contains a little pitfall: Usually the estimate $\hat{{\boldsymbol \theta}}_{T}(u_0)$ is defined for $u_0=0$ or $u_0=1$ in a different way due to edge-effects. This means that also $\tilde{\mathcal{L}}_{T}(u_0, {\boldsymbol \theta})$ looks different, that is one would usually prefer a uniform convergence result for $u_0 \in (0,1)$ which is more difficult to prove.

Even more interesting and challenging is a uniform convergence result with a rate of convergence. For time varying AR(p)-processes this is stated for a different likelihood in Theorem~\ref{uniform-convergence-Whittle}. We mention that such a result usually requires an exponential bound and maximal inequalities which need to be tailored to the specific model at hand.

We now state the corresponding result on asymptotic normality in case of second order smoothness.  $\nabla$ denotes the derivatives with respect to the $\theta_i$, i.e. $\nabla := \big( \frac {\partial}
{\partial \theta_{i}} \big)_{i=1,\ldots,d}$.

\begin{theorem} \label{NonlinearCLT}
Let ${\boldsymbol \theta}_0 :={\boldsymbol \theta}_0 (u_{0})$. Suppose that $\mathcal{L}_{T}(u_0, {\boldsymbol \theta})$, $\tilde{\mathcal{L}}_{T}(u_0, {\boldsymbol \theta})$ and $\mathcal{L} (u_{0},{\boldsymbol \theta})$ are twice continuously differentiable in ${\boldsymbol \theta}$ with nonsingular matrix
$\Gamma (u_0) := \nabla^{2} \mathcal{L} (u_{0},{\boldsymbol \theta}_0)$. Let further
\begin{equation*} \label{}
\sqrt{b\,T} \; \nabla \tilde{\mathcal{L}}_{T}(u_0, {\boldsymbol \theta}_0) \dconv \mathcal{N} \big(0, V (u_0)\big)
\end{equation*}
with some sequence $b=b_{T}$ where $b \rightarrow 0$ and $b\,T \rightarrow \infty$ (the definition of $b$ is part of the definition
of the likelihood -- it is usually some bandwidth) and
\begin{equation*}
\sup_{{\boldsymbol \theta} \in \Theta}
\big|\nabla^{2} \tilde{\mathcal{L}}_{T}(u_{0},{\boldsymbol \theta}) -
\nabla^{2} \mathcal{L} (u_{0},{\boldsymbol \theta})\big| \pconv 0.
\end{equation*}
If in addition
\begin{equation} \label{BiasConditionB}
\sqrt{bT} \Big(\Gamma (u_0)^{-1} \, \nabla \mathcal{B}_{T}(u_0, {\boldsymbol \theta}_0) - \frac {b^{2}} {2}\,{\boldsymbol \mu}^{0}(u_{0}) \Big) = o_{p}(1)
\end{equation}
with some ${\boldsymbol \mu}^{0}(\cdot)$ (to be specified below -- cf.(\ref{DetailBias})) and
\begin{equation} \label{ConditionB2}
\sup_{{\boldsymbol \theta} \in \Theta}
\big|\nabla^{2} \mathcal{B}_{T}(u_0, {\boldsymbol \theta})\big| \pconv 0
\end{equation}
then
\begin{equation} \label{AssertionCLT}
\sqrt{b\,T} \; \Big(\hat{{\boldsymbol \theta}}_{T}(u_0) - {\boldsymbol \theta}_0 (u_0) + \frac {b^{2}} {2}  {\boldsymbol \mu}^{0} (u_0)\Big) \dconv \mathcal{N} \Big(0, \Gamma (u_0)^{-1} \,V(u_0) \,\Gamma(u_0)^{-1} \Big).
\end{equation}
\end{theorem}

\proofbegin The usual Taylor-expansion of $\nabla \mathcal{L}_{T}(u_0, {\boldsymbol \theta})$ around ${\boldsymbol \theta}_0$ yields
\begin{equation} \label{TaylorNonlinear}
\sqrt{b\,T} \; \big(\hat{{\boldsymbol \theta}}_{T}(u_0)  - {\boldsymbol \theta}_0   +
\Gamma(u_0)^{-1} \nabla   \mathcal{B}_{T}(u_0, {\boldsymbol \theta}_0) \big)  = - \sqrt{b\,T} \; \Gamma(u_0)^{-1} \,\nabla \tilde{\mathcal{L}}_{T}(u_0, {\boldsymbol \theta}_0) + o_{p}(1).
\end{equation}
The result then follows immediately.

\begin{remark} \label{}
(i) Again the first two conditions are conditions on the stationary process $\tilde{X}_{t}(u_0)$ with (fixed) parameter ${\boldsymbol \theta} (u_0)$ and the stationary likelihood / minimum-distance function $\tilde{\mathcal{L}}_{T}(u_0, {\boldsymbol \theta})$ which are usually known from existing results on stationary processes.\\[6pt]
(ii) Of course an analogous result also holds under different smoothness conditions and with other rates
than $b^{2}$ in (\ref{BiasConditionB}) and (\ref{AssertionCLT}).\\[6pt]
(iii) Under additional regularity conditions one can usually prove that the same expansion as in (\ref{TaylorNonlinear}) also holds
for the moments, leading to
\begin{equation} \label{AsymptoticBias}
\mean \, \hat{{\boldsymbol \theta}}_{T}(u_0) = {\boldsymbol \theta}_0 (u_0) - \frac {b^{2}} {2} \, {\boldsymbol \mu}^{0} (u_0) + o (b^{2})
\end{equation}
and
\begin{equation} \label{AsymptoticVariance}
\var \big( \hat{{\boldsymbol \theta}}_{T}(u_0) \big) = \frac {1} {bT} \,\Gamma (u_0)^{-1} \,V(u_0) \,\Gamma(u_0)^{-1} + o \Big( \frac {1} {bT} \Big)
\end{equation}
(note that (\ref{TaylorNonlinear}) is a stochastic expansion which does not automatically imply
these moment relations). The proof of these properties is usually not easy. \hfill $\Box$
\end{remark}
\begin{example} [Kernel-type local likelihoods] \label{BiasExpressionLocalLikelihood}
We now return to the local conditional likelihood (\ref{localConditionalLikelihood}) as a special case and provide some heuristics on how to calculate the above terms (in particular the bias ${\boldsymbol \mu}^{0} (u_0)$).  We stress that in the concrete situation where a specific model is given the exact proof usually goes along the same lines but the details may be quite challenging.

Suppose that the local likelihood of the stationary process $\tilde{X}_t (u_0)$ converges in probability to
\begin{equation*} \label{}
\mathcal{L} (u_0, {\boldsymbol \theta}) := \lim_{T \rightarrow \infty} \tilde{\mathcal{L}}_{T}(u_0, {\boldsymbol \theta}) = \lim_{t \rightarrow \infty} \mean \, \tilde{\ell}_t (u_0,{\boldsymbol \theta}).
\end{equation*}

Usually we have $X_{t,T} = \tilde{X}_t (t/T) + O_{p}(T^{-1})$ and
\begin{equation*} \label{}
\mean \,\nabla \ell_{t,T} ({\boldsymbol \theta}) = \mean \,\nabla \tilde{\ell}_t \big(\frac {t} {T},{\boldsymbol \theta}\big) + o\big((bT)^{-1/2}\big)
= \nabla \mathcal{L} \big(\frac {t} {T}, {\boldsymbol \theta}\big) + o\big((bT)^{-1/2}\big)
\end{equation*}
uniformly in $t$. A Taylor-expansion then leads in the case $b^{3}= o((bT)^{-1/2})$ with the symmetry of the kernel $K$ to
\begin{align} \label{ExpectationLocalConditionalLikelihood}
\mean \, \nabla \mathcal{L}_{T}(u_0, {\boldsymbol \theta}) & =
\frac{1}{bT} \sum^{T}_{t=1} \, K \Big(\frac{u_0-t/T}{b}\Big) \, \nabla  \mathcal{L} \big(\frac {t} {T}, {\boldsymbol \theta}\big) + o\big((bT)^{-1/2}\big) \nonumber\\
& =  \nabla \mathcal{L} (u_0, {\boldsymbol \theta}) + \Big[\frac {\partial} {\partial u} \nabla \mathcal{L} (u_0, {\boldsymbol \theta}) \Big]
\, \frac{1}{bT} \sum^{T}_{t=1} \, K \Big(\frac{u_0-t/T}{b}\Big) \; \Big(\frac {t} {T} - u_0\Big)  + \nonumber\\
& \qquad \qquad \qquad \!\! + \frac {1} {2} \Big[\frac {\partial^{2}} {\partial u^{2}} \nabla \mathcal{L} (u_0, {\boldsymbol \theta}) \Big]
\, \frac{1}{bT} \sum^{T}_{t=1} \, K \Big(\frac{u_0-t/T}{b}\Big) \;\Big(\frac {t} {T} - u_0\Big)^{2}  + o\big((bT)^{-1/2}\big) \nonumber\\
& = \nabla \mathcal{L} (u_0, {\boldsymbol \theta}) +  \frac {1} {2} \, b^{2} \, d_K \, \frac {\partial^{2}} {\partial u^{2}} \nabla \mathcal{L} (u_0, {\boldsymbol \theta}) + o\big((bT)^{-1/2}\big)
\end{align}
with $d_K := \int x^{2} K(x) \,dx$. Since $\mean \, \nabla \tilde{\mathcal{L}}_{T}(u_0, {\boldsymbol \theta}) = \nabla \mathcal{L} (u_0, {\boldsymbol \theta}) + o\big((bT)^{-1/2}\big)$
this leads with (\ref{BiasConditionB}) to the bias term
\begin{equation} \label{DetailBias}
{\boldsymbol \mu}^{0}(u_{0})  =  d_K \, \Gamma (u_0)^{-1} \,  \frac {\partial^{2}} {\partial u^{2}} \nabla \mathcal{L} \big(u, {\boldsymbol \theta}_0 (u_0)\big) \rfloor_{u=u_{0}} =: d_k \, {\boldsymbol \mu}(u_{0})
\end{equation}
Let $\boldsymbol {\theta}_0 := \boldsymbol {\theta}_0 (u_0)$. If the model is correctly specified it usually can be shown that $\nabla \tilde{\ell}_t (u_0,{\boldsymbol \theta}_0)$
is a martingale difference sequence and the condition of the Lindeberg martingale central
limit theorem are fulfilled leading to
\begin{equation*} \label{}
\sqrt{b\,T} \; \nabla \tilde{\mathcal{L}}_{T}(u_0, {\boldsymbol \theta}_0) \dconv \mathcal{N} \Big(0, v_K \, \mean  \big( \nabla \tilde{\ell}_t (u_0,{\boldsymbol \theta}_0)\big) \big( \nabla \tilde{\ell}_t (u_0,{\boldsymbol \theta}_0)\big)'\Big)
\end{equation*}
with $v_K = \int K(x)^{2} dx$. Furthermore, if the model is correctly specified we usually have
\begin{equation*} \label{}
\mean  \big( \nabla \tilde{\ell}_t (u_0,{\boldsymbol \theta}_0)\big)\big( \nabla \tilde{\ell}_t (u_0,{\boldsymbol \theta}_0)\big)' = \nabla^{2} \mathcal{L} (u_{0},{\boldsymbol \theta}_0) = \Gamma (u_0)
\end{equation*}
that is
\begin{equation} \label{AssertionCLT2}
\sqrt{b\,T} \; \Big(\hat{{\boldsymbol \theta}}_{T}(u_0) - {\boldsymbol \theta}_0 (u_0) + \frac {b^{2}} {2} \, d_K \, \Gamma (u_0)^{-1} \,  \frac {\partial^{2}} {\partial u^{2}} \nabla \mathcal{L} (u_0, {\boldsymbol \theta}_0)\Big) \dconv \mathcal{N} \Big(0, v_K \,\Gamma(u_0)^{-1} \Big).
\end{equation}
If we are able to prove in addition the formulas (\ref{AsymptoticBias}) and (\ref{AsymptoticVariance}) on the asymptotic bias and variance we obtain the same formula for the asymptotic mean squared error as in (\ref{AsymptoticMSE}) with $\tau(u_0)={\rm tr} \{\Gamma(u_0)^{-1}\}$ and $\mu (u_0)^{2}$ replaced by $\|{\boldsymbol \mu}(u_{0})\|^{2}$ where
${\boldsymbol \mu} (u_{0}) =  \Gamma (u_0)^{-1} \,  \frac {\partial^{2}} {\partial u^{2}}
\nabla \mathcal{L} (u_0, {\boldsymbol \theta}_0)$. As in Remark~\ref{MinimalMSE} this leads to the optimal segment length and the optimal mean squared error. The implications for non-rescaled processes are the same as in Remark~\ref{ImplicationsNonRescaled}.
\hfill $\Box$
\end{example}

We now present three examples where the above results have been proved explicitly.

\begin{example} [Local Whittle estimates] \label{LocalWhittleEstimatesBiasVariance}
The first example are local Whittle estimates on segments $\hat{\boldsymbol{\theta}}_{T}^{W}(u_0)$ obtained by
minimizing $\mathcal{L}_{T}^{W}(u_0, \boldsymbol{\theta})$ (cf.(\ref{NonparametricWhittleLikelihood})). In case
of a tvAR(p)-process $\hat{\boldsymbol{\theta}}_{T}^{W}(u_0)$ is exactly the local Yule-Walker estimate defined
in (\ref{LocalYuleWalkerEst}) with the covariance-estimates given in
(\ref{TaperedCovarianceEstimate}). $\mathcal{L}_{T}^{W}(u, \boldsymbol{\theta})$ is not exactly a local
conditional likelihood as defined in (\ref{localConditionalLikelihood}) but approximately (in the
same sense as ${\hat c}_{T}(u_0,k)$ from (\ref{TaperedCovarianceEstimate}) is an approximation to
the kernel covariance estimate). For that reason the above heuristics also applies to this estimate
and can be made rigorous.

In Dahlhaus and Giraitis (1998), Theorem 3.1 and 3.2, bias and asymptotic normality of $\hat{\boldsymbol{\theta}}^{W}_{T}(u_0)$ have been derived rigorously including a derivation of the variance and the mean squared error as given in (\ref{AsymptoticBias}) and (\ref{AsymptoticVariance}) (i.e. not only the stochastic expansion in (\ref{TaylorNonlinear})). We mention that therefore also the results on the optimal kernel and bandwidth in (\ref{OptimalKernel})  and (\ref{OptimalBandwidth}) apply to this situation.

In the present situation we have (cf. Dahlhaus and Giraitis (1998), (3.7)))
\begin{equation*} \label{}
{\mathcal L}(u,\boldsymbol{\theta})=\frac{1}{4 {\pi}} \int_{-\pi}^{\pi} \Big\{ \log 4
{\pi}^2 f_{\boldsymbol{\theta}}(\lambda)
+\frac{f(u,\lambda)}{f_{\boldsymbol{\theta}}(\lambda)}\Big\} \, d
\lambda.
\end{equation*}
Therefore
\begin{equation*} \label{}
 \frac {\partial^{2}} {\partial u^{2}} \nabla {\mathcal L}(u_0,\boldsymbol{\theta})
 = \frac{1}{4 {\pi}} \int_{-\pi}^{\pi}   \nabla f_{\boldsymbol{\theta}}(\lambda)^{-1} \frac {\partial^{2}} {\partial u^{2}} f(u_0,\lambda) \, d
\lambda
\end{equation*}

and in the correctly specified case where $f(u,\lambda)=f_{{\boldsymbol \theta}_0 (u)}(\lambda)$
\begin{equation*} \label{}
\Gamma (u_0) = \nabla^{2} \mathcal{L} (u_{0},{\boldsymbol \theta}_0) = \frac{1}{4{\pi}} \int_{-\pi}^{\pi} \big(\nabla\log f_{\boldsymbol{\theta}_0}\big)
\big(\nabla \log f_{\boldsymbol{\theta}_0}\big)' \,
d\lambda
\end{equation*}
leading to the asymptotic bias ${\boldsymbol \mu} (u_{0})$ in (\ref{DetailBias}) and the asymptotic variance in the central limit theorem (\ref{AssertionCLT2}). A uniform convergence result for $\hat{\boldsymbol{\theta}}^{W}_{T}(u_0)$ is stated in Theorem~\ref{uniform-convergence-Whittle}.
\hfill $\Box$
\end{example}

\begin{example} [tvAR(p)-processes] \label{LocalARestimatesDahGir}
In the special case of a Gaussian tvAR(p)-process the exact results for the local Yule-Walker estimates
(\ref{LocalYuleWalkerEst}) follow as a special case from the above results on local Whittle estimates (see also Section 2 in Dahlhaus and Giraitis, 1998, where
tvAR(p)-processes are discussed separately). In that case we have with $R(u)$
and $r(u)$ as in (\ref{TheoreticalYuleWalker}) that $\Gamma (u) = \frac {1} {\sigma^{2}(u)} R(u)$. Furthermore
\begin{equation*} \label{}
\nabla \mathcal{L} \big(u, {\boldsymbol \theta}\big) = \frac {1} {\sigma^{2}} \big[ R(u)\, {\boldsymbol \alpha} + r(u) \big]
\end{equation*}
which implies
\begin{equation*} \label{}
{\boldsymbol \mu}(u_{0})  =   R(u_0)^{-1} \, \Big[ \Big(\frac {\partial^{2}} {\partial u^{2}} R(u)\Big) {\boldsymbol \alpha}(u_0) + \Big(\frac {\partial^{2}} {\partial u^{2}} r(u)\Big) \Big]_{u=u_0}.
\end{equation*}
We conjecture that exactly the same asymptotic results hold for the conditional likelihood estimate obtained by minimizing
\begin{flalign*} \label{}
\quad  &\mathcal{L}_{T}^{C}(u_0, {\boldsymbol \theta}) :=
\frac{1}{T} \sum^{T}_{t=1}  \;\frac{1}{\,b}\, K \Big(\frac{u_0-t/T}{b}\Big) \; \bigg[\frac{1}{2} \, \textrm{log}
\big( 2\pi\,\sigma^2  \big) + \frac{1}{2\,\sigma^2} \;\Big(
X_{t,T} + \sum_{j=1}^{p}\alpha_{j}\; X_{t-j,T} \Big)^{\!2} \bigg].&\Box
\end{flalign*}
\end{example}

We now introduce derivative processes. The key idea in the proofs of Dahlhaus and Giraitis (1998) is to use at time $u_0 \in (0,1)$ the stationary approximation $\tilde{X}_{t}(u_0)$ (there denoted by $Y_{t}$) to the original process $X_{t,T}$ and to calculate the bias resulting from the use of this approximation. As in Dahlhaus and Subba Rao (2006) we now extend this idea leading to the Taylor-type expansion (\ref{ExpansionDerivativeProcesses}) which is an expansion of the original process in terms of (usually ergodic) stationary processes called derivative processes. This expansion is a powerful tool since all techniques for stationary processes including the ergodic theorem may be applied for the local investigation of the nonstationary process $X_{t,T}$. The use of this expansion and of derivative processes in general leads to a general structure of the proofs and simplifies the derivations a lot.

We start with the simple example of a tvAR(1)-process since in this case everything can be calculated directly.
Then $X_{t,T}$ is defined by $X_{t,T} + \alpha_{1} (t/T)  X_{t-1,T} = \varepsilon_{t}, t\in\znum$ and
the stationary approximation $\tilde{X}_{t} (u_0)$ at time $u_0=t_0/n$ by
$\tilde{X}_{t} (u_0) + \alpha_{1} (u_0) \tilde{X}_{t} (u_0)= \varepsilon_{t}, t\in\znum$. Repeated plug-in yields
under suitable regularity conditions (for a rigorous argument see the proof of Theorem 2.3 in Dahlhaus (1996a))
\begin{align}
X_{t,T} & = \sum_{j=0}^{\infty} (-1)^{j} \, \Big[\prod_{k=0}^{j-1} \alpha_{1} \Big(\frac {t-k} {T}\Big) \Big] \, \varepsilon_{t-j}
= \sum_{j=0}^{\infty} (-1)^{j} \, \alpha_{1}\Big(\frac {t} {T}\Big)^{j} \, \varepsilon_{t-j} \,+\, O_{p}\Big(\frac {1} {T}\Big)\label{Derivative_tvAR0}\\[6pt]
& = \tilde{X}_{t} \Big(\frac {t} {T}\Big) \,+\, O_{p}\Big(\frac {1} {T}\Big)
= \tilde{X}_{t} (u_0) + \Big( \frac {t} {T} - u_0\Big) \,\frac {\partial \tilde{X}_{t} (u)} {\partial u}\rfloor_{u=u_{0}} + O_{p}\Big(\frac {1} {T}\Big). \label{Derivative_tvAR1}
\end{align}
We have in the present situation
\begin{equation*} \label{}
\frac {\partial \tilde{X}_{t} (u)} {\partial u} = \sum_{j=0}^{\infty} (-1)^{j} \, \frac {\partial \alpha_{1}(u)^{j}} {\partial u} \, \varepsilon_{t-j} = \sum_{j=0}^{\infty} (-1)^{j} \big[ j \,\alpha_{1}(u)^{j-1} \alpha_{1}(u)'\big] \, \varepsilon_{t-j}
\end{equation*}
that is $\frac {\partial \tilde{X}_{t} (u)} {\partial u}$ is a stationary ergodic process in $t$ with
$\Big|\frac {\partial \tilde{X}_{t} (u)} {\partial u}\Big| \le \sum_{j=1}^{\infty} j \rho^{j-1} | \varepsilon_{t-j} |$
where $|\rho| < 1$. In the same way we have
\begin{align} \label{ExpansionDerivativeProcesses}
X_{t,T} =
\tilX_{t}(u_{0}) & + \Big(\frac{t}{T} - u_{0}\Big)\frac{\partial \tilX_{t}(u)}
{\partial u}\rfloor_{u=u_{0}} +
\frac{1}{2}\Big(\frac{t}{T} - u_{0}\Big)^{\!2}
\frac{\partial^{2}\tilX_{t}(u)}{\partial u^{2}}\rfloor_{u=u_{0}} \nonumber\\
& + O_{p}\Big(\big(\frac{t}{T} - u_{0}\big)^{\!3} + \frac{1}{T}\Big)
\end{align}

with the second order derivative process $
\frac{\partial^{2}\tilX_{t}(u)}{\partial u^{2}}\rfloor_{u=u_{0}}$ which is defined analogously. It is not difficult to prove existence and uniqueness in a rigorous sense.

For general tvAR(p)-processes the same results holds -- however, it is difficult in that case to write the derivative process in explicit form. It is interesting to note that the
derivative process fulfills the equation
\begin{eqnarray*} \label{FunctionalRelationDerivative}
\frac{\partial \tilde{X}_{t}(u)}{\partial u} + \sum_{j=1}^{p} \Big(
\alpha_{j}(u)\, \frac{\partial \tilde{X}_{t-j}(u)}{\partial u} + \alpha_{j}^{\prime}(u) \, \tilde{X}_{t-j}(u) \Big) = \frac{\partial
\sigma(u)}{\partial u}\, \varepsilon_{t}.
\end{eqnarray*}
where $\alpha_{j}^{\prime}(u)$ denotes the derivative of $\alpha_{j}(u)$ with respect to
$u$. This is formally obtained by differentiating both sides of equation
(\ref{ARpApproximation}). Furthermore, it can be shown that this equation
system uniquely defines the derivative process.

We are convinced that the expansion (\ref{ExpansionDerivativeProcesses}) and equation systems like  (\ref{FunctionalRelationDerivative}) can be established for several other locally stationary time series models. As mentioned above the important point is that (\ref{ExpansionDerivativeProcesses}) is an expansion in terms of stationary processes.

In the next example we show how derivative processes are used for deriving the properties of local likelihood estimates.

\begin{example} [tvARCH-processes] \label{ExampleTVARCH2} The definition of the processes $X_{t,T}$ and $\tilde{X}_t (u_0)$ has
been given above in (\ref{eq:tvarch-p})  and (\ref{eq:tvarch-rescale}) and of the local likelihood in (\ref{LocalLikelihoodTVArch})
and (\ref{localConditionalLikelihood}). In Dahlhaus and Subba Rao (2006), Theorem 2 and 3, consistency and asymptotic normality have been established for the resulting estimate and in particular (\ref{AssertionCLT2}) has been proved. Derivative processes play a major role in the proofs and we briefly indicate how they are used. First, existence and uniqueness of the derivative
processes have been proved including the Taylor-type expansion for the process $X_{t,T}^{2}$:
\begin{align}
\label{eq:cor-tay}
X_{t,T}^{2} =
\tilX_{t}(u_{0})^{2} & + \Big(\frac{t}{T} - u_{0}\Big)\, \frac{\partial \tilX_{t}(u)^{2}}
{\partial u}\rfloor_{u=u_{0}}  +
\frac{1}{2}\Big(\frac{t}{T} - u_{0}\Big)^{\!2}\,
\frac{\partial^{2}\tilX_{t}(u)^{2}}{\partial u^{2}}\rfloor_{u=u_{0}} \nonumber\\[2pt]
 & + O_{p}\Big(\big(\frac{t}{T} - u_{0}\big)^{\!3} + \frac{1}{T}\Big)
\end{align}
(in this model we are working with $X_{t,T}^{2}$ rather than $X_{t,T}$ since $X_{t,T}^{2}$ is uniquely determined).
Furthermore, $\frac{\partial \tilX_{t}(u)^{2}}{\partial u}$
is almost surely the
unique solution of the  equation
\begin{eqnarray}
\label{eq:tvarch-der}
\frac{\partial\tilX_{t}(u)^{2}}{\partial u}
=  \Big(\alpha_{0}^{\prime}(u) +
\sum_{j=1}^{\infty} \alpha_{j}^{\prime}(u)\, \tilX_{t-j}(u)^{2} +
\sum_{j=1}^{\infty}\alpha_{j}(u)\, \frac{\partial \tilX_{t-j}(u)^{2}}{\partial u}
\Big)\, Z_{t}^{2}
\end{eqnarray}
which can formally be obtained by differentiating  (\ref{eq:tvarch-rescale}). By taking the second derivative of this expression we obtain a similar expression for the second derivative $\frac{\partial^{2}\tilX_{t}(u)^{2}}{\partial u^{2}}$ etc.

A key step in the above proofs is the derivation of (\ref{BiasConditionB}) and of the bias
term ${\boldsymbol \mu}^{0}(\cdot)$ in this situation. We briefly sketch this. We have with ${\boldsymbol \theta}_0 = {\boldsymbol \theta}_0(u_0)$
\begin{equation*} \label{}
\nabla \mathcal{B}_{T}(u_0, {\boldsymbol \theta}_0) =
\frac{1}{bT} \sum^{T}_{t=1} \, K \Big(\frac{u_0-t/T}{b}\Big) \,   \big( \nabla \ell_{t,T} ({\boldsymbol \theta}_0) -
\nabla \tilde{\ell}_t (u_0,{\boldsymbol \theta}_0) \big).
\end{equation*}
First $\nabla \ell_{t,T} ({\boldsymbol \theta}_0)$ is replaced by $\nabla \tilde{\ell}_{t} (t/T,{\boldsymbol \theta}_0)$ where we omit details (this works since $X_{t,T}^{2}$ is approximately the same as $\tilde{X}_{t}^{2}(t/T)$). Then a Taylor-expansion is applied:
\begin{align}
\label{eq:l-tayl-u}
 \nabla \tilde{\ell}_{t}(\frac{t}{T},{\boldsymbol \theta}_0) & -
\nabla \tilde{\ell}_{t}(u_{0},{\boldsymbol \theta}_0) =
\Big(\frac{t}{T} - u_{0}\Big)\,
\frac{\partial \nabla \tilde{\ell}_{t}(u,{\boldsymbol \theta}_0) }{\partial u}
\rfloor_{u=u_{0}}\nonumber\\[4pt]
& \quad + \frac {1} {2}\Big(\frac{t}{T} - u_{0}\Big)^{\!2} \,
\frac{\partial^{2}\nabla \tilde{\ell}_{t}(u,{\boldsymbol \theta}_0) }
{\partial u^{2}}
\rfloor_{u=u_{0}}
  +  \frac {1} {6}\, \Big(\frac{t}{T} - u_{0}\Big)^{\!3} \,
\frac{\partial^{3}\nabla \tilde{\ell}_{t}(u,{\boldsymbol \theta}_0) }{\partial u^{3}}
\rfloor_{u=\tilde{U}_{t}}
\end{align}
with a random variable $\tilde{U}_{t}\in (0,1]$. The breakthrough now is that $\frac{\partial \nabla \tilde{\ell}_{t}(u,{\boldsymbol \theta}_0) }{\partial u}$ can be written explicitly in terms of the derivative process $\frac{\partial \tilX_{t}(u)^{2}}
{\partial u}$ and of the process $\tilX_{t}(u)^{2}$, that is we obtain with the formula for the total derivative
\begin{eqnarray*}
\frac{\partial \nabla \tilde{\ell}_{t}(u,{\boldsymbol \theta}_0) }{\partial u}
=\sum_{j=0}^{p}\bigg(\frac{\partial }{\partial \tilde{X}_{t-j}(u)^{2}}
\bigg[\frac{\nabla w_{t}(u,{\boldsymbol \theta}_{0})}{w_{t}(u,{\boldsymbol \theta}_{0})} -
\frac{\tilde{X}_{t}(u)^{2}\nabla w_{t}(u,{\boldsymbol \theta}_{0})}{w_{t}(u,{\boldsymbol \theta}_{0})^{2}}\bigg] \times \frac{\partial \tilde{X}_{t-j}(u)^{2}}{\partial u}\bigg),
\end{eqnarray*}
where  $w_{t}(u,{\boldsymbol \theta}) = c_{0}({\boldsymbol \theta}_{0})  + \sum_{j=1}^{\infty}c_{j}({\boldsymbol \theta})\tilde{X}_{t-j}(u)^{2}$
(the same holds true for the higher order terms).
In particular $\frac{\partial \nabla \tilde{\ell}_{t}(u,{\boldsymbol \theta}_0) }{\partial u}$ is a stationary process with constant mean. Due to the symmetry of the kernel we therefore obtain after some lengthy but straightforward calculations
\begin{equation} \label{BiasConditionB2}
\sqrt{bT} \Big(\Gamma (u_0)^{-1} \, \nabla \mathcal{B}_{T}(u_0, {\boldsymbol \theta}_0) - \frac {b^{2}} {2}\,d_K \, \Gamma (u_0)^{-1} \,  \frac {\partial^{2}} {\partial u^{2}} \nabla \mathcal{L} \big(u, {\boldsymbol \theta}_0 \big) \rfloor_{u=u_{0}} \Big) = o_{p}(1).
\end{equation}

A very simple example is the tvARCH$(0)$ process
\begin{equation*}
X_{t,T} = \sigma_{t,T}Z_{t}, \quad \sigma_{t,T}^{2} =
\alpha_{0}(\frac{t}{T}).
\end{equation*}
In this case $\frac{\partial \tilX_{t}(u)^{2}}{\partial u}
= \alpha_{0}^{\prime}(u) \, Z_{t}^{2}$ and we have
\begin{equation*}
\frac{\partial^{2}\nabla \mathcal{L}(u,{\boldsymbol \alpha}_{u_{0}})
}{\partial u^{2}}\rfloor_{u=u_{0}} =
-\frac{1}{2}\frac{\alpha_{0}^{\prime \prime}(u_{0})}{\alpha_{0}(u_{0})^{2}}
\quad \textrm{ and } \quad \Sigma(u_{0}) =
\frac{1}{2\,\alpha_{0}(u_{0})^{2}}
\end{equation*}
that is $\mu(u_{0}) = - \alpha_{0}^{\prime \prime}(u_{0})$. This is another example which illustrates how the bias is linked to the
nonstationarity of the process - if the process were stationary the
derivatives of $\alpha_{0}(\cdot)$ would be zero causing the bias also to
be zero. The formula (\ref{OptimalBandwidth}) for the optimal bandwidth leads in this case to
\begin{equation*}
b_{opt}(u_0)= \bigg[\frac {2 v_K} {d_K^{2}}\bigg]^{1/5}
\bigg[\frac{\alpha_{0}(u_{0})}{\alpha_{0}^{\prime\prime}(u_{0})} \bigg]^{2/5} T^{-1/5}
\end{equation*}
leading to a large bandwidth if $\alpha_{0}^{\prime\prime}(u_{0})$ is
small and vice versa. As in Remark~\ref{ImplicationsNonRescaled} this can be ``translated'' to the non-rescaled case.
\hfill $\Box$
\end{example}

\begin{example} [tvGARCH-processes] \label{ExampleTVGARCH}
A tvGARCH$(p,q)$-process satisfies the following representation
\begin{align}
\label{eq:tvgarch}
X_{t,T} &= \sigma_{t,T}\,Z_{t} \nonumber\\
\textrm{ where } \quad \sigma_{t,T}^{2} &= \alpha_{0}(\frac{t}{T}) + \sum_{j=1}^{p}\alpha_{j}(\frac{t}{T})\,X_{t-j,T}^{2} +
\sum_{i=1}^{q}\beta_{i}(\frac{t}{T})\,
\sigma_{t-i,T}^{2}, \qquad \;
\end{align}
where $\{Z_{t}\}$ are iid random variables with $\mean Z_{t} = 0$ and $\mean Z_{t}^{2}=1$. The corresponding stationary approximation at
time $u_{0}$ is given by
\begin{align}
\label{eq:tvgarch-rescale}
  \tilX_{t}(u_{0})  & = \sigma_{t}(u_{0})\, Z_{t}\quad \textrm{for } t \in \znum \nonumber\\
\textrm{ where } \quad   \sigma_{t}(u_{0})^{2}   & = \alpha_{0}(u_{0}) +
\sum_{j=1}^{p}\alpha_{j}(u_{0})\, \tilX_{t-j}(u_{0})^{2} + \sum_{i=1}^{q}\beta_{i}(u_{0})\,\sigma_{t-i}(u_{0})^{2}.
\end{align}
Under the condition that $\sup_{u}\big(\sum_{j=1}^{p}\alpha_{j}(u) + \sum_{i=1}^{q}\beta_{i}(u)\big) < 1$ Subba Rao (2006), Section 5, has shown that
$X_{t,T}^{2} = \tilX_{t}(u_{0})^{2} + O_{p}(|\frac{t}{T} - u_{0}| + \frac{1}{T})$. To obtain estimators of the parameters $\{\alpha_{j}(\cdot)\}$
and $\{\beta_{i}(\cdot)\}$ an approximation of the conditional quasi-likelihood is used, which is constructed as if the
innovations $\{Z_{t}\}$ were Gaussian. As the infinite past is unobserved, an observable approximation of the conditional quasi-likelihood is
\begin{equation}
\label{eq:garch-likelihood}
\ell_{t,T}({\boldsymbol \theta}) = \frac{1}{2} \log w_{t,T}({\boldsymbol \theta}) +
\frac{X_{t,T}^{2}}{2 \,w_{t,T}({\boldsymbol \theta})}
\textrm{ with }  w_{t,T}({\boldsymbol \theta}) = c_{0}({\boldsymbol \theta}) + \sum_{j=1}^{t-1}c_{j}({\boldsymbol \theta})X_{t-j,T}^{2},
\end{equation}
where a recursive formula for $c_{j}({\boldsymbol \theta})$ in terms of the parameters of interest,
$\{\alpha_{j}\}$ and $\{\beta_{i}\}$, can be found in Berkes et.al. (2003).
Given that the derivatives of the time varying GARCH parameters exist we can
formally differentiate (\ref{eq:tvgarch-rescale}) to obtain
\begin{align*}
\frac{\partial \tilde{X}_{t}(u)^{2}}{\partial u} &= \frac{\partial \sigma_{t}(u)^{2}}{\partial u}Z_{t}^{2} \\
\frac{\partial \sigma_{t}(u)^{2}}{\partial u} &= \alpha_{0}^{\prime}(u) +
\sum_{j=1}^{p}\Big(\alpha_{j}^{\prime}(u)\,\tilde{X}_{t-j}(u)^{2} +
\alpha_{j}(u)\,\frac{\partial \tilde{X}_{t-j}(u)^{2}}{\partial u} \Big)  \\
& \qquad \qquad \!+ \sum_{i=1}^{q}\Big(\beta_{i}^{\prime}(u)\,\sigma_{t-i}(u)^{2}
+ \beta_{i}(u)\,\frac{\partial \sigma_{t-i}(u)^{2}}{\partial u} \Big).
\end{align*}
Subba Rao (2006) has shown that one can represent the above as a state-space representation which almost surely has a unique
solution which is the derivative of $\tilde{X}_{t}(u)^{2}$ with respect to $u$.
Thus $X_{t,T}^{2}$ satisfies the expansion in (\ref{eq:cor-tay}).
Moreover, Fryzlewicz and
Subba Rao (2011) show geometric $\alpha$-mixing of the tvGARCH process. Using these results and
under some technical assumptions it can be shown that Theorem \ref{NonlinearConsistency}(i) and Theorem
\ref{NonlinearCLT} hold for the local approximate conditional quasi-likelihood estimator. In particular,
a result analogous to (\ref{BiasConditionB2}) holds true, where
\begin{equation*}
\mathcal{L}(u,{\boldsymbol \theta}) = \mean \Big( \log \big(c_{0}({\boldsymbol \theta}) + \sum_{j=1}^{\infty}
c_{j}({\boldsymbol \theta})\tilde{X}_{t-j}(u) \big)\Big)
 + \mean\bigg( \frac{\tilde{X}_{t}(u)^{2}}{c_{0}({\boldsymbol \theta}) + \sum_{j=1}^{\infty}c_{j}({\boldsymbol \theta})\tilde{X}_{t-j}(u)^{2}}\bigg).
\end{equation*}
Amado and Ter\"{a}svirta (2011) investigate parametric tvGARCH-models where the time varying parameters are modeled with the logistic transition function -- see Section~\ref{DeepExample}.6.
\hfill $\Box$
\end{example}

Similar methods as described in this section have also been applied in Koo and Linton (2010) who investigate semiparametric estimation of locally stationary diffusion models. They also prove a central limit theorem with a bias term as in (\ref{AssertionCLT}). In their proofs they use the stationary approximation $\tilde{X}_{t}(u_0)$ and the Taylor-type expansion (\ref{ExpansionDerivativeProcesses}). Vogt (2011) investigates nonlinear nonparametric models allowing for locally stationary regressors and a regression function that changes smoothly over time.

\section{A general definition, linear processes and time varying spectral densities} \label{SectionSpectralDensity}

The intuitive idea for a general definition is to require that
locally around each rescaled time point $u_0$ the process $\{X_{t,T}\}$ can be approximated by a stationary process $\{\tilde{X}_t (u_0)\}$ in a stochastic sense by using the property (\ref{AR-LocalApproximation}) (cf. Dahlhaus and Subba Rao, 2006). Vogt (2011) has formalized this by requiring that for each $u_0$ there exists a stationary process $\tilde{X}_{t}(u_0)$ with
\begin{equation} \label{GeneralDefinitionLSP}
\|X_{t,T} - \tilde{X}_{t}(u_0)\big| \le \Big(\Big|\frac{t}{T}-u_0\Big| + \frac{1}{T} \Big) \, U_{t,T} (u_0)
\end{equation}
where $U_{t,T} (u_0)$ is a positive stochastic process fulfilling some uniform moment conditions. However up to now no general theory exists based on such a general definition.

In the following we move on towards a general theory for linear locally stationary processes.
In some cases we even assume Gaussianity or use Gaussian likelihood methods and their approximations. In this situation a fairly  general theory can be derived in which
parametric and nonparametric inference problems, goodness of fit tests, bootstrap procedures etc
can be treated in high generality. We use a general definition tailored for linear processes which however implies (\ref{GeneralDefinitionLSP}).

\textbf{Definition of linear locally stationary processes:} We give this definition in
terms of the time varying MA($\infty$) representation
\begin{equation*} \label{}
X_{t,T} = \mu\big(\frac{t}{T}\big) + \sum^{\infty}_{j = -\infty} a_{t,T} (j) \, \varepsilon_{t-j} \qquad
\mbox{where} \qquad  a_{t,T}(j) \approx \,a(\frac {t} {T},j)
\end{equation*}
with coefficient functions $a(\cdot,j)$ which need to fulfill additional regularity function
(dependent on the result to be proved -- details are provided below). In several papers of
the author instead the time varying spectral representation
\begin{equation} \label{} \label{LSFrepresentation}
X_{t,T} = \mu\big(\frac{t}{T}\big) + \frac{1}{\sqrt{2\pi}} \int^{\pi}_{-\pi} \exp (i \lambda
t) \,A_{t,T} (\lambda) \,d \xi (\lambda) \qquad \mbox{where} \qquad A_{t,T}(\lambda) \approx A(\frac {t} {T},\lambda)
\end{equation}
with the time varying transfer function $A(\cdot,\lambda)$ was used. Both representations are
basically equivalent -- see the derivation of (\ref{tvSpectralRepresentation}). In the results
presented below we will always use the formulation ``Under suitable regularity conditions \ldots''
and refer the reader to the original paper. We conjecture however that all results can be reproved
under Assumption~\ref{DefLocStat}. We emphasize that this is not an easy task since in most
situations it means to transfer the proof from the frequency to the time domain. In that case it
would be worthwhile to require only martingale differences $\varepsilon_{t}$ since also some
nonlinear processes admit such a representation.

Let
\begin{equation} \label{DefinitionBoundedVariation}
V(g) = \sup \Big\{ \sum^{m}_{k=1} \big| g(x_{k})-g(x_{k-1})\big| : 0
\leq x_{o} < \ldots< x_{m} \leq 1, \, m \in \mathbf N \Big\}
\end{equation}
be the total variation of $g$.

\begin{assumption} \label{DefLocStat} The sequence of stochastic processes
$X_{t,T}$ has a representation
\begin {equation} \label{TheLSProcess}
X_{t,T} = \mu\big(\frac{t}{T}\big) + \sum^{\infty}_{j = -\infty} a_{t,T} (j) \, \varepsilon_{t-j}
\end{equation}
where $\mu$ is of bounded variation and the $\varepsilon_{t}$ are iid with $E \varepsilon_{t} = 0$,
$E \varepsilon_{s} \varepsilon_{t} = 0$ for $s \not= t$, $E \varepsilon_{t}^{2} = 1$. Let
\begin{equation*}
\ell (j) :=
\begin{cases}
1, & |j| \leq 1 \\
|j| \log^{1+\kappa} |j|, &|j| > 1
\end{cases}
\end{equation*}
for some $\kappa > 0$ and
\begin {equation} \label{2}
\sup_{t} \big| a_{t,T} (j)\big| \leq \frac{K}{\ell(j)} \quad \mbox{(with $K$ indep. of $T$).}
\end{equation}
Furthermore we assume that there exist functions $a(\cdot, j) : (0,1] \rightarrow \mathbf R $ with
\begin{equation} \label{2a}
\sup_{u} \big|a(u,j)\big| \leq \frac{K}{\ell(j)},\quad
\end{equation}
\begin{equation} \label{3}
\sup_{j} \sum^{T}_{t=1} \Big|a_{t,T} (j) - a(\frac{t}{T}, j) \Big| \leq K,
\end{equation}
\begin{equation} \label{4}
V\big(a (\cdot, j)\big) \leq \frac{K}{\ell(j)}.\quad
\end{equation}
\end{assumption}

The above assumptions are weak in the sense that only bounded variation is required for the
coefficient functions. In particular for local results stronger smoothness assumptions have to be
imposed -- for example in addition for some $i$
\begin{equation} \label{2c}
\sup_{u}\Big|\frac{ \partial^{i} \mu(u)}{ \partial u^{i}}\Big| \leq K, \qquad \qquad
\end{equation}
\begin{equation} \label{2b}
\qquad \qquad \sup_{u} \Big|\frac{ \partial^{i} a(u,j)}{ \partial u^{i}}\Big| \leq
\frac{K}{\ell (j)} \quad \textrm{ for } j=0,1,\ldots
\end{equation}
and instead of (\ref{3}) the stronger assumption
\begin{equation} \label{3a}
\quad  \sup_{t,T} \Big|a_{t,T} (j) - a (\frac{t}{T}, j)\Big| \leq
\frac{K}{T\,\ell (j)}.
\end{equation}
The construction with $a_{t,T} (j)$ and $a(\frac{t}{T}, j)$ looks complicated at first glance. The
function $a(\cdot, j)$ is needed for rescaling and to impose necessary smoothness conditions while
the additional use of $a_{t,T} (j)$ makes the class rich enough to cover interesting cases such as
tvAR-models (the reason for this in the AR(1)-case can be understood from (\ref{Derivative_tvAR0})). Cardinali and Nason (2010) created the term
\textit{close pair} for $\big(a(\frac{t}{T}, j),a_{t,T} (j)\big)$. Usually, additional moment conditions on $\varepsilon_{t}$ are required.

It is straightforward to construct the stationary approximation and the derivative processes. We have
\begin{equation*} \label{}
\tilde{X}_{t} (u) := \mu(u) + \sum^{\infty}_{j = -\infty} a(u,j)
\, \varepsilon_{t-j}
\end{equation*}
and
\begin{equation*} \label{}
\frac{\partial^{i}\tilX_{t}(u)}{\partial u^{i}} = \frac{\partial^{i} \mu(u)}{\partial u^{i}} + \sum^{\infty}_{j = -\infty} \frac{\partial^{i} a(u,j)}{\partial u^{i}}
\, \varepsilon_{t-j}
\end{equation*}
and it is easy to prove (\ref{GeneralDefinitionLSP}) and more general the expansion (\ref{ExpansionDerivativeProcesses}). We define the time varying spectral
density by
\begin{equation} \label{tvSpectralDensity}
f(u,\lambda) := \frac{1}{2\pi}\big|A(u,\lambda)\big|^{2}
\end{equation}
where
\begin {equation} \label{tvTransferFunction}
A(u,\lambda) := \sum^{\infty}_{j = -\infty} a(u,j) \, \exp (-i \lambda j),
\end {equation}
and the time varying covariance of lag $k$ at rescaled time $u$ by
\begin{equation} \label{tvCovariance}
c(u,k) : =  \int^{\pi}_{-\pi} f(u, \lambda) \, \exp (i \lambda k) \, d
\lambda = \sum^{\infty}_{j=-\infty} a (u, k + j) \, a (u, j).
\end{equation}
$f(u,\lambda)$ and $c(u,k)$ are the spectral density and the covariance function of the stationary
approximation $\tilde{X}_{t}(u)$. Under Assumption~\ref{DefLocStat} and (\ref{3a}) it can be shown
that
\begin{equation} \label{tvCovarianceProperty}
\textrm{cov} \big(X_{[uT],T}, X_{[uT] + k,T}\big) = c(u,k) + O(T^{-1})
\end{equation}
uniformly in $u$ and $k$ -- therefore we  call $c(u,k)$ also the time varying covariance of the
processes $X_{t,T}$. In Theorem~\ref{UniquenessSD} we show that $f(u,\lambda)$ is the uniquely
defined time varying spectral density of $X_{t,T}$.

\begin{example} \label{ExampleLocStat}
(i) A simple example of a process $X_{t,T}$ which fulfills the above assumptions is $X_{t,T} = \mu
(\frac{t}{T}) + \phi (\frac{t}{T} ) Y_{t}$ where $Y_{t} = \Sigma_{j}\, a(j) \, \varepsilon_{t-j}$
is stationary with $ |a(j)| \leq K/\ell (j)$ and $\mu$ and $\phi$ are of bounded variation. If
$Y_{t}$ is an $AR(2)$-process with complex roots close to the unit circle then $Y_{t}$ shows a
periodic behavior and $\phi (\cdot)$ may be regarded as a time varying amplitude modulating
function of the process $X_{t,T}$. $\phi(\cdot)$ may either be
parametric or nonparametric.\\[6pt]
(ii) The tvARMA(p,q) process
\begin{equation} \label{9}
\sum^{p}_{j=0} \alpha_{j} (\frac{t}{T}) \, X_{t-j,T} = \sum^{q}_{k=0}
\beta_{k} (\frac{t}{T}) \, \sigma (\frac{t-k}{T}) \, \varepsilon_{t-k}
\end{equation}
where $\varepsilon_{t}$ are iid with $E \varepsilon_{t} = 0$ and $E \varepsilon_{t}^2 < \infty$ and
all $\alpha_{j} (\cdot)$, $\beta_{k}(\cdot)$ and $\sigma^2(\cdot)$ are of bounded variation with $\alpha_0(\cdot)\equiv\beta_0(\cdot) \equiv 1$ and
$\sum^{p}_{j=0} \alpha_{j} (u) z^{j} \not= 0$ for all $u$ and all $|z| \leq 1 + \delta$ for some
$\delta >0$, fulfills Assumption~\ref{DefLocStat}. If the parameters are differentiable with
bounded derivatives then also (\ref{2c})-(\ref{3a}) are fulfilled (for i=1). The time varying
spectral density is
\begin{equation} \label{SpectralDensityARMA}
f(u,\lambda) = \frac{\sigma^{2} (u)}{2\pi} \frac{|\sum^{q}_{k=0}
\beta_{k} (u) \exp (i \lambda k)|^{2}}{| \sum^{p}_{j=0} \alpha_{j}
(u) \exp (i \lambda j) |^{2}}.
\end{equation}
This is proved in Dahlhaus and Polonik (2006). $\alpha_j (\cdot)$ and $\beta_k (\cdot)$ may either
be parametric or nonparametric. \hfill $\Box$
\end{example}

The time varying MA($\infty$)-representation (\ref{TheLSProcess}) can easily be transformed into a time
varying spectral representation as used e.g. in Dahlhaus (1997, 2000). If the $\varepsilon_{t}$ are assumed to be stationary then there exists a Cram\'er
representation (cf. Brillinger (1981))
\begin{equation} \label{11}
\varepsilon_{t} = \frac{1}{\sqrt{2\pi}} \int^{\pi}_{-\pi} \exp (i
\lambda t) \, d \xi (\lambda)
\end{equation}
where $\xi (\lambda)$ is a process with mean $0$ and orthonormal increments. Let
\begin{equation}\label{12}
A_{t,T} (\lambda) := \sum^{\infty}_{j= -\infty} a_{t,T} (j) \,
\exp(-i \lambda j).
\end{equation}
Then
\begin{equation} \label{tvSpectralRepresentation}
X_{t,T} = \frac{1}{\sqrt{2\pi}} \int^{\pi}_{-\pi} \exp (i \lambda
t) \,A_{t,T} (\lambda) \,d \xi (\lambda).
\end{equation}
(\ref{3a}) now implies
\begin{equation} \label{tvSpectralRepresentationCond1}
\sup_{t,\lambda} \big|A_{t,T} (\lambda) - A(\frac{t}{T}, \lambda)\big|
\leq K T^{-1}
\end{equation}
which was assumed in the above cited papers. Conversely, if we start with (\ref{tvSpectralRepresentation}) and (\ref{tvSpectralRepresentationCond1}) then we can conclude from adequate smoothness conditions on $A(u,\lambda)$ to the conditions of
Assumption~\ref{DefLocStat}.

We now state a uniqueness property of our spectral representation. The Wigner-Ville spectrum for
fixed T (cf. Martin and Flandrin (1985)) is
\begin{equation*} \label{}
f_T(u,\lambda):=\frac{1}{2\pi}\sum_{s= -\infty}^{\infty} \cov(X_{[uT-s/2],T},
X_{[uT+s/2],T})\,\exp(-i\lambda s)
\end{equation*}
with $X_{t,T}$ as in (\ref{TheLSProcess})  (either with the coefficient extended as
constants for $u \notin [0,1]$  or set to $0$). Below we prove that
$f_T(u,\lambda)$ tends in squared mean to $f(u,\lambda)$ as defined in (\ref{tvSpectralDensity}).
Therefore it is justified to call $f(u,\lambda)$ the time varying spectral density of the process.

\begin{theorem} \label{UniquenessSD}
If $X_{t,T}$ is locally stationary and fulfills Assumption~\ref{DefLocStat}  and (\ref{2b}) for all
$j$ then we have for all $u\in(0,1)$
\begin{equation*} \label{}
\int_{-\pi}^{\pi} \big|f_T(u,\lambda)-f(u,\lambda)\big|^2 d\lambda=o(1).
\end{equation*}
\end{theorem}
\proofbegin The result was proved in Dahlhaus (1996b) under a different set of conditions. It
is not very difficult to prove the result also under the present conditions. \proofend

As a consequence the time varying spectral density $f(u,\lambda)$ is uniquely defined. If in addition the process $X_{t,T}$ is non-Gaussian, then even $A(u,\lambda)$ and therefore also
the coefficients $a(u,j)$ are uniquely determined which may be proved similarly by considering
higher-order spectra. Since $\mu(t/T)$ is the mean of the process it is also uniquely determined. This is
remarkable since in the non-rescaled case time varying processes do not have a unique spectral density or a unique time varying spectral representation
(cf. Priestley (1981), Chapter 11.1; M\'{e}lard and Herteleer-de Schutter (1989)). $f(u,\lambda)$ from Theorem~\ref{UniquenessSD} has been called instantaneous spectrum (in particular for tvAR-process -- c.f. Kitagawa and Gersch (1985)). The above theorem gives a theoretical justification for this definition.

There is a huge benefit from having a unique time varying spectral density. We now give
an example for this. We derive the limit of the Kullback-Leibler information for Gaussian processes
and show that it depends on $f(u,\lambda)$. Replacing this by a spectral estimate will lead to a
quasi likelihood for parametric models similar to the Whittle likelihood for stationary processes. Without a unique spectral density such a construction were
not possible.

Consider the exact Gaussian maximum likelihood estimate
\begin{equation*} \label{GaussianMLE}
\hat{\eta}_T^{ML}:=\argmin_{\eta\in\Theta_\eta} {\mathcal L}_T^{E}(\eta)
\end{equation*}
where $\eta$ is a finite-dimensional parameter (as in (\ref{BlockWhittleEstimate})) and
\begin{equation} \label{ExactGaussianLikelihood}
{\mathcal L}_T^{E}(\eta) = \frac{1}{2}\log(2
\pi)+\frac{1}{2T}\log\det {\Sigma}_{\eta}+\frac{1}{2T}
  ({\bf X}-{\bf \mu}_{\eta})'{\Sigma}_{\eta}^{-1}({\bf X}-{\bf
\mu}_{\eta})
\end{equation}
with ${\bf X}=(X_{1,T},\ldots,X_{T,T})'$, ${\bf
\mu}_{\eta}=\big(\mu_{\eta}
(1/T),\ldots,\mu_{\eta}(T/T)\big)'$ and $\Sigma_{\eta}$ being the
covariance matrix of the model. Under certain regularity conditions
$\hat{\eta}_T^{ML}$ will converge to
\begin{equation} \label{ThetaNull}
\eta_0:=\argmin_{\eta\in\Theta_\eta} {\mathcal L}(\eta)
\end{equation}
where
\begin{equation*} \label{}
{\mathcal L}(\eta):=\lim_{T\rightarrow\infty} \mean \,
   {\mathcal L}_T^{E} (\eta).
\end{equation*}
If the model is correct, then typically ${\eta}_0$ is the true parameter value.
Otherwise it is some ``projection'' onto the parameter space. It is therefore important to
calculate $\mathcal {L}(\eta)$ which is equivalent to the calculation of the
Kullback-Leibler information divergence.
\begin{theorem} \label{LimitLikelihood} Let $X_{t,T}$ be a locally stationary process with true mean- and spectral density curves
$\mu(\cdot)$, $f(u,\lambda)$ and model curves $\mu_{\eta}(\cdot)$,
$f_{\eta}(u,\lambda)$ respectively. Under suitable regularity conditions we have
\begin{align*}
  {\mathcal L} (\eta) &= \lim_{T \rightarrow \infty} \, \mean \, {\mathcal L}_T^{E}(\eta)\\
&=\frac{1}{4\pi}\int_0^1\int_{-\pi}^{\pi} \Big\{\log4{\pi}^2 f_{\eta}(u,\lambda)+
      \frac{f(u,\lambda)}{f_{\eta}(u,\lambda)}\Big\} \, d \lambda \, du + \frac{1}{4 {\pi}}\int_0^1\frac{\big(\mu_{
\eta}(u)-\mu(u)\big)^2}{f_{\eta}(u,0)} \,du.
\end{align*}
\end{theorem}
\proofbegin See Dahlhaus (1996b), Theorem~3.4.

The Kullback-Leibler information divergence for stationary processes is obtained from this as a special case \big(cf. Parzen (1983)\big).

\begin{example} \label{ModelStationary}
Suppose that the model is stationary, i.e.
$f_{\eta}(\lambda):=f_{\eta}(u,\lambda)$ and
$m:=\mu_{\eta}(u)$ do not depend on $u$. Then
\begin{equation*} \label{}
{\mathcal L}(\eta)=\frac{1}{4 {\pi}} \int_{-\pi}^{\pi} \Big\{ \log 4
{\pi}^2 f_{\eta}(\lambda)
+\frac{\int_0^1 f(u,\lambda) \, du}{f_{\eta}(\lambda)}\Big\} \, d
\lambda+\frac{1}{4 {\pi}}\,
f_{\eta}(0)^{-1} \int_0^1 \! \big(m-\mu(u)\big)^2 du
\end{equation*}
i.e. $m_0=\int_0^1\mu(u)du$, and $f_{\eta_0}(\lambda)$ give the best approximation to
the time integrated true spectrum $\int_0^1 f(u,\lambda) \, du$. These are the values which are ``estimated'' by the MLE or a quasi-MLE if a stationary model is fitted to locally stationary data. \hfill $\Box$
\end{example}

Given the form of ${\mathcal L}(\eta)$ as in Theorem~\ref{LimitLikelihood} we can
now suggest a quasi-likelihood criterion
\begin{equation*} \label{QuasiLikelihood}
  {\mathcal L}_{T}^{QL} (\eta) =
\frac{1}{4\pi}\int_0^1\int_{-\pi}^{\pi} \Big\{\log4{\pi}^2 f_{\eta}(u,\lambda)+
      \frac{ \hat{f}(u,\lambda)}{f_{\eta}(u,\lambda)}\Big\} \, d \lambda \, du + \frac{1}{4 {\pi}}\int_0^1\frac{\big(\mu_{
\eta}(u)- \hat{\mu}(u)\big)^2}{f_{\eta}(u,0)} \,du
\end{equation*}
where $\hat{f}(u,\lambda)$ and $\hat{\mu}(u)$ are suitable nonparametric estimates of $f
(u,\lambda)$ and $\mu(u)$ respectively. The block Whittle likelihood ${\mathcal L}_T^{BW}(\eta)$ in (\ref{BlockWhittleLikelihood}) and the generalized Whittle likelihood
${\mathcal L}_{T}^{GW} (\eta)$ in (\ref{GeneralizedWhittleLikelihood}) are of this
form.

We now calculate the Fisher information matrix
\begin{equation*} \label{}
\Gamma:=\lim_{T\rightarrow\infty} T \,{\mean}_{\eta_0}
\, \big(\nabla {\mathcal L}_T^{E}(\eta_0)\big) \, \big(\nabla {\mathcal {L}}_T^{E}(\eta_0)\big)'
\end{equation*}
in order to study efficiency of parameter estimates (see also Theorem~\ref{AsymptoticPropertiesGWhittle}).

\begin{theorem} \label{LimitFisherInfo}
Let $X_{t,T}$ be a locally stationary process with correctly specified mean curve
$\mu_{\eta}(u)$ and time varying spectral density
$f_{\eta}(u,\lambda)$. Under suitable regularity conditions we have
\begin{equation*}
\Gamma = \frac{1}{4{\pi}}\int_0^1\int_{-\pi}^{\pi} \big(\nabla\log f_{\eta_0}\big)
\big(\nabla \log f_{\eta_0}\big)' \,
d\lambda \, du + \frac{1}{2{\pi}} \int_0^1 \big(\nabla \mu_{\eta_0} (u)\big)
\big(\nabla\mu_{\eta_0} (u)\big)'
f_{\eta_0}^{-1} (u,0) \, du.
\end{equation*}
\end{theorem}
\proofbegin See Dahlhaus (1996b), Theorem 3.6.

\medskip

We now briefly discuss how the time varying spectral density can be estimated. Following the
discussion in the last section we start with a classical ``stationary'' smoothed periodogram estimate
on a segment. Let $I_{T} (u,\lambda)$ be the tapered periodogram on a segment
of length $N$ about $u$ as defined in (\ref{PeriodogramOnSegment}). Even in the stationary case $I_{T} (u,\lambda) $ is not a consistent estimate of the spectrum and
we have to smooth it over neighboring frequencies. Let therefore
\begin{equation} \label{tvSDestimate1}
\hat{f}_T (u,\lambda) := \frac{1}{b_{f}} \int K_{f} \left( \frac{\lambda
- \mu}{b_{f}} \right) I_{T} (u,\mu) \, d\mu
\end{equation}
where $K_{f}$ is a symmetric kernel with $\int K_{f} (x) \, dx = 1$ and $b_{f}$ is the bandwidth in
frequency direction. Theorem~\ref{} below shows that the estimate is implicitly also a kernel
estimate in time direction with kernel
\begin{equation}
K_{t} (x) : = \left\{ \int^{1}_{0} h(x)^{2} dx \right\}^{-1} h(x+
1/2)^{2}, \qquad x \in [-1/2, 1/2]
\label{2.1}
\end{equation}
and bandwidth $b_{t} : = N/T$, that is the estimate behaves like a kernel estimates with two
convolution kernels in frequency and time direction. We mention that an asymptotically equivalent
estimate is the kernel estimate
\begin {equation} \label{tvSDestimate2}
\tilde{f}_{T}(u,\lambda)  :=  \frac{2 \pi}{T^{2}} \sum^{T}_{t=1} \sum^{T}_{j=1} \int
\frac{1}{b_{t}} K_{t} \Big(\frac{u-t/T}{b_{t}}\Big) \;\frac{1}{b_{f}} K_{f}
\Big(\frac{\lambda-\lambda_j}{b_{f}}\Big) \; J_{T}
 \Big(\frac{t}{T}, \lambda_j \Big)
\end{equation}
with the pre-periodogram  $J_{T} (u, \lambda)$ as defined in (\ref{Preperiodogram}). One may also
replace the integral in frequency direction in (\ref{tvSDestimate1}) by a sum over the Fourier
frequencies.

\begin{theorem}
Let $X_{t,T}$ be a locally stationary process with $\mu (\cdot) \equiv 0$. Under suitable
regularity conditions we have
\begin{flalign*} \label{}
(i)\;\;    & \quad \mean I_{T}(u,\lambda) = f(u,\lambda) + \frac{1}{2} \,b_{t}^{2}
\int^{1/2}_{-1/2} x^{2} K_{t} (x) \, dx \, \frac{\partial^{2}}{\partial u^{2}} f(u,\lambda)  + o
(b_{t}^{2}) + \, O
\Big( \frac{\log (b_{t} T)}{b_{t} T}\Big);&\\[6pt]
(ii)\; & \quad \mean \hat{f}_T (u,\lambda) =  f(u,\lambda) + \frac{1}{2} \, b_{t}^{2}
\int^{1/2}_{-1/2}
x^{2} K_{t}(x) \, dx \, \frac{\partial^{2}}{\partial u^{2}} f (u,\lambda)&\\[6pt]
& \hspace*{3.7cm} + \frac{1}{2} \, b_{f}^{2} \int^{1/2}_{-1/2} x^{2} K_{f} (x) \,dx
\frac{\partial^{2}}{\partial \lambda^{2}} f(u,\lambda) + o \Big( b_{t}^{2} + b_{f}^{2} + \frac{\log
(b_{t} T)}{b_{t} T}\Big);&\\[6pt]
(iii) & \quad \var \big(\hat{f}_T(u,\lambda)\big) = \big(b_{t} b_{f} T\big)^{-1} \, 2 \pi f(u,\lambda)^{2}
\int^{1/2}_{-1/2} \!K_{t} (x)^{2} dx \int^{1/2}_{-1/2} \!K_{f} (x)^{2} dx \; \big(1 +
\delta_{\lambda 0} \big). &
\end{flalign*}
\end{theorem}
\proofbegin A sketch of the proof can be found in Dahlhaus (1996c), Theorem 2.2.

\medskip

Note, that the first bias term of $\hat{f}$ is due to nonstationarity while the second is due to
the variation of the spectrum in frequency direction.

As in Remark~\ref{MinimalMSE}  one may now minimize the relative mean squared error
RMSE$(\hat{f}) := E\big(\hat{f}(u,\lambda) / f (u,\lambda) - 1\big)^{2}$ with respect to $b_{f}$,
$b_{t}$ (i.e. $N$), $K_{f}$ and $K_{t}$ (i.e. the data taper $h$). This has been done in Dahlhaus
(1996c), Theorem 2.3. The result says that with
\begin{equation*} \label{}
\Delta_{u} := \frac{\partial^{2}}{\partial u^{2}} f(u,\lambda)
\Big/f(u,\lambda) \quad \mbox{and} \quad \Delta_{\lambda} := \frac{\partial^{2}}{\partial \lambda^{2}}
f(u,\lambda) \Big/f(u,\lambda)
\end{equation*}
the optimal RMSE is obtained with
\begin{equation*} \label{}
b_{t}^{\mbox{opt}} = T^{-1/6} (576 \pi)^{1/6} \left(
\frac{\Delta_{\lambda}}{\Delta^{5}_{u}} \right)^{1/12}, \qquad b_{f}^{\mbox{opt}} = T^{-1/6} (576 \pi)^{1/6} \left( \frac{\Delta_{u}}{
\Delta_{\lambda}^{5}} \right)^{1/12}
\end{equation*}
and optimal kernels $K_{t}^{\mbox{opt}} (x) = K_{f}^{\mbox{opt}} (x) = 6 \,\big(1/4 - x^{2}\big)$
with optimal rate $T^{-2/3}$.

The relations $b_{t} = N/T$ and (\ref{2.1}) immediately lead to the optimal segment length and the
optimal data taper $h$. The result of Theorem~\ref{}  is quite reasonable: If the degree of
nonstationarity is small, then $\Delta_{u}$ is small and $b_{t}^{\mbox{opt}}$ gets large. If the
variation of $f$ is small in frequency direction, then $\Delta_{\lambda}$ is small and
$b_{t}^{\mbox{opt}}$ gets smaller (more smoothing is put in frequency direction than in time
direction). This is another example, how the bias due to nonstationarity can be quantified with the
approach of local stationarity and balanced with another bias term and a variance term. Of course
the data-adaptive choice of the bandwidth parameters remains to be solved. Asymptotic normality of
the estimates can be derived from Theorem~\ref{CLTtimedependent} (cf. Dahlhaus (2009),
Example~4.2).

Rosen et.al. (2009) estimate the logarithm of the local spectrum by
using a Bayesian mixture of splines. They assume that the log spectrum on a partition of the data is a mixture of individual log spectra and use a mixture of smoothing splines with time varying mixing weights to estimate the evolutionary log spectrum. Guo et.al. (2003) use a smoothing spline ANOVA to estimate the time varying log spectrum.

\section{Gaussian likelihood theory for locally stationary processes}
\label{SectionLikelihood}

The basics of the likelihood theory for univariate stationary processes were laid by Whittle (1953,
1954). His work was much later taken up and continued by many others. Among the large number of
papers we mention the results of Dzhaparidze (1971) and Hannan (1973) for univariate time series,
Dunsmuir (1979) for multivariate time series and e.g. Hosoya and Taniguchi (1982) for misspecified
multivariate time series. A general overview over this likelihood theory and in particular
Whittle estimates for stationary models may be found in the monographs Dzhaparidze (1986) and
Taniguchi and Kakizawa (2000).

From a practical point of view the most famous outcome of this theory is the Whittle likelihood
\begin{equation} \label{stationaryWhittle}
\frac{1}{4\pi} \int^{\pi}_{-\pi} \biggl\{
\log 4\pi^{2} f_{\eta} (\lambda) + \frac{I_{T}(\lambda)}{f_{\eta}
(\lambda)} \biggr\} \, d \lambda
\end{equation}
as an approximation of the negative log Gaussian likelihood (\ref{ExactGaussianLikelihood}) where $I_{T}(\lambda)$
is the periodogram. This likelihood has been used also beyond the classical framework
 -- for example by Mikosch et al. (1995) for linear
processes where the innovations have heavy tailed distributions, by Fox and Taqqu (1986) for long
range dependent processes and by Robinson (1995) to construct semiparametric estimates for long
range dependent processes.

The outcome of this likelihood theory goes far beyond the construction of the Whittle likelihood.
Its technical core is the theory of Toeplitz matrices and in particular the approximation of the
inverse of a Toeplitz matrix by the Toeplitz matrix of the inverse function. It is essentially this
approximation which leads from the ordinary Gaussian likelihood to the Whittle likelihood. Beyond
that the theory can be used to derive the convergence of experiments for Gaussian stationary
processes in the H\'{a}jek-Le Cam sense, construct the properties of many tests and derive the
properties of the exact MLE and the Whittle estimate (cf. Dzhaparidze (1986); Taniguchi and Kakizawa (2000)).

For locally stationary processes it turns out that this likelihood theory can be generalized in a
nice way such that the classical likelihood theory for stationary processes arises as a
special case. Technically speaking this is achieved by a generalization of Toeplitz
matrices tailored especially for locally stationary processes (the matrix $U_{T} (\phi)$ defined in
(\ref{DefinitionMatrixU})).

Some results coming from this theory have already been stated in
Section~\ref{SectionSpectralDensity}, namely the limit of the Kullback-Leibler information
divergence in Theorem~\ref{LimitLikelihood} and the limit of the Fischer information in
Theorem~\ref{LimitFisherInfo}.  We now describe further results. We start with a decomposition of
the periodogram leading to a Whittle-type likelihood. We have
\begin{align}
I_{T}(\lambda) & =   \frac{1}{2\pi T} \Big| \sum^{T}_{r=1} X_{r} \exp (-i
\lambda r) \Big|^{2} \nonumber \\ \label{transformcovariances}
  & =  \frac{1}{2\pi} \sum^{T-1}_{k=-(T-1)} \Biggl( \frac{1}{T}
\sum^{T-|k|}_{t=1} X_{t} X_{t+|k|} \Biggr) \exp (-i \lambda k) \\
  & =  \frac{1}{T} \sum^{T}_{t=1} \frac{1}{2\pi} \mathop{{\sum}_{k}}_{1
\leq [t+0.5+k/2], [t+0.5 - k/2] \leq T} \!\!\!\!\!\!\!\!\! X_{[t+0.5+k/2],T}
X_{[t+0.5-k/2],T} \exp(-i\lambda k) \nonumber \\
  & = \frac{1}{T} \sum^{T}_{t=1} J_{T}(\frac{t}{T}, \lambda), \label{periodogramequality}
\end{align}
where the so-called pre-periodogram
\begin{equation} \label{Preperiodogram}
J_{T} (u, \lambda)\! := \!\frac{1}{2\pi} \mathop{{\sum}_{k}}_{1 \leq [uT+0.5 + k/2],
[uT+0.5-k/2] \leq T} \!\!\!\!\!\!\!\! X_{[uT+0.5 + k/2],T} X_{[uT+0.5 - k/2],T} \, \exp (-i \lambda k)
\end{equation}
may be regarded as a local version of the periodogram at time $t$. While the ordinary periodogram
$I_{T}(\lambda)$ is the Fourier transform of the covariance estimator of lag $k$ over the whole
segment \big(see (\ref{transformcovariances})\big) the pre-periodogram just uses the pair
$X_{[t+0.5+k/2]} X_{[t+0.5-k/2]}$ as a kind of ``local estimator'' of the covariance of lag $k$ at
time $t$ \big(note that $[t+0.5+k/2] - [t+0.5-k/2] =k$\big). The pre-periodogram was introduced by
Neumann and von Sachs (1997) as a starting point for a wavelet estimate of the time varying
spectral density. The above decomposition means that the periodogram is the average of the
pre-periodogram over time.

If we replace $I_{T} (\lambda)$ in (\ref{stationaryWhittle})  by the above average of the
pre-periodogram and afterwards replace the model spectral density $f_{\eta}
(\lambda)$ by the time varying spectral density $f_{\eta} (u,\lambda)$ of a
nonstationary model, we obtain the \underline{generalized Whittle
likelihood}
\begin{equation} \label{GeneralizedWhittleLikelihood}
{\mathcal L}_{T}^{GW} (\eta) :=  \frac{1}{T}
\sum^{T}_{t=1} \frac{1}{4\pi} \int^{\pi}_{-\pi} \bigg\{\!\log 4 \pi^{2} f_{\eta}
\Big(\frac{t}{T}, \lambda \Big)  + \frac{J_{T}
(\frac{t}{T}, \lambda)}{f_{\eta} (\frac{t}{T}, \lambda)} \bigg\} \, d \lambda.
\end{equation}
If the fitted model is stationary, i.e. $f_{\eta}(u, \lambda) =
f_{\eta} (\lambda)$ then (due to (\ref{periodogramequality})) the above likelihood
is identical to the Whittle likelihood and we obtain the classical Whittle estimator. Thus the
above likelihood is a true generalization of the Whittle likelihood to nonstationary processes. In
Theorem~\ref{Theorem:ApproximationLikelihood} we show that this likelihood is a very close
approximation to the Gaussian log likelihood for locally stationary processes. In particular (we
conjecture that) it is a better approximation than the block Whittle likelihood ${\mathcal L}_T^{BW}(\eta)$ from (\ref{BlockWhittleLikelihood}).

We now briefly state the asymptotic normality result in the parametric case.  An example is the
tvAR(2)-model with polynomial parameter curves from Section~\ref{DeepExample}.4.
Let
\begin{equation}
\hat{\eta}_{T}^{GW} := \argmin_{\eta\in \Theta_\eta} {\mathcal L}_{T}^{GW} (\eta)
\end{equation}
be the corresponding quasi likelihood estimate, $\hat{\eta}_{T}^{ML}$ be the
Gaussian MLE defined in (\ref{ExactGaussianLikelihood}), and $\eta_0$ as in
(\ref{ThetaNull}) i.e. the model may be misspecified.

\begin{theorem} \label{AsymptoticPropertiesGWhittle} Let $X_{t,T}$ be a locally stationary process. Under suitable regularity conditions
we have in the case $\mu (\cdot)=\mu_{\eta} (\cdot)= 0$
\begin{equation*} \label{}
\sqrt{T}\big(\hat{\eta}_{T}^{GW} - \eta_{0}\big) \dconv {\cal N}\big(0, \Gamma^{-1}
V \Gamma^{-1}\big) \quad \mbox{and} \quad  \sqrt{T}\big(\hat{\eta}_{T}^{ML} - \eta_{0}\big) \dconv
{\cal N}\big(0, \Gamma^{-1} V \Gamma^{-1}\big)
\end{equation*}
with
\begin{equation*}
\Gamma_{ij} =  \frac{1}{4\pi} \int^{1}_{0} \int^{\pi}_{-\pi}
\big(f - f_{\eta_{0}}\big) \, \nabla_{\! ij} f_{\eta_{0}}^{-1} \, d\lambda \, du
+ \frac{1}{4\pi} \int^{1}_{0} \int^{\pi}_{-\pi}
\big(\nabla_{\! i} \log
f_{\eta_{0}}\big) \big(\nabla_{\! j} \log f_{\eta_{0}}\big) \, d\lambda \,du
\end{equation*}
and
\begin{equation*}
V_{ij} = \frac{1}{4\pi} \int^{1}_{0} \int^{\pi}_{-\pi}  f \,
\big(\nabla_{\! i} f_{\eta}^{-1}\big) \, f \, \big(\nabla_{\!j} f_{\eta}^{-1}\big) \, d\lambda \,
du.
\end{equation*}

If the model is correctly specified then $V=\Gamma$ and $\Gamma$ is the same as in Theorem~\ref{LimitFisherInfo} -- that is both estimates are asymptotically Fisher-efficient. Even more the sequence of experiments is locally asymptotically normal (LAN) and both estimates are locally asymptotically minimax.
\end{theorem}
\proofbegin See Dahlhaus (2000), Theorem 3.1. LAN and LAM has been proved for the MLE in Dahlhaus (1996b), Theorem 4.1 and 4.2 -- these results together with the LAM-property of the generalized Whittle estimate also follow from the technical lemmas in Dahlhaus (2000) (cf. Remark~3.3 in that paper). \proofend

The corresponding result in the multivariate case and in the case $\mu (\cdot)\neq 0$ or
$\mu_{\eta} (\cdot) \neq 0$ can be found in Dahlhaus (2000), Theorem 3.1.

A deeper investigation of ${\mathcal L}_{T}^{GW} (\eta)$ reveals that it can be
derived from the Gaussian log-likelihood by an approximation of the inverse of the covariance
matrix. Let $\underline{X} = \big(X_{1,T}, \ldots, X_{T,T}\big)'$, $\underline{\mu} = \big(\mu
(\frac{1}{T}), \ldots, \mu (\frac{T}{T})\big)'$, and $\Sigma_{T}(A,B)$ and $U_{T}(\phi)$ be
$T\!\times\! T$ matrices with $(r,s)$-entry
\begin{equation}
\Sigma_{T}(A,B)_{r,s} = \frac {1} {2\pi} \int^{\pi}_{-\pi} \exp \big(i \lambda (r-s)\big) \,
A_{r,T} (\lambda) \,B_{s,T} (-\lambda) d \lambda
\end{equation}
and
\begin{equation} \label{DefinitionMatrixU}
U_{T} (\phi)_{r,s} = \int^{\pi}_{-\pi} \exp \big(i \lambda (r-s)\big) \; \phi
\Big(\frac{1}{T} \Big[ \frac{r+s}{2} \Big]^{\ast}\!, \lambda
\Big) d \lambda
\end{equation}
$(r,s = 1, \ldots T)$ where the functions $A_{r,T}(\lambda)$, $B_{r,T}(\lambda)$, $\phi
(u,\lambda)$ fulfill certain regularity conditions \big($A_{r,T}(\lambda)$ $B_{r,T}(\lambda)$ are
transfer functions or their derivatives as defined in (\ref{12})\big). $[x]^{*} = [x]$ denotes the
largest integer less or equal to $x$ (we have added the
* to discriminate the notation from brackets). Direct calculation shows that
\begin{equation}
{\mathcal L}_{T}^{GW} (\eta) =
\frac{1}{4\pi} \frac{1}{T} \sum_{t=1}^{T} \int^{\pi}_{-\pi} \log \big[
4 \pi^{2} f_{\eta} (\frac{t}{T}, \lambda)\big] \,d \lambda + \frac{1}{8 \pi^{2} T} (\underline{X} -
\underline{\mu}_{\eta})' \, U_{T} (f_{\eta}^{-1})\,(\underline{X} -
\underline{\mu}_{\eta}).
\end{equation}
Furthermore,  the exact Gaussian likelihood is
\begin{equation}
{\mathcal L}_{T}^{E} (\eta) := \frac{1}{2} \log (2\pi) + \frac{1}{2T}
\log \det \Sigma_{\eta}  + \frac{1}{2T} (\underline{X} -
\underline{\mu}_{\eta})' \Sigma_{\eta}^{-1} (\underline{X} -
\underline{\mu}_{\eta})
\end{equation}
where $ \Sigma_{\eta} = \Sigma_{T} (A_{\eta},
A_{\eta})$.

Proposition~\ref{Prop:ApproximationMatrices} below states that $U_{T} (\frac{1}{4\pi^{2}}
f_{\eta}^{-1})$ is an approximation of $\Sigma_{\eta}^{-1}$. Together
with the generalization of the Szeg\"{o} identity in Proposition~\ref{Prop:SzegoeIdentity} this implies
that ${\mathcal L}_{T}^{GW}$ is an approximation of ${\mathcal L}_{T}^{E}$ (see
Theorem~\ref{Theorem:ApproximationLikelihood}). If the model is stationary, then
$A_{\eta}$ is constant in time and $\Sigma_{\eta} = \Sigma_{T}
(A_{\eta}, A_{\eta})$ is the Toeplitz matrix of the spectral density
$f_{\eta} (\lambda) = \frac {1} {2\pi} |A_{\eta}|^{2}$ while $U_{T}
(\frac{1}{4\pi^{2}} f_{\eta}^{-1})$ is the Toeplitz matrix of $\frac{1}{4 \pi^{2}}
f_{\eta}^{-1}$. This is the classical matrix-approximation leading to the
Whittle likelihood (cf. Dzhaparidze, 1986).

\begin{proposition} \label{Prop:ApproximationMatrices}
Under suitable regularity conditions we have for each $\varepsilon > 0$ for the Euclidean norm
\begin{equation}
\frac{1}{T} \,\big\| \Sigma_{T} (A,
A)^{-1} - U_{T} \big(\{ 2
\pi A \bar{A}' \}^{-1}\big) \big\|_2^{2} = O (T^{-1 + \varepsilon})
\end{equation}
and
\[
\frac{1}{T} \, \big\| U_{T} (\phi)^{-1} - U_{T} \big(\{ 4 \pi^{2} \phi \}^{-1}\big)
\big\|_2^{2} = O(T^{-1+\varepsilon}).
\]
\end{proposition}
\proofbegin See Dahlhaus (2000), Proposition~2.4.

\medskip

By using the above approximation it is possible to prove the following generalization of the Szeg\"{o}
identity (cf. Grenander and Szeg\"{o} (1958), Section 5.2) to  locally stationary processes.

\begin{proposition} \label{Prop:SzegoeIdentity}
Under suitable regularity conditions we have with $f(u,\lambda) = \frac {1} {2\pi} |A (u,\lambda)|^{2}$
for each $\varepsilon >0$
\[
\frac{1}{T} \log \det \Sigma_{T} (A,A) = \frac{1}{2\pi} \int^{1}_{0}
\int^{\pi}_{-\pi} \log \big[2 \pi f(u,\lambda)\big] \,d \lambda \,du +
O(T^{-1+\varepsilon}).
\]
If $A = A_{\eta}$ depends on a parameter $\eta$ then the $O(T^{-1+
\varepsilon})$ term is uniform in $\eta$.
\end{proposition}
\proofbegin See Dahlhaus (2000), Proposition~2.5.

\medskip

In certain situations the right hand side can be written in the form $\int^{1}_{0} \log
\big(2 \pi \sigma^{2}(u) \big)\,du$ where $\sigma^{2}(u)$ is the one step prediction error at time $u$.

The mathematical core of the above results consists of the derivation of properties of products of matrices
$\Sigma_{T} (A,B)$, $\Sigma_{T}(A,A)^{-1}$ and $U_{T} (\phi)$. These properties are derived in
Dahlhaus (2000) in Lemma A.1, A.5, A.7 and A.8. These results are generalizations of corresponding
results in the stationary case proved by several authors before.

We now state the properties of the different likelihoods.

\begin{theorem}  \label{Theorem:ApproximationLikelihood}
Under suitable regularity conditions we have for $k = 0,1,2$

 (i)
\[ \sup_{\eta \in \Theta_\eta} \big| \nabla^{k} \big\{ {\mathcal L}_{T}^{GW} (\eta) - {\mathcal L}_{T}^{E} (\eta) \big\} \big| \pconv 0,
\]

(ii)
\[ \sup_{\theta \in \Theta_\eta} \big|\nabla^{k} \big\{ {\mathcal L}_{T}^{GW} (\eta) -
{\mathcal L} (\eta) \big\} \big| \pconv 0,
\]

(iii)
\[
\sup_{ \eta \in \Theta_\eta} \big| \nabla^{k} \big\{ {\mathcal L}_{T}^{E}
(\eta) - {\mathcal L} (\eta) \big\} \big| \pconv 0.
\]
\end{theorem}

\proofbegin See Dahlhaus (2000), Theorem 3.1.

\medskip

Under stronger assumptions one may also conclude that
$\hat{\eta}_{T}^{GW}\!\!-\hat{\eta}_{T}^{ML} = O_{p}(T^{-1+\varepsilon})$
which means that $\hat{\eta}_{T}^{GW}$ is a close approximation of the MLE.  A
sketch of the proof is given in Dahlhaus (2000), Remark 3.4.

\begin{remark} \label{}
It is interesting to compare the generalized Whittle estimate $\hat{\eta}_{T}^{GW}$
and its underlying approximation $U_{T}(\frac{1}{4\pi^{2}} f_{\eta}^{-1})$ of
$\Sigma^{-1}_{\eta}$ with the block Whittle estimate
$\hat{\eta}_{T}^{BW}$ defined in (\ref{BlockWhittleLikelihood}). There some
overlapping block Toeplitz matrices are used as an approximation which we regard as worse. A
similar result as in Proposition~\ref{Prop:ApproximationMatrices} has been proved in Lemma 4.7 of
Dahlhaus (1996a) for this approximation. We conjecture that also a similar result as in
Theorem~\ref{Theorem:ApproximationLikelihood} with ${\mathcal L}_{T}^{BW} (\eta)$
can be proved and even more that
$\hat{\eta}_{T}^{BW}\!\!-\hat{\eta}_{T}^{ML} = O_{p}\big(\frac {N}
{T^{1-\varepsilon}}+ \frac {1} {N}\big)$ (this is more a vague guess than a solid conjecture) which
means that the latter approximation and presumably also the estimate
$\hat{\eta}_{T}^{BW}$ are worse. It would be interesting to have more rigorous
results and a careful simulation study with a comparison of both estimates. \hfill $\Box$
\end{remark}
We now remember the generalized Whittle likelihood from (\ref{GeneralizedWhittleLikelihood}) which was
\begin{equation*} \label{*}
{\mathcal L}_{T}^{GW} (\eta) =  \frac{1}{T}
\sum^{T}_{t=1} \frac{1}{4\pi} \int^{\pi}_{-\pi} \bigg\{\!\log 4 \pi^{2} f_{\eta}
\Big(\frac{t}{T}, \lambda \Big)  + \frac{J_{T}
(\frac{t}{T}, \lambda)}{f_{\eta} (\frac{t}{T}, \lambda)} \bigg\} \, d \lambda.
\end{equation*}
Contrary to the true Gaussian likelihood this is a sum over time and the summands can be
interpreted as a local log likelihood at time point $t$. We therefore define
\begin {equation} \label{LocalWhittleLikelihoodAtTimet}
\ell^{\,*}_{t,T}\! (\boldsymbol{\theta} ) := \frac{1}{4\pi} \int^{\pi}_{-\pi} \Big\{
\log 4{\pi}^2 {f_{\boldsymbol{\theta} }(\lambda)} +  \frac {J_{T} (\frac{t}{T}, \lambda)}
{f_{\boldsymbol{\theta}}(\lambda)}  \Big\} \; d\lambda.
\end{equation}
(to avoid confusion we mention that we use the notation $\eta$ for a finite dimensional parameter which determines the whole curve, that is $\boldsymbol{\theta}(\cdot)=\boldsymbol{\theta}_\eta(\cdot)$ and $f_{\eta} (u, \lambda) = f_{\boldsymbol{\theta}_\eta (u)} (\lambda)$). We now can construct all nonparametric estimates (\ref{LocalKernelEstimate})--(\ref{FullConditionalLikelihood}) with $\ell_{t,T}\! \big(\boldsymbol{\theta} \big)$
replaced by $\ell^{\,*}_{t,T}\! \big(\boldsymbol{\theta} \big)$ leading in each of
the $5$ cases to an alternative local quasi likelihood estimate.

The parametric estimator (\ref{FullConditionalLikelihood}) with this
local likelihood is the estimate $\hat{\eta}_{T}^{GW}$ from above.
The orthogonal series estimator (\ref{OrthogonalSeriesEstimate}) with $\ell^{\,*}_{t,T}\! \big(\boldsymbol{\theta} \big)$ has been investigated for a truncated wavelet series
expansion together with nonlinear thresholding in Dahlhaus and Neumann (2001). The method
is fully automatic and adapts to different smoothness classes. It is shown that the usual rates of
convergence in Besov classes are attained up to a logarithmic factor. The nonparametric estimator
(\ref{GlobalNonparametricEstimate}) with $\ell^{\,*}_{t,T}\! \big(\boldsymbol{\theta} \big)$ is studied in Dahlhaus and Polonik (2006). Rates of convergence, depending on the
metric entropy of the function space, are derived. This includes in particular maximum likelihood
estimates derived under shape restriction. The main tool for deriving these results is the so
called empirical spectral processes discussed in the next section. The kernel estimator (\ref{LocalKernelEstimate}) with $\ell^{\,*}_{t,T}\! \big(\boldsymbol{\theta} \big)$ has been investigated in Dahlhaus (2009), Example 3.6. Uniform convergence has been proved in Dahlhaus and Polonik (2009), Section~4 (see also Example~\ref{ExampleLocalWhittle} and Theorem~\ref{uniform-convergence-Whittle} below). The local polynomial fit
(\ref{LocalPolynomialFit}) has not been investigated yet in combination with this likelihood.

The whole topic needs a more careful investigation -- both theoretically and from a practical
point including simulations and data-examples.

\section{Empirical spectral processes} \label{SectionEmpiricalSpectralProcesses}

We now emphasize the relevance of the empirical spectral process for linear locally
stationary time series. The theory of empirical processes does not only play a major role in
proving theoretical results for statistical methods but also provides a deeper understanding of
many techniques and the arising problems. The theory was first developed for stationary processes
(c.f. Dahlhaus~(1988), Mikosch and Norvaisa~(1997), Fay and Soulier~(2001)) and then extended to
locally stationary processes in Dahlhaus and Polonik (2006,2009) and Dahlhaus (2009). The empirical
spectral process is indexed by classes of functions. Basic results that later lead to several
statistical applications are a functional central limit theorem, a maximal exponential inequality
and a Glivenko-Cantelli type convergence result. All results use conditions based on the metric
entropy of the index class. Many results stated earlier in this article have been proved by using
these techniques.

The empirical spectral process is defined by
\begin {equation*}\label{}
E_{T} (\phi) := \sqrt{T} \, \Big(F_{T} (\phi) - F(\phi)\Big)
\end{equation*}
where
\begin {equation}\label{A6}
F (\phi) := \int^{1}_{0} \int^{\pi}_{-\pi} \phi (u,\lambda) \, f(u,\lambda) \, d \lambda \, du
\end{equation}
is the generalized spectral measure and
\begin{align}\label{A7}
F_{T} (\phi) & :=  \frac{1}{T} \sum^{T}_{t=1} \int^{\pi}_{-\pi}
\phi(\frac{t}{T}, \lambda) \, J_{T} (\frac{t}{T}, \lambda)\, d
\lambda
\end{align}
the empirical spectral measure with the pre-periodogram as defined in (\ref{Preperiodogram}).

We first give an overview of statistics that can be written in the form $F_{T} (\phi)$ - several of
them have already been discussed earlier in this article ($K_T$ always denotes a kernel function).
\begin{tabbing} \label{}
\hspace*{0.2cm} \= \hspace{0.6cm} \= \hspace{5.8cm} \=
\hspace{4.7cm} \= \hspace{3cm}\kill
\> 1. \> $\phi (u,\lambda) = K_{T}(u_{0}\! - \! u) \, \cos (\lambda
k)$
\> local covariance estimator \> (\ref{KernelCovarianceEstimate}) a.s.; Remark~\ref{RemarkRelatedEstimates}\\[6pt]
\> 2. \> $\phi (u,\lambda) = K_{T}(u_{0}\! - \! u) \,
K_{T}(\lambda_{0} \! - \! \lambda)$
\>  spectral density estimator \> (\ref{tvSDestimate2}) a.s.; Remark~\ref{RemarkRelatedEstimates}\\[6pt]
\> 3. \> $\phi (u,\lambda) =  K_{T}(u_{0}\! - \! u) \, \nabla
f_{\boldsymbol{\theta}_0}(u,\lambda)^{-1}$
\> $\nabla \mathcal{L}_{T}^{GW}(u_0, \boldsymbol{\theta}_0), \boldsymbol{\theta}_0=\boldsymbol{\theta}_0 (u_0)$ \> Example~\ref{ExampleLocalWhittle}\\[6pt]
\> 4. \> $\phi (u,\lambda) \approx  K_{T}(u_{0}\! - \! u) \, \nabla
f_{\boldsymbol{\theta}_0}(u,\lambda)^{-1}$
\> local least squares \> Ex.~\ref{ExampleNonlinearModels}; Rem.~\ref{RemarkRelatedEstimates}\\[6pt]
\> 5. \> $\phi (u,\lambda) = \nabla f_{\eta_0}(u,\lambda)^{-1}$
\> param. Whittle estimator \> Example 3.7 in\\
\> \> \> \hspace*{3.10cm} Dahlhaus and Polonik (2009)\\[6pt]
\> 6. \> $\phi (u,\lambda) = \big(I_{[0,u_0]} (u) \!- \!u_0 \big)
I_{[0,\lambda_0]} (\lambda)$
\> testing stationarity \> Example~\ref{Testing Stationarity}\\[6pt]
\> 7. \> $\phi (u,\lambda) = \cos(\lambda k)$
\> stationary covariance \> Remark~\ref{stationary case}\\[6pt]
\> 8. \> $\phi (u,\lambda) = \nabla f_{\eta_0}(\lambda)^{-1}$
\> stat. Whittle estimator \> Remark~\ref{stationary case}\\[6pt]
\> 9. \> $\phi (u,\lambda) = K_{T}(\lambda_{0} \! - \! \lambda)$
\> stationary spectral density \> Remark~\ref{stationary case}
\end{tabbing}
Examples 1-4 and 9 are examples with index functions $\phi_{T}$ depending on $T$. More complex examples are nonparametric maximum likelihood
estimation under shape restrictions (Dahlhaus and Polonik, 2006), model selection with a sieve
estimator (Van Bellegem and Dahlhaus, 2006) and wavelet estimates (Dahlhaus and Neumann, 2001).
Moreover $F_T (\phi)$ occurs with local polynomial fits (Kim, 2001; Jentsch, 2006) and several
statistics suitable for goodness of fit testing. These applications are quite involved.

However, applications are limited to quadratic statistics, that is the empirical spectral measure is usually of no help in dealing with nonlinear models.  Furthermore, for linear processes the empirical process only applies without further modification to the (score function and the Hessian of the) likelihood ${\mathcal L}_{T}^{GW} (\eta)$ and its local variant ${\mathcal L}_{T}^{GW} (u,\boldsymbol{\theta})$ and the local Whittle likelihood ${\mathcal L}_{T}^{W} (u,\boldsymbol{\theta})$. It also applies to the exact likelihood ${\mathcal L}_{T}^{E} (\eta)$ after proving $\nabla  {\mathcal L}_{T}^{GW}
(\eta_0) - \nabla {\mathcal L}_{T}^{E} (\eta_0) = o_p (T^{-1/2})$ \big(see also Theorem~\ref{Theorem:ApproximationLikelihood} (i)\big) and the conditional likelihoods ${\mathcal L}_{T}^{C} (\eta)$ and ${\mathcal L}_{T}^{C} (u,\boldsymbol{\theta})$ in the tvAR-case (see Remark~\ref{RemarkRelatedEstimates} - in the general case this is not clear yet). For the block Whittle likelihood ${\mathcal L}_{T}^{BW} (\eta)$ it may also be applied after establishing $\nabla  {\mathcal L}_{T}^{GW}
(\eta_0) - \nabla {\mathcal L}_{T}^{BW} (\eta_0) = o_p (T^{-1/2})$. However, this is also not clear yet.

We first state a central limit theorem for $E_{T} (\phi)$ with index functions $\phi$ that do not
vary with $T$. We use the assumption of bounded variation in both components of $\phi (u,\lambda)$.
Besides the definition in (\ref{DefinitionBoundedVariation}) we need a definition in two
dimensions. Let
\begin {align*} \label{DefinitionBoundedVariation2dim}
V^{2}(\phi) = \sup \Big\{ \sum^{\ell,m}_{j,k=1} & | \phi (u_{j},\lambda_{k}) - \phi
(u_{j-1},\lambda_{k}) - \phi
(u_{j},\lambda_{k-1}) + \phi (u_{j-1},\lambda_{k-1}) |\, :\\
&  0 \leq u_{0} < \ldots< u_{\ell} \leq 1; \; - \pi \leq \lambda_{0} < \ldots< \lambda_{m} \leq
\pi; \; \ell, m \in \nnum \Big\}.
\end{align*}
\noindent For simplicity we set
\begin {equation*} \label{}
\|\phi\|_{\infty,V} := \sup_{u} V \big(\phi (u,\cdot)\big), \quad \|\phi\|_{V,\infty} :=
\sup_{\lambda} V \big( \phi (\cdot, \lambda)\big),
\end{equation*}
\begin {equation*} \label{}
\|\phi\|_{V,V} := V^{2}(\phi) \quad \textrm{and} \quad \|\phi\|_{\infty,\infty} := \sup_{u,\lambda}
|\phi(u,\lambda)|.
\end{equation*}

\begin{theorem} \label{clt} Suppose Assumption~\ref{DefLocStat}
holds and let $\phi_{1}, \ldots, \phi_{d}$ be functions with
 $\|\phi_{j}\|_{\infty,V}$,
$\|\phi_{j}\|_{V,\infty}$, $\|\phi_{j}\|_{V,V}$ and $\|\phi_{j}\|_{\infty,\infty}$ being finite
$(j=1,\ldots,d)$. Then
\[
\big(E_{T} (\phi_{j})\big)_{j=1,\ldots, d} \dconv \big(E
(\phi_{j})\big)_{j=1,\ldots,d}
\]
where $\big(E (\phi_{j})\big)_{j=1,\ldots, d}$ is a Gaussian random vector with mean $0$ and
\begin{align} \label{cov-structure}
\cov \big(E (\phi_{j}), E(\phi_{k}) \big) & = 2\pi \int^{1}_{0}
 \int^{\pi}_{-\pi} \phi_{j}
(u,\lambda) \, [ \phi_{k} (u,\lambda) + \phi_{k} (u,- \lambda) ] \,
f^2 (u,\lambda) \, d \lambda \, du & \;\\ \nonumber
& + \kappa_{4} \int^{1}_{0}  \Big( \int^{\pi}_{-\pi} \phi_{j} (u,
\lambda_{1}) f (u, \lambda_{1}) \, d \lambda_{1} \Big) \Big(
\int^{\pi}_{-\pi} \phi_{k} (u, \lambda_{2}) f (u, \lambda_{2}) \, d
\lambda_{2} \Big) \, du.
\end{align}
\end{theorem}
\proofbegin See Dahlhaus and Polonik (2009), Theorem 2.5.

\begin{remark} [Stationary processes/Model misspecification by stationary models] \label{stationary case}
The classical central limit theorem for the weighted periodogram in the stationary case can be
obtained as a corollary: If $\phi (u,\lambda)=\tilde{\phi}( \lambda)$ is time-invariant then
\begin {equation} \label{SpectralMeasureStatCase}
F_{T} (\phi) = \int^{\pi}_{-\pi} \tilde{\phi}( \lambda) \, \frac {1} {T} \sum^{T}_{t=1} \, J_{T}
(\frac{t}{T}, \lambda) \, d \lambda = \int^{\pi}_{-\pi} \tilde{\phi}( \lambda) \, I_{T}(\lambda) \,
d \lambda
\end{equation}
(see(\ref{periodogramequality})) that is $F_{T} (\phi)$ is the classical spectral measure in the
stationary case with the following applications:

\begin{romanlist}
\item $\phi (u,\lambda)=\tilde{\phi}( \lambda) = \cos \lambda k$ is the empirical covariance
    estimator of lag $k$;
\item $\phi (u,\lambda)=\tilde{\phi}( \lambda) = \frac {1} {4\pi} \nabla
    f_{\boldsymbol{\theta}}^{-1}(\lambda)$
    is the score function of the Whittle likelihood.
\end{romanlist}

\noindent Theorem~\ref{clt} gives the asymptotic distribution for these examples - both in the
stationary case and in the misspecified case where the true underlying process is only locally
stationary. If $\phi (u,\lambda)=\tilde{\phi}( \lambda)$ is a kernel we obtain an
estimate of the spectral density whose asymptotic distribution is a special case of
Theorem~\ref{CLTtimedependent} below (also in the misspecified case). \hfill $\Box$
\end{remark}

We now state a central limit theorem for $F_{T} (\phi_{T}) - F(\phi_{T})$ with index functions
$\phi_{T}$ depending on $T$. In addition we extend the hitherto definitions to tapered data
\begin{equation*} \label{}
X_{t,T}^{(h_T)} := h_T\Big(\frac{t}{T}\Big)\cdot X_{t,T}
\end{equation*}
where $h_T: (0,1] \to [0,\infty)$ is a data taper (with $h_{T} (\cdot) = I_{(0,1]} (\cdot)$ being
the nontapered case). The main reason for introducing data-tapers is to include segment estimates -
see the discussion below. As before the empirical spectral measure is defined by
\begin{align}\label{xA17}
F_{T} (\phi) & =  F_{T}^{(h_T)} (\phi) := \frac {1} {T}
\sum^{T}_{t=1} \int^{\pi}_{-\pi} \phi(\frac{t}{T}, \lambda) \,
J_{T}^{(h_{T})} (\frac{t}{T}, \lambda)\, d \lambda
\end{align}
now with the tapered pre-periodogram
\begin{equation} \label{xA18}
J^{\left(h_{T}\right)}_{T} \left(\frac{t}{T}, \lambda\right) =
\frac{1}{2\pi} \sum_{k: 1\leq [t+1/2 \pm k/2] \leq T}
X^{(h_{T})}_{[t+1/2 + k/2],T} \,X^{(h_{T})}_{[t+1/2 - k/2],T} \exp (-i
\lambda k)
\end{equation}
\big(we mention that in some cases a rescaling may be necessary for $J^{\left(h_{T}\right)}_{T}
\left(u , \lambda\right)$ to become a pre-estimate of $f(u,\lambda)$ - an obvious example for this
is $h_{T} (u) = (1/2) \, I_{(0,1]} (u)$\big).

$F(\phi)$ is the theoretical counterpart of $F_T (\phi)$
\begin {equation}\label{xA19}
F (\phi) = F^{(h_T)} (\phi) := \int^{1}_{0}  h_T^{2} (u) \int^{\pi}_{-\pi} \phi (u,\lambda)\,
f(u,\lambda) \, d \lambda \, du.
\end{equation}

Note that (\ref{periodogramequality}) also holds with a data-taper, that is
\begin{equation*} \label{}
\frac{1}{T} \sum^{T}_{t=1}
J_{T}^{(h_T)}\big(\frac{t}{T},
\lambda\big) = \frac {H_{2,T}} {T} \,I_{T}^{(h_T)}(\lambda)
\end{equation*}
with the tapered periodogram
\begin{equation} \label{tapered periodogram}
I_{T}^{(h_T)}(\lambda) :=  \frac{1}{2\pi H_{2,T}} \, \bigg|
\sum^{T}_{s=1} X_{s}^{(h_T)} \exp (-i \lambda s) \bigg|^{2} \quad \mbox{where} \quad H_{2,T} := \sum_{t=1}^{T} h_T\Big(\frac {t} {T}\Big)^{2}.
\end{equation}

An important special case is $h_{T}^{(u_0)} (\frac {t} {T}) := k\big(\frac {u_0-t/T} {b_T}\big)$
with bandwidth $b_T$ and $k$ having compact support on $[-\frac {1} {2},\frac {1} {2}]$. If $b_T :=
N/T$ then $I_{T}^{(h_T)}(\lambda)=I_{T} (u_0,\lambda)$ with $I_{T} (u_0,\lambda)$ as in
(\ref{PeriodogramOnSegment}). If in addition $\phi(u, \lambda) = \psi(\lambda)$ we obtain
\begin{equation*} \nonumber
F_T (\phi)  = \int^{\pi}_{-\pi} \psi(\lambda) \, \Big(\frac{1}{T} \sum^{T}_{t=1}
J_{T}^{(h_T)}\big(\frac{t}{T},
\lambda\big)\Big) \, d\lambda = \frac {H_{2,T}} {T} \,\int^{\pi}_{-\pi} \psi(\lambda) \,
I_{T}^{(h_T)}(\lambda) \, d\lambda.
\end{equation*}
For example for $\psi(\lambda):= \exp {i\lambda k}$ this is exactly  $\frac {H_{2,T}} {T} \, {\hat
c}_{T}(u_0,k)$ with the tapered covariance estimate from (\ref{TaperedCovarianceEstimate}). In this
case $\frac {H_{2,T}} {T}$ is proportional to $b_T$.

The last example suggests to use  $\frac {1} {H_{2,T}}$ instead of $\frac {1} {T}$ in (\ref{xA17})
as a norming constant. However, this is not always the right choice (as can be seen from case (ii)
in Remark~\ref{Typical applications}).

It turns out that in the above situation the rate of converge of the empirical spectral measure
becomes $\sqrt{T} / \rho_{2}^{(h_T)} (\phi_T)$ where
\begin{align*}\label{}
\rho_{2}^{(h_T)} (\phi) & := \Big(\int_{0}^{1} h_T^{4} (u) \,
\int_{-\pi}^{\pi} \phi (u,\lambda)^{2} \,d\lambda \; du\Big)^{1/2}.
\end{align*}
Therefore we can embed this case into the situation treated in the last section by studying the
convergence of
\begin{equation*}\label{}
E_{T}^{(h_T)} \Big( \frac {\phi_{T}} {\rho_{2}^{(h_T)} (\phi_T)}
\Big) = \frac {\sqrt{T}} {\rho_{2}^{(h_T)} (\phi_T)} \, \Big(F_{T}
(\phi_{T}) - F^{(h_T)}(\phi_{T})\Big).
\end{equation*}
Furthermore, let
\begin{align} \label{cov-structure2}
c_{E}^{(h_T)} (\phi_{j}, \phi_{k}) & := 2\pi \int^{1}_{0} \! h_T^{4}
(u) \int^{\pi}_{-\pi} \phi_{j} (u,\lambda) \, \big[ \phi_{k}
(u,\lambda) + \phi_{k} (u,- \lambda) \big] \, f^2 (u,\lambda) \, d
\lambda \, du
\\ \nonumber
& \; + \kappa_{4} \int^{1}_{0} \! h_T^{4} (u) \Big(
\int^{\pi}_{-\pi} \!\phi_{j} (u, \lambda_{1}) f (u, \lambda_{1})  d
\lambda_{1} \Big) \! \Big( \int^{\pi}_{-\pi} \!\phi_{k} (u,
\lambda_{2}) f (u, \lambda_{2})  d \lambda_{2} \Big)  du.
\end{align}

\medskip
\begin{theorem} \label{CLTtimedependent} Suppose that $X_{t,T}$ is a locally stationary process and suitable regularity
conditions hold. If the limit
\begin {equation} \label{ConvergenceMoments}
\Sigma_{j,k} := \lim_{T \rightarrow \infty} \, \frac {c_{E}^{(h_T)} (\phi_{Tj}, \phi_{Tk})}
{\rho_{2}^{(h_T)} (\phi_{Tj}) \, \rho_{2}^{(h_T)} (\phi_{Tk})}
\end{equation}
exists for all $j,k=1,\ldots,d$ then

\begin {equation} \label{CLT with kernels}
\Bigg(\frac {\sqrt{T}} {\rho_{2}^{(h_T)} (\phi_{Tj})} \, \Big(F_{T} (\phi_{Tj}) -
F^{(h_T)}(\phi_{Tj})\Big)\Bigg)_{\!j=1,\ldots, d} \dconv \;\mathcal{N}(0,\Sigma).
\end{equation}
\end{theorem}

\medskip

\begin{remark} [Bias] \label{TvBias}
In addition we have the bias term
\begin {equation*}\label{}
\frac {\sqrt{T}} {\rho_2^{(h_T)}(\phi_{T})} \, \Big(F^{(h_T)} (\phi_{T}) -
\lim_{T \rightarrow\infty} F^{(h_T)} (\phi_{T})\Big).
\end{equation*}
The magnitude of this bias depends on the smoothness of the time varying spectral density. In this
section we usually require conditions such that this bias is of lower order. This is different in
Section~\ref{NonlinearProcAndDerivative} where
the bias has explicitly been investigated.
\hfill $\Box$
\end{remark}

\begin{remark} [Typical applications] \label{Typical applications}
A typical application of this result is the case of kernel type local estimators which can be
constructed by using kernels, data-tapers or a combination of both:
\begin{tabbing} \label{}
\hspace*{0.7cm} \= \hspace{1cm} \= \hspace{6cm} \= \hspace{3cm}\kill
\> (i) \> $\phi_{T} (u,\lambda) = \frac{1}{b_T} K
(\frac{u_0-u}{b_T})
\;\psi (\lambda)$  \>  $h_{T} (\cdot) = I_{(0,1]} (\cdot)$\\[6pt]
\> (ii) \> $\phi_{T} (u,\lambda) = \frac{1}{b_T} K
(\frac{u_0-u}{b_T}) \;\psi (\lambda)$  \>  $h_{T} (u) =
I_{[u_0-b_T/2,u_0+b_T/2]} (u)$\\[6pt]
\> (iii) \> $\phi_{T} (u,\lambda) = \psi (\lambda)$  \>  $h_{T}
(\frac {t} {T}) = k\big(\frac {u_0-t/T} {b_T}\big)$
\end{tabbing}
where $K(\cdot)$ and $k(\cdot)$ are kernel functions and $b_T$ is the bandwidth. If $K(\cdot)=
k(\cdot)^{2}$ then the resulting estimates all have the same asymptotic properties - see below.
Dependent on the function $\psi (\lambda)$ this leads to different applications: If we set
$\psi(\lambda)=\cos(\lambda k)$ the estimate (iii) is the estimate ${\hat c}_{T}(u_0,k)$ from
(\ref{TaperedCovarianceEstimate}) and (i) is ``almost'' the estimate $\tilde{c}_{T} (u_0,k)$ from
(\ref{KernelCovarianceEstimate}) (for $k$ even it is exactly the same, for $k$ odd the difference
can be treated with the methods mentioned in Remark~\ref{}).

We now show how Theorem~\ref{CLTtimedependent} leads to the asymptotic distribution for these
estimates:

(i) If $K(\cdot)$ and $\psi (\cdot)$ are of bounded variation and $b_T \rightarrow 0$, $b_T T
\rightarrow \infty$ then the regularity conditions of Theorem~\ref{}  are fulfilled (see Dahlhaus
(2009), Remark 3.4). Furthermore,
\begin{equation} \label{statement3}
\rho_{2}^{(h_T)} (\phi_{T}) = \rho_{2} (\phi_{T}) = \Big( \frac {1}
{b_T} \int K^{2}(x) \, dx \int |\psi (\lambda)|^{2} \, d\lambda
\Big)^{1/2} \approx b_T^{-1/2}.
\end{equation}
For  $f(\cdot,\lambda)$ continuous at $u_0$ we have
\begin{align*} \label{cov-structure2}
c_{E}^{(h_T)} (\phi_{Tj}, \phi_{Tk}) &  \sim \frac {1} {b_T} \int K^{2}(x) \, dx \, \Big[ 2\pi \,
\int^{\pi}_{-\pi} \psi_{j} (\lambda) \, \big[ \psi_{k} (\lambda) +
\psi_{k} (-
\lambda) \big] \, f^2 (u_0,\lambda) \, d \lambda\\
\nonumber
& \;\; + \kappa_{4} \, \Big( \int^{\pi}_{-\pi} \!\psi_{j} (
\lambda_{1}) f (u_0, \lambda_{1}) \, d \lambda_{1} \Big) \! \Big(
\int^{\pi}_{-\pi} \!\psi_{k} (\lambda_{2}) f (u_0, \lambda_{2}) \, d
\lambda_{2} \Big) \Big] =: \frac {1} {b_T} \Gamma_{jk}
\end{align*}

that is (\ref{ConvergenceMoments}) is also fulfilled and we obtain from Theorem~\ref{CLTtimedependent}
\begin{equation} \label{CLTkernelestimate}
\sqrt{b_T T} \, \Big(F_{T} (\phi_{Tj}) -
F^{(h_T)}(\phi_{Tj})\Big)_{j=1,\ldots,d} \dconv
\;\mathcal{N}(0,\Gamma).
\end{equation}

(ii) The additional taper $h_{T} (u) = I_{[u_0-b_T/2,u_0+b_T/2]} (u)$ implies that we use only data
from the interval $[u_0-b_T/2,u_0+b_T/2]$. We obtain in this case
\begin{equation*}\label{}
\rho_{2}^{(h_T)} (\phi_T) =  \Big(\int_{0}^{1}  \frac{1}{b_T^{2}} K
\Big(\frac{u_0-u}{b_T}\Big)^{2} \, du \int_{-\pi}^{\pi} |\psi (\lambda)|^{2}
\, d\lambda \Big)^{1/2},
\end{equation*}
i.e. we have the same $\rho_{2}^{(h_T)} (\phi_T)$ as above. Furthermore, $c_{E}^{(h_T)} (\phi_{T},
\phi_{T})$ is the same. Thus we obtain the same asymptotic distribution and the same rate of
convergence.

(iii) If $K(\cdot)= k(\cdot)^{2}$ we obtain in this case
\begin{align*}\label{}
\frac {1} {b_T} \, \rho_{2}^{(h_T)} (\phi_T) & = \Big(\int_{0}^{1}
\frac{1}{b_T^{2}} K \Big(\frac{u_0-u}{b_T}\Big)^{2} \, du \int_{-\pi}^{\pi}
|\psi (\lambda)|^{2} \, d\lambda \Big)^{1/2}
\end{align*}
i.e. we obtain again the same expression. Furthermore, $\frac {1} {b_T^{2}} \, c_{E}^{(h_T)}
(\phi_{Tj}, \phi_{Tk})$ is the same as $c_{E}^{(h_T)} (\phi_{Tj}, \phi_{Tk})$ above. Thus we have
again the same asymptotic distribution and the same rate of convergence. \hfill $\Box$
\end{remark}

\begin{example} [Curve estimation by local quasi likelihood estimates]
\label{ExampleLocalWhittle} $\;$\\[3pt]
Local Whittle estimates on a segment where defined in (\ref{NonparametricWhittleEstimates}) and
discussed in Example~\ref{LocalWhittleEstimatesBiasVariance} (the bias was heuristically derived in
Example~\ref{BiasExpressionLocalLikelihood}). We now consider the presumably equivalent local quasi
likelihood estimate defined by %
\begin{equation} \label{GeneralizedKernelEstimate}
\hat{\boldsymbol{\theta}}_{T}^{GW}\!(u_0) := \argmin_{\boldsymbol{\theta} \in \Theta}
\mathcal{L}_{T}^{GW}(u_0, \boldsymbol{\theta})
\end{equation}
with
\begin {equation} \label{GeneralizedKernelLikelihood}
\mathcal{L}_{T}^{GW}(u_0, \boldsymbol{\theta}) := \frac{1}{4\pi} \frac{1}{T} \sum^{T}_{t=1}
\;\frac{1}{\,b_{T}}\, K \Big(\frac{u_0-t/T}{b_{T}}\Big) \; \int^{\pi}_{-\pi} \Big\{ \log 4{\pi}^2
f_{\boldsymbol{\theta}} (\lambda) + \frac {J_{T} (\frac{t}{T}, \lambda)} {f_{\boldsymbol{\theta}}
(\lambda)} \Big\} \, d\lambda.
\end {equation}
(this is a combination of (\ref{LocalKernelEstimate}) and (\ref{LocalWhittleLikelihoodAtTimet})).
The asymptotic normality of the estimate $\hat{\boldsymbol{\theta}}_{T}^{GW}\!(u_0)$ is derived
in Dahlhaus (2009), Example 3.6. Key steps in the proof are the fact that both the score function
and the Hessian matrix can be written in terms of the empirical spectral process leading to a
rather simple proof. For example
\begin{align} \label{score function Whittle}
\sqrt{b_T T} & \;\nabla_i \mathcal{L}_{T}\big(u_0,\boldsymbol{\theta}_{0}(u_0)\big)
= \sqrt{b_T T} \Big( F_{T} (\phi_{T,u_0,i}) - F (\phi_{T,u_0,i})\Big) + o_p(1)
\end{align}
where $\phi_{T,u_0,i} (v,\lambda) := \frac{1}{\,b_{T}} K (\frac{u_0-v}{b_{T}}) \; \frac {1} {4 \pi}
\nabla_i f_{\boldsymbol{\theta}}^{-1}(\lambda)_{|\boldsymbol{\theta}=\boldsymbol{\theta}_{0}(u_0)}$. Theorem~\ref{CLTtimedependent} then immediately gives the asymptotic normality of the score function and after some additional considerations also asymptotic normality of $\hat{\boldsymbol{\theta}}_{T}^{GW}\!(u_0)$. For details see Dahlhaus (2009), Example 3.6.

The above estimate corresponds to case (i) in Remark~\ref{Typical applications}. Case (iii) in Remark~\ref{Typical applications} leads instead to the tapered Whittle estimate $\hat{\boldsymbol{\theta}}_{T}^{W}(u_0)$ on the segment, since for $h_{T}^{(u_0)} (\frac {t} {T}) := k\big(\frac {u_0-t/T} {b_T}\big)$ we have $I_{T}^{(h_T)}(\lambda)=I_{T} (u_0,\lambda)$ with $I_{T} (u_0,\lambda)$ as in (\ref{PeriodogramOnSegment}). This estimate has the same asymptotic properties provided $k(\cdot)^{2}=K(\cdot)$. It's asymptotic properties can now also be derived by using Theorem~\ref{CLTtimedependent}.
\hfill $\Box$
\end{example}

\begin{remark} [Related estimates] \label{RemarkRelatedEstimates}%
Many estimates are only approximately of the form discussed above - \underline{for example} the sum statistic
\begin{equation} \label{A7sum}
F_{T}^{\Sigma} (\phi) :=  \frac{2 \pi}{T^{2}} \sum^{T}_{t=1}
\sum_{j=1}^{T} \phi \Big(\frac{t}{T}, \lambda_{j}\Big) \, J_{T}^{(h_T)}
\Big(\frac{t}{T}, \lambda_{j}\Big)
\end{equation}
where $\lambda_{j} = \frac {2 \pi j} {T}$ - or representations in terms of the
Fourier-coefficients. Important examples of related estimates are the spectral density estimate
(\ref{tvSDestimate2}), the covariance estimates (\ref{KernelCovarianceEstimate})
and (\ref{KernelCovarianceEstimate2}) and the score function of the local least squares
tvAR(p)-estimate from Example~\ref{ExampleNonlinearModels}. We mention that the central limit theorem in
Theorem~\ref{CLTtimedependent} also holds for several modified estimators. Details and proofs can be
found in Dahlhaus (2009), Section 4. \hfill $\Box$
\end{remark}
We now briefly mention the exponential inequality. Since this is a non-asymptotic result it holds
regardless whether $\phi$ depends on $T$. Let $\rho_{2,T} (\phi): = \big( \frac{1}{T}
\sum^{T}_{t=1} \int^{\pi}_{-\pi} \phi (\frac{t}{T}, \lambda)^{2} d\lambda \big)^{1/2}$.
\begin{theorem} [Exponential inequality] \label{Exp-Ineq}
Under suitable regularity conditions we have for all $\eta>0$
\begin {equation}
\label{exp-ineq1}
 P\Big(\,\Big|\sqrt{T} \, \big(F_{T}
(\phi) - \mean F_{T} (\phi) \big) \Big| \geq \eta \Big) \leq c_{1} \exp \Big(-c_{2} \sqrt{
\frac{\eta}{\rho_{2,T} (\phi)}} \; \Big)
\end{equation}
with some constants $c_{1}$, $c_{2} > 0$ independent of $T$.
\end{theorem}
This result is proved in Dahlhaus and Polonik (2009), Theorem 2.7. There exist several versions of
this result - for example in the Gaussian case it is possible to omit the $\sqrt{\cdot}$ in
(\ref{exp-ineq1}) or to prove a Bernstein-type inequality which is even stronger (cf. Dahlhaus and
Polonik, 2006, Theorem 4.1).

Subsequently, a maximal inequality, i.e. an exponential inequality for $\sup_{\phi \in
\Phi}|E_T(\phi)|$ has been proved in Dahlhaus and Polonik (2009), Theorem 2.9 under conditions on
the metric entropy of the corresponding function class $\Phi$. We refer to that paper for details.

With the maximal inequality tightness of the empirical spectral process can be proved leading to a
functional central limit theorem for the empirical spectral process indexed by a function class
(cf. Dahlhaus and Polonik (2009), Theorem 2.11). Furthermore a Glivenko Cantelli type result for
the empirical spectral process can be obtained (Theorem 2.12).

Other applications of the maximal inequality are for example uniform rates of convergence for
different estimates. As an example we now state a uniform convergence result for the local quasi
likelihood estimate $\hat{\boldsymbol{\theta}}_{T}^{GW}\!(u_0)$
from (\ref{GeneralizedKernelEstimate}).

\begin{theorem} \label{uniform-convergence-Whittle}
Let $X_{t,T}$ be a locally stationary process with $\mu (\cdot) \equiv 0$. Under suitable
regularity conditions (in particular under the assumption that $f_{\boldsymbol{\theta}}(\lambda)$
is twice differentiable in $\boldsymbol{\theta}$ with uniformly Lipschitz continuous derivatives in
$\lambda$) we have for $b_T T >> (\log T)^{6}$
\begin {equation*}
\sup_{u_0 \,\in \,[b_{T}/2\,,\,1-b_{T}/2]} \big\|
 \hat{\boldsymbol{\theta}}^{GW}_{T}(u_0) - \boldsymbol{\theta}_{0} (u) \big\|_{2}= O_{p} \Big(\frac {1}
{\sqrt{b_{T} T}} + b_{T}^{2}\Big)\,,
\end {equation*}
that is for $b_{T} \sim T^{-1/5}$ we obtain the uniform rate $O_{p}(T^{-2/5})$.
\end{theorem}

\proofbegin The result has been proved in Dahlhaus and Polonik (2009), Theorem 4.1.

\begin{example} [Testing for stationarity] \label{Testing Stationarity} $\;$
Another application of the maximal inequality is the derivation of a functional central limit for
the empirical spectral process. A possible application is a test for stationarity. We briefly
present the idea - although we clearly mention that the construction below is finally not
successful. The idea for a test of stationarity is to test whether the time varying spectral
density $f(u,\lambda)$ is constant in $u$. This is for example achieved by the test statistic
\begin{equation}\label{test statistic stationarity}
\sqrt{T} \sup_{u \in[0,1]} \sup_{\lambda \in [0,\pi]} \, \Big| F_{T}
(u,\lambda) - u \, F_{T} (1,\lambda) \Big|
\end{equation}
where
\begin{equation*}\label{A7sum}
F_{T} (u,\lambda) :=  \frac{1}{T} \sum^{[uT]}_{t=1}
\int_{0}^{\lambda}  J_{T} (\frac{t}{T}, \mu) \, d\mu
\end{equation*}
is an estimate of the integrated time frequency spectral density $F(u,\lambda):=\int_{0}^{u}
\int_{0}^{\lambda} f(v,\mu) \, d\mu dv$, and
\begin{equation*}\label{A7sum}
u \, F_{T} (1,\lambda) =  u \int_{0}^{\lambda}
I_{T} ( \mu) \, d\mu
\end{equation*}
is the corresponding estimate of $F(u,\lambda)$ under the hypothesis of stationarity where
$f(v,\mu)=f(\mu)$. Under the hypothesis of stationarity we have
\begin{equation*}\label{A7sum}
F (u,\lambda) - u \, F (1,\lambda) = \int_{0}^{1} \int_{0}^{\lambda}
\big(I_{[0,u]} (v) - u \big)  f(\mu) \,d\mu \, dv = 0
\end{equation*}
and therefore
\begin{equation*}\label{A7sum}
\sqrt{T} \Big(F_{T} (u,\lambda) - u \, F_{T} (1,\lambda) \Big) = E_{T}
(\phi_{u,\lambda})
\end{equation*}
with $\phi_{u,\lambda}(v,\mu) = \big(I_{[0,u]} (v) - u \big)\, I_{[0,\lambda]} (\mu)$.  We now need
functional convergence of $E_{T} (\phi_{u,\lambda})$. Convergence of the finite dimensional
distributions follows from Theorem~\ref{clt} above. Tightness and therefore the functional
convergence follows from Theorem 2.11 of Dahlhaus and Polonik (2009). As a consequence we obtain
under the null hypothesis
\begin{equation*} \label{}
\sqrt{T} \, \Big( F_{T} (u,\lambda) - u \, F_{T} (1,\lambda)\Big)_{u
\in [0,1], \lambda \in [0,\pi]} \dconv E \big(u,\lambda\big)_{u \in
[0,1], \lambda \in [0,\pi]}
\end{equation*}
where $E \big(u,\lambda\big)$ is a Gaussian process. If $\kappa_4=0$ (Gaussian case) and $f(\mu)=c$
it can be shown that this is the Kiefer-M\"{u}ller process. However, for general $f$ it is a difficult
and unsolved task to calculate or estimate the limit distribution and in particular the
distribution of the test statistic in (\ref{test statistic stationarity}). This may be done by
transformations (like $U_p$ - or $T_p$ - type transforms) and/or by finding an adequate bootstrap
method.

We mention that Paparoditis (2009, 2010) has given two different solutions of this testing problem. \hfill $\Box$
\end{example}

\section{Additional topics and further references} \label{SectionAdditionalTopics}

This section gives an overview over additional topics with further references. We concentrate on work which uses the infill asymptotic approach of local stationarity. Even in this case it is not possible to give a complete overview.

\underline{1. Locally stationary wavelet processes:} There exists a large number of papers on the use of wavelets for modeling locally stationary processes. The first type of application is to estimate the parameter curves via the use of wavelets. This has been mentioned a few times in the above presentation (cf. (\ref{OrthogonalSeriesEstimate})).

A breakthrough for the application of wavelets to nonstationary processes was the introduction of ``locally stationary wavelet processes'' by Nason et.al. (2000). This class is somehow the counterpart to the representation (\ref{LSFrepresentation}) for locally stationary processes. It  also uses a rescaling argument - thus making all methods for these processes accessible to a meaningful asymptotic theory. Locally stationary wavelet processes are processes with
the wavelet representation
\begin{equation} \label{} \label{LSWrepresentation}
X_{t,T} = \mu\big(\frac{t}{T}\big) + \sum_{j=1}^{\infty} \sum_{k=-\infty}^{\infty} w_{j,k;T} \,\psi_{j,k-t} \,\xi_{j,k}
\end{equation}
where $\{\xi_{j,k}\}$ are a collection of uncorrelated random variables with mean $0$ and variance $1$, the $\{\psi_{j,t}\}$ are a set of discrete nondecimated wavelets (compactly supported oscillatory vectors with
support proportional to $2^{j}$), and $\{w_{j,k;T}\}$ are a collection of amplitudes that are smooth in a particular way as a function of $k$. The smoothness of $\{w_{j,k;T}\}$ controls the degree of local stationarity of $X_{t,T}$.  The spectrum is linked to the process by $\{w_{j,k;T}\} \approx S_{j} \big(\frac {k} {T}\big)$. Nason et.al. (2000) also define the ``evolutionary wavelet spectrum'' and show how this can be estimated by a smoothed wavelet periodogram. In addition this leads to an estimate of the local covariance. An introduction to LSW-processes and an overview on early results for such processes can be found in Nason and von Sachs (1999). Fryzlewicz and Nason (2006) suggest the use of a Haar-Fisz method for the estimation of evolutionary wavelet spectra by combining Haar wavelets and the variance stabilizing Fisz transform. Van Bellegem and von Sachs (2008) consider wavelet processes whose spectral density function changes very quickly in time. By using a wavelet-type transform of the autocovariance function with respect to so-called autocorrelation wavelets they propose a pointwise adaptive estimator of the time varying autocovariance and the time varying spectrum.

Furthermore, several papers mentioned below use the framework of LSW-processes.

\medskip
\underline{2. Multivariate locally stationary processes:} We first mention that in particular the Gaussian likelihood theory for locally stationary processes from Section~\ref{SectionLikelihood} also holds for multivariate processes -- see Dahlhaus (2000).

Beyond that Chiann and Morettin (1999, 2005) investigate the estimation of time varying coefficients of a linear system where the input and output are locally stationary processes. They study different estimation techniques in the frequency- and time domain.

Ombao et al. (2001) analyze bivariate nonstationary time series. They use SLEX functions (time-localized generalization of the Fourier waveform) and propose a method that automatically segments the time series into approximately stationary blocks and selects the span to be used to obtain the smoothed estimates of the time varying spectra and coherence. Ombao et.al. (2005) use  the SLEX framework to build a family of multivariate models that can explicitly characterize the time varying spectral and coherence properties of a multivariate time series. Ombao and Van Bellegem (2008) estimate the time varying coherence by using time-localized linear filtering. Their method automatically selects via tests of homogeneity the optimal window width for estimating local coherence.

Motta et.al. (2011) propose a locally stationary factor model for large cross-section and time dimensions. Factor loadings are estimated by the
eigenvectors of a nonparametrically estimated covariance matrix. Eichler et. al. (2011) investigate dynamic factor modeling of locally stationary processes. They estimate the common components of the dynamic factor model by the eigenvectors of an estimator of the time varying spectral density matrix. This can also be seen as a time varying principal components approach in the frequency domain.

Cardinali and Nason (2010) introduce the concept of costationary of two locally stationary time series where some linear combination of the two processes is stationary. They show that costationarity imply a error-correction type of formula in which changes in the variance of one series are reflected by simultaneous balancing changes in the other. Sanderson et.al. (2010) propose a new method of measuring the dependence between non-stationary time series based on a wavelet coherence between two LSW-processes.

\medskip
\underline{3. Testing of locally stationary processes -- in particular tests for stationarity:}
Among the \linebreak large literature on testing there is a considerable part devoted to testing of stationarity. Tests of stationarity have already been proposed and theoretically investigated before the framework of local stationarity was created. In that cases the theoretical investigations mainly consisted in the investigation of the asymptotic distribution of the test statistics under the hypothesis of stationarity.

Priestley and Subba Rao (1969) proposed testing the homogeneity of a set of evolutionary spectra evaluated at different points of time. For Gaussian processes and for the purpose of a change point detection Picard (1985) developed a test based on the difference between spectral distribution functions estimated on different parts of the series and evaluated using a supremum type statistic. Giraitis and Leipus (1992) generalized this approach to the case of linear processes. Von Sachs and Neumann (2000) developed a test of stationarity based on empirical wavelet coefficients estimated using localized versions of the periodogram. Paparoditis (2009) developed a nonparametric test for stationarity against the alternative of a smoothly time varying spectral structure based
on a local spectral density estimate. He also investigated the power under the fixed alternative of a locally stationary processes. Paparoditis (2010) tested the assumption of stationarity by evaluating the supremum over time of an $L_2$-distance between the local periodogram over a rolling segment and an
estimator of the spectral density obtained using the entire time series at hand. The critical values of a supremum type test are obtained using a stationary bootstrap procedure. Dwivedi and Subba Rao (2011) construct a Portmanteau type test statistic for testing stationarity of a time series by using the property that the discrete Fourier transforms of a time series at the canonical frequencies are asymptotically uncorrelated if and only if
the time series is second-order stationary.

Tests of general hypothesis are derived in Sakiyama and Taniguchi (2003) who test parametric composite hypothesis by the Gaussian likelihood ratio test, the Wald test and the Lagrange multiplier test. Sergides and Paparoditis (2009) develop tests of the hypothesis that the time varying spectral density has a semiparametric structure. The test introduced is based on a $L_2$-distance measure in the spectral domain. As a special case they test for the presence of a tvAR model. A bootstrap procedure is applied to approximate more accurately the distribution of the test statistic under the null hypothesis. Preu{\ss} et. al. (2011) also test semiparametric hypotheses. Their method is based on an empirical version of the $L_2$-distance between the true time varying spectral density and its best approximation under the null hypothesis.

Zhou and Wu (2010) construct simultaneous confidence tubes for time varying
regression coefficients in functional linear models. Using a Gaussian approximation result for non-stationary multiple time series, they show that the constructed simultaneous confidence tubes have asymptotically correct nominal coverage probabilities.

\medskip
\underline{4. Bootstrap methods for locally stationary processes:} Bootstrap methods are in particular needed to derive the asymptotic distribution of test statistics. A time domain local block bootstrap procedure for locally stationary processes has been proposed by Paparoditis and Politis (2002) and by Dowla et al. (2003). Sergides and Paparoditis (2008) develop a method to bootstrap the local periodogram. Their method generates pseudo local periodogram ordinates by combining a parametric time and nonparametric
frequency domain bootstrap approach. They first fit locally a time varying
autoregressive model to capture the essential characteristics of the underlying process. A locally calculated non-parametric correction in the frequency domain is then used so as to improve upon the locally parametric autoregressive fit. Kreiss and Paparoditis (2011) propose a nonparametric bootstrap method by generating pseudo time series which mimic the local second and fourth order moment structure of the underlying process. They prove a bootstrap central limit theorem for a general class of preperiodogram based statistics.

\medskip
\underline{5. Model misspecification and model selection:} Model selection criteria have been heuristically suggested many times for time varying processes -- c.f. Ozaki and Tong (1975); Kitagawa and Akaike (1978) and Dahlhaus (1996b, 1997) among others -- in all papers AIC-type criteria have been suggested for different purposes.

Van Bellegem and Dahlhaus (2006) consider semiparametric estimation and estimate the Kullback-Leibler distance between the semiparametric model and the true process. They use this estimate then as a model selection criterion. Hirukawa et.al. (2008) propose a generalized information criterion based on nonlinear functionals of the time varying spectral density. Chandler (2010) investigates how time varying parameters affect order selection.

Another interesting aspect is that many results of this paper also hold under model - misspecification -- for example Theorem~\ref{AsymptoticPropertiesGWhittle} and the corresponding result for the Block Whittle estimate from (\ref{BlockWhittleEstimate}). An important example is the case where a stationary model is fitted and the underlying process in truth is only locally stationary - see Example~\ref{ModelStationary} and the more detailed discussion for stationary Yule-Walker estimates in Dahlhaus (1997), Section 5.

\medskip
\underline{6. Likelihood theory and large deviations:} Local asymptotic normality (LAN) is derived in the parametric Gaussian case in Dahlhaus (1996b) and Dahlhaus (2000) (cf. Remark~3.3 in that paper). A nonparametric LAN-result is derived in Sakiyama and Taniguchi (2003) and a LAN result under non-Gaussianity in Hirukawa and Taniguchi (2006). In both papers the results are applied to asymptotically optimal estimation and testing. For some statistics also the asymptotic distribution under contiguous alternatives is derived. Tamaki (2009) studies second order asymptotic efficiency of appropriately modified maximum likelihood estimators for Gaussian locally stationary processes.

Large
deviations principles for quadratic forms of locally stationary processes are derived in Zani (2002) including applications to local spectral density and covariance estimation. Wu and Zhou (2011) obtain an invariance principle for non-stationary vector-valued stochastic processes. They show that the partial sums of non-stationary processes can be approximated on a richer probability space by sums of independent Gaussian random vectors.

\medskip
\underline{7. Recursive estimation:} Recursive estimation algorithms are of the form
\begin {equation} \label{estalgo}
\widehat{\theta}_{t} = \widehat{\theta}_{t-1} + \lambda_{t} \;
\psi({\ensuremath {\boldsymbol X}_{t}},\widehat{\theta}_{t-1})
\end{equation}
where ${\ensuremath {\boldsymbol X}_{t}}=(X_1,\dots,X_t)'$. The
recursive structure yields an update of the estimate as soon as
the next observation becomes available and the estimate therefore
is particularly of importance in an online situation. For
stationary processes the algorithm is used with $\lambda_t \sim
1/t$ while in nonstationary situations one uses a nondecreasing
$\lambda$ (constant stepsize case) that is the estimate puts
stronger weights on the last observations.

Adaptive estimates of the above type have been investigated over
the last 30 years in different scientific communities: by system
theorists under the name ``recursive identification methods'' (cf.
Ljung (1977); Ljung and S\"{o}derstr\"{o}m (1983)), in the stochastic
approximation community (cf. Benveniste, M\'{e}tivier and Priouret
(1990); Kushner and Yin (1997)), in the neural network community
under the name ``back-propagation algorithm'' (cf. White (1992) or
by Haykin (1994)), and in applied sciences, particularly for
biological and medical applications (cf. Schack and Grieszbach
(1994)).

The properties of recursive estimation algorithms have rigorously
been investigated in many papers under the premise that the underlying true
process is stationary. However, for nonstationary processes and
the constant stepsize case there did not exist for a long time a reasonable framework to study theoretically the properties of these algorithms. This has changed with the concept of locally stationary processes with it's infill asymptotics which now allows for theoretical investigations of these algorithms.

In Moulines et.al. (2005) the properties of recursive estimates of tvAR-processes have been investigated in the framework of locally stationary processes. The asymptotic properties of the estimator have been proved including a minimax result. In Dahlhaus and Subba Rao (2007) a recursive  algorithm for estimating the parameters of a tvARCH-process has been proposed. Again the asymptotic properties of the estimator have been proved.

\medskip
\underline{8. Inference for the mean curve:} Modeling the time varying mean of a locally stationary process is an important task which has not been discussed in this overview. In principle nearly all known techniques from nonparametric regression may be used such as kernel estimates, local polynomial fits, wavelet estimates or others. The situation is however much more challenging since the ``residuals'' are in this case a locally stationary process which usually is modeled at the same time.

In general the topic needs more investigation. Dahlhaus (1996a, 1996b, 1997, 2000) and Dahlhaus and Neumann (2001) contain also results where the mean is time varying and/or estimated. A more detailed investigation is contained in Tunyavetchakit (2010) in the context of time varying AR(p)-processes where the mean curve is estimated in parallel and the optimal segment length is determined similar to (\ref{NoptTVAR}).

\medskip

\underline{9. Piecewise constant models:} Davis et.al. (2005) consider the problem of modeling a class of nonstationary time series using piecewise constant AR-processes. The number and locations of the piecewise AR segments, as well as the orders of the respective AR processes, are determined by the minimum
description length principle. The best combination is then determined by a genetic algorithm. In Davis et.al. (2008) to general parametric time series models for the segments and illustrate the method with piecewise GARCH-models, stochastic volatility and generalized state space models.

Locally constant parametric models have also been considered in a non-asymptotic approach by Mercurio and Spokoiny (2004) and others where the so-called small modeling bias condition is used to determine the length of the interval of time homogeneity and to fit the parameters -- for more details see also Spokoiny (2010).

\medskip
\underline{10. Long memory processes:} Beran (2009) and Palma and Olea (2010) have extended the concept of local stationarity to long-range dependent
processes. While Beran (2009) uses a nonparametric approach with a local least squares estimate similar to (\ref{LocalKernelEstimate}) Palma and Olea (2010) use a parametric approach and use the block Whittle likelihood from (\ref{BlockWhittleLikelihood}). Both papers then investigate the asymptotic properties. Roueff and von Sachs (2011) use a local log-regression
wavelet estimator of the time-dependent long memory parameter and study it's asymptotic properties.

\medskip
\underline{11. Locally stationary random fields:} Fuentes (2001) studies different methods for
locally stationary isotropic random fields with parameters varying across space. In particular she uses local Whittle estimates. Eckley et.al.~(2010) propose the modeling and analysis of image texture by using an extension
of a locally stationary wavelet process model for lattice processes. They construct estimates of a spatially localized spectrum and a localized autocovariance which are then used to characterize textures in a multiscale and spatially adaptive way. Anderes and Stein (2011) develop a weighted local likelihood estimate for the parameters that govern the local spatial dependency of a locally stationary random field.

\medskip
\underline{12. Discrimination Analysis:} Discrimination Analysis for locally stationary processes based on the Kullback-Leibler divergence as a classification criterion has been investigated in Sakiyama and Taniguchi (2004) and for multivariate processes in Hirukawa (2004). Huang et.al. (2004) propose a discriminant scheme based on the SLEX-library and a discriminant criterion that is also related to the Kullback-Leibler divergence. Chandler and Polonik (2006) develop methods for the discrimination of locally stationary processes based on the shape of different features. In particular they use shape measures of the variance function as a criterion for discrimination and apply their method to the discrimination of earthquakes and explosions. Fryzlewicz and Ombao (2009) use a bias-corrected non-decimated wavelet transform for classification in the framework of LSW-processes.

\medskip
\underline{13. Prediction:} Fryzlewicz et.al. (2003) address the problem of how to forecast non-stationary time series by means of non-decimated wavelets. Using the class of LSW-processes they introduce a new predictor based on wavelets and derive the prediction equations as a generalization of the Yule-Walker equations. Van Bellegem and von Sachs (2004) apply locally stationary processes to the forecasting of several economic data sets such as returns and exchange rates.

\medskip
\underline{14. Finance:} There is a growing interest in finance for models with time varying parameters. An overview on locally stationary volatility models is given in Van Bellegem (2011). A general discussion on local stationary in different areas of finance can be found in Guegan (2007) -- see also Taniguchi et. al. (2008). For example, many researchers are convinced that the observed slow decay of the sample autocorrelation function of absolute stock returns is not a long memory effect but due to nonstationary changes in the unconditional variance (c.f. Mikosch and St$\breve{a}$ric$\breve{a}$ (2004), St$\breve{a}$ric$\breve{a}$ and Granger (2005), Fryzlewicz et. al. (2006)) leading for example to GARCH-models with time varying parameters.

References for work on tvGARCH-models have been given in Section~\ref{NonlinearProcAndDerivative}. Other work on applications of locally stationary processes in finance is for example the work on optimal portfolios with locally stationary returns of assets by Shiraishi and Taniguchi (2007). Hirukawa (2006) uses locally stationary processes for a clustering problem of stock returns. Fryzlewicz (2005) models some stylized facts of financial log returns by LSW-processes. Fryzlewicz et. al. (2006) consider a locally stationary model for financial log-returns and propose a wavelet thresholding algorithm for volatility estimation, in which Haar wavelets are combined with the variance-stabilizing Fisz transform.

\medskip
\underline{15. Further topics:} Robinson (1989) uses also the infill asymptotics approach in his work on \underline{nonparametric regression} with time varying coefficients. Orbe et al. (2000) estimate nonparametrically a time varying coefficients model allowing for seasonal and smoothness constraints. Orbe et.al. (2005) estimate the time varying coefficients under shape restrictions over and for locally stationary regressors. Chiann and Morettin (2005) investigate the estimation of coefficient curves in time varying linear systems.

Estimation of time varying \underline{quantile curves} for nonstationary processes has been done in Draghicescu et.al. (2009) and Zhou and Wu (2009). Specification tests of time varying quantile curves have been investigated in Zhou (2010).

\bigskip

\newpage

\noindent{\Large \bf References}
\begin{description}
\baselineskip1.2em
\itemsep-0.05cm

\item Amado, C., and Ter\"{a}svirta, T. (2011). Modelling volatility with variance decomposition. CREATES Research Paper 2011-1, Aarhus University.

\item Anderes, E.B. and Stein, M.L. (2011). Local likelihood estimation for nonstationary random fields. {\em Journal of Multivariate Analysis} {\bf 102}, 506--520.

\item Benveniste, A., M\'{e}tivier, M. and Priouret, P. (1990). {\em Adaptive Algorithms and Stochastic Approximations}. Springer Verlag, Berlin.

\item Beran, J. (2009). On parameter estimation for locally stationary long-memory processes. {\em J. Statist. Plann. Inference} {\bf 139}, 900--915.

\item Berkes, I., Horv\'ath, L. and Kokoskza, P. (2003). GARCH processes: structure and estimation. {\em Bernoulli} {\bf 9}, 201--207.

\item Brillinger, D.R. (1981). {\em Time Series: Data Analysis and Theory.} Holden Day, San Francisco.

\item Brockwell P.J., and Davis R.A. (1991). {\em Time Series: Theory and Methods}, 2nd ed. Springer-Verlag, New York.

\item Cardinali, A. and Nason, G. (2010). Costationarity of locally stationary time series. {\em Journal of Time Series Econometrics} {\bf  2}, No. 2, Article 1. DOI: 10.2202/1941-1928.1074

\item Chandler, G. (2010). Order selection for heteroscedastic autoregression: A study on concentration. {\em Statistics and Probability Letters} {\bf 80}, 1904--1910.

\item Chandler, G. and Polonik, W. (2006). Discrimination of locally stationary time series based on the excess mass functional.  {\em J. Amer. Statist. Assoc.} {\bf 101}, 240--253.

\item Chiann, C. and Morettin, P. (1999). Estimation of time varying linear systems. {\em Statist. Inference Stoch. Proc.} {\bf 2}, 253--285.

\item Chiann, C. and Morettin, P. (2005). Time-domain estimation of time-varying linear systems. {\em J. Nonpar. Statist.} {\bf 17}, 365--383.

\item Dahlhaus, R. (1988). Empirical spectral processes and their applications to time series analysis. {\em Stoch. Proc. Appl.} {\bf 30}, 69--83.

\item Dahlhaus, R. (1996a). On the Kullback-Leibler information divergence for locally stationary processes. {\em Stoch. Proc. Appl.} {\bf 62}, 139--168.

\item Dahlhaus, R. (1996b). Maximum likelihood estimation and model selection for locally stationary processes. {\em J. Nonpar. Statist.} {\bf 6}, 171--191.

\item Dahlhaus, R. (1996c). Asymptotic statistical inference for nonstationary processes with evolutionary spectra. In: {\em Athens Conference on Applied Probability and Time Series Vol II}. (P.M. Robinson and M. Rosenblatt, eds.), 145--159, Lecture Notes in Statistics 115, Springer, New York.

\item Dahlhaus, R. (1997). Fitting time series models to nonstationary processes. {\em Ann. Statist.} {\bf 25}, 1--37.

\item Dahlhaus, R. (2000). A likelihood approximation for locally stationary processes. {\em Ann. Statist.} {\bf 28}, 1762--1794.

\item Dahlhaus, R. (2009). Local inference for locally stationary time series based on the empirical spectral measure. {\em J. Econometrics} {\bf 151}, 101--112.

\item Dahlhaus, R. and Giraitis, L. (1998). On the optimal segment length for parameter estimates for locally stationary time series. {\em J. Time Series Anal.} {\bf 19}, 629--655.

\item Dahlhaus, R., Neumann, M.H. and von Sachs, R. (1999). Nonlinear wavelet estimation of time-varying autoregressive processes {\em Bernoulli} {\bf 5}, 873--906.

\item Dahlhaus, R. and Neumann, M.H. (2001). Locally adaptive fitting of semiparametric models to nonstationary time series. {\em Stoch. Proc. and Appl.} {\bf 91}, 277--308.


\item Dahlhaus, R. and Polonik, W. (2006). Nonparametric quasi maximum likelihood estimation for Gaussian locally stationary processes. {\em Ann. Statist.} {\bf 34}, 2790--2824.

\item Dahlhaus, R. and Polonik, W. (2009). Empirical spectral processes for locally stationary time series. {\em Bernoulli} {\bf 15}, 1--39.

\item Dahlhaus, R. and Subba Rao, S. (2006). Statistical inference for locally stationary ARCH models. {\em Ann. Statist.} {\bf 34}, 1075--1114.

\item Dahlhaus R. and Subba Rao S. (2007). A recursive online algorithm for the estimation of time-varying ARCH parameters. {\em Bernoulli} {\bf 13}, 389--422.

\item Davis, R.A. and Lee, T., and Rodriguez-Yam, G. (2005). Structural break estimation for nonstationary time series models. {\em J. Amer. Statist. Assoc.} {\bf 101}, 223--239.

\item Davis, R.A., Lee, T., and Rodriguez-Yam, G. (2008). Break detection for a class of nonlinear time series models. {\em J. Time Ser. Anal.} {\bf 29}, 834--867.


\item Dowla, A., Paparoditis, E. and Politis, D.N. (2003). Locally stationary processes and the local bootstrap. In: {\em Recent Advances and Trends in Nonparametric Statistics} (Eds. M. G. Akritas and D. N. Politis). Elsevier Science B.V., Amsterdam, 437--445.

\item Draghicescu, D., Guillas, S. and Wu, W.B. (2009). Quantile curve estimation and visualization for non-stationary time series. {\em J. Comput. Graph. Statist.} {\bf 18}, 1--20.

\item Dunsmuir, W. (1979). A central limit theorem for parameter estimation in stationary vector time series and its application to models for a signal observed with noise. {\em Ann. Statist.} {\bf 7}, 490--506.

\item Dzhaparidze, K. (1971). On methods for obtaining asymptotically efficient spectral parameter estimates for a stationary Gaussian process with rational spectral density. {\em Theory Probab. Appl.} {\bf 16}, 550--554.

\item Dzhaparidze, K. (1986). {\em Parameter Estimation and Hypothesis Testing in Spectral Analysis of Stationary Time Series.} Springer Verlag, New York.

\item Dwivedi, Y. and Subba Rao, S. (2011). A test for second-order stationarity of a time series based on the discrete Fourier transform. {\em J. Time Series Anal.} {\bf 32} 68--91.

\item Eckley, I.A., Nason, G.P. and Treloar, R.L. (2010). Locally stationary wavelet fields with application to the modelling and analysis of image texture. {\em Appl. Statist.} {\bf 59}, 595--616.

\item Eichler, M., Motta, G. and von Sachs, R. (2011). Fitting dynamic factor models to non-stationary time series. {\em J. Econometrics} {\bf 163}, 51-70.

\item Fay, G. and Soulier, P. (2001) The periodogram of an i.i.d. sequence. {\em Stoch. Proc. Appl.} {\bf 92}, 315--343.

\item Fox, R. and Taqqu, M.S. (1986). Large-sample properties of parameter estimates for strongly dependent stationary Gaussian time series. {\em Ann. Statist.} {\bf 14}, 517--532.

\item Fryzlewicz, P. (2005) Modelling and forecasting financial log-returns as locally stationary wavelet processes. {\em J. Appl. Statist.} {\bf 32}, 503--528.

\item Fryzlewicz, P. and Nason, G. P. (2006). Haar-Fisz estimation of evolutionary wavelet spectra. {\em J. R. Statist. Soc. B} {\bf 68}, 611--634.

\item Fryzlewicz, P. and Ombao, H. (2009). Consistent classification of nonstationary time series using stochastic wavelet representations.  {\em J. Amer. Statist. Assoc.} {\bf 104}, 299-312.

\item Fryzlewicz, P., Sapatinas, T. and Subba Rao, S. (2006). A Haar-Fisz technique for locally stationary volatility estimation. {\em Biometrika} {\bf 93}, 687--704.

\item Fryzlewicz, P., Sapatinas, T. and Subba Rao, S. (2008). Normalised least-squares estimation in time-varying ARCH models. {\em Ann. Statist.} {\bf 36}, 742--786.

\item Fryzlewicz, P. and Subba Rao, S. (2011). On mixing properties of ARCH and time-varying ARCH processes. {\em Bernoulli} {\bf 17}, 320--346.

\item Fryzlewicz, P., Van Bellegem, S. and von Sachs, R. (2003). Forecasting non-stationary time series by wavelet process modeling. {\em Ann. Inst. Statist. Math.} {\bf 55}, 737--764.

\item Fuentes, M. (2001). A high frequency kriging approach for non-stationary environmental processes. {\em Environmetrics} {\bf 12}, 469--483.

\item Giraitis, L. and Leipus, R. (1992). Testing and estimating in the change-point problem of the spectral function. {\em Lithuanian Mathematical Journal} {\bf 32}, 15--29.

\item Granger, C.W.J. and Hatanaka, M. (1964). {\em Spectral Analysis of Economic Time Series}. Princeton University Press. Princeton, New Jersey.

\item Grenander, U. and Szeg\"o, G. (1958). {\em Toeplitz Forms and their Applications.} University of California Press, Berkeley.

\item Grenier, Y. (1983). Time dependent ARMA modelling of nonstationary signals. {\em IEEE Trans. Acoust. Speech Signal Process.} {\bf 31}, 899--911.

\item Gu\'{e}gan D. (2007). Global and local stationary modelling in finance: Theory and empirical evidence. Preprint, Centre d'Economique de la Sorbonne.

\item Guo, W., Dai, M., Ombao, H.C. and von Sachs, R. (2003). Smoothing spline ANOVA for time-dependent spectral analysis. {\em J. Amer. Statist. Assoc.} {\bf 98}, 643--652.

\item Hannan, E.J. (1973). The asymptotic theory of linear time series models. {\em J. Appl. Prob.} {\bf 10}, 130--145.

\item Hallin, M. (1986). Nonstationary q-dependent processes and time-varying moving average models: invertibility properties and the forecasting problem. {\em Adv. Appl. Probab.} {\bf 18}, 170--210.

\item Hirukawa, J. (2004). Discriminant analysis for multivariate non-Gaussian locally stationary processes. {\em Scientiae Mathematicae Japonicae Online} {\bf 10}, 235--258.

\item Hirukawa, J. (2006). Cluster analysis for non-Gaussian locally stationary processes. {\em Intern. J. Theor. Appl. Finance} {\bf 9}, 113--132.

\item Hirukawa, J., Kato, H.S., Tamaki, K. and Taniguchi, M. (2008). Generalized information criteria in model selection for locally stationary processes. {\em J. Japan Statist. Soc.} {\bf 38}, 157--171.

\item Hirukawa, J. and Taniguchi, M. (2006). LAN theorem for non-Gaussian locally stationary processes and its applications. {\em J. Statist. Planning Infer.} {\bf 136}, 640--688.

\item Hosoya, Y. and Taniguchi, M. (1982). A central limit theorem for stationary processes and the parameter estimation of linear processes {\em Ann. Statist.} {\bf 10}, 132--153.

\item Huang, H.-Y. Ombao, H.C. and Stoffer, D.S. (2004). Discrimination and Classification of Nonstationary Time Series Using the SLEX Model. {\em J. Amer. Statist. Assoc.} {\bf 99}, 763--774.

\item Jentsch, C. (2006). Asymptotik eines nicht-parametrischen Kernsch\"{a}tzers f\"{u}r zeitvariable autoregressive Prozesse. Diploma thesis, University of Braunschweig.

\item Kayhan, A., El-Jaroudi, A. and Chaparro, L. (1994). Evolutionary periodogram for nonstationary signals. {\em IEEE Trans. Signal Process.} {\bf 42}, 1527--1536.

\item Kim, W. (2001). Nonparametric kernel estimation of evolutionary autoregressive processes. Discussion paper 103. Sonderforschungsbereich 373, Berlin.

\item Kitagawa, G. and Akaike, H. (1978). A Procedure for The Modeling of Non-Stationary Time Series. {\em Ann. Inst. Statist. Math.} {\bf 30 B}, 351--363.

\item Kitagawa, G. and Gersch, W. (1985). A smoothness priors time-varying AR coefficient modeling of the nonstationary covariance time series. {\em IEEE Trans. Automat. Control.} {\bf 30}, 48-56.

\item Koo, B. and Linton, O. (2010). Semiparametric estimation of locally stationary diffusion models. LSE STICERD Research Paper No. EM/2010/551.

\item Kreiss, J.-P. and Paparoditis, E. (2011). Bootstrapping Locally Stationary Processes. Technical report.

\item Kushner, H. J. and Yin, G.G. (1997). {\em Stochastic Approximation Algorithms and Applications}. Springer Verlag, New York.

\item Ljung, L. (1977). Analysis of recursive stochastic algorithms. {\em IEEE Trans. Automatic Control} {\bf 22}, 551--575.

\item Ljung, L. and S\"{o}derstr\"{o}m, T. (1983). {\em Theory and Practice of Recursive Identification}.  MIT Press, Cambridge, MA.

\item Martin, W. and Flandrin, P. (1985). Wigner-Ville spectral analysis of nonstationary processes. {\em IEEE Trans. Acoust. Speech Signal Process.} {\bf 33}, 1461--1470.

\item M\'{e}lard, G. and A. Herteleer-de-Schutter, A. (1989). Contributions to evolutionary spectral theory. {\em J. Time Series Anal.} {\bf 10} 41--63.

\item Mercurio, D. and Spokoiny, V. (2004). Statistical inference for time-inhomogenous volatility models. {\em Ann. Statist.} {\bf 32}, 577--602.

\item Mikosch, T., Gadrich, T., Kl\"uppelberg, C., Adler, R.J. (1995). Parameter estimation for ARMA models with infinite variance innovations. {\em Ann. Statist.} {\bf 23}, 305--326.

\item Mikosch, T. and Norvaisa, R. (1997). Uniform convergence of the empirical spectral distribution function. {\em Stoch. Proc. Appl.} {\bf 70}, 85--114.

\item Mikosch, T., and C. St$\breve{a}$ric$\breve{a}$, C. (2004). Nonstationarities in financial time series, the long-range dependence, and the IGARCH effects. {\em The Review of Economics and Statistics} {\bf 86}, 378--390.

\item Motta, G., Hafner, C.M. and von Sachs, R. (2011). Locally stationary factor models: Identification and nonparametric estimation. {\em Econometric Theory} {\bf 27}, 1279--1319 doi:10.1017/S026646661100005

\item Moulines, E., Priouret, P. and Roueff, F. (2005). On recursive estimation for locally stationary time varying autoregressive processes. {\em Ann. Statist.} {\bf 33}, 2610--2654.

\item Nason, G.P. and von Sachs, R. (1999). Wavelets in time series analysis. {\em Phil. Trans. R. Soc. Lond. A} {\bf 357}, 2511--2526.

\item Nason, G.P., von Sachs, R. and Kroisandt, G. (2000). Wavelet processes and adaptive estimation of evolutionary wavelet spectra. {\em J. Royal Statist. Soc. B} {\bf 62}, 271--292.

\item Neumann, M.H. and von Sachs, R. (1997). Wavelet thresholding in anisotropic function classes and applications to adaptive estimation of evolutionary spectra. {\em Ann. Statist.} {\bf 25}, 38--76.

\item Ombao, H.C., Raz, J.A., von Sachs, R. and Malow, B.A. (2001). Automatic statistical analysis of bivariate nonstationary time series. {\em J. Amer. Statist. Assoc.} {\bf 96}, 543--560.

\item Ombao, H.C., von Sachs, R. and Guo, W. (2005). The SLEX analysis of multivariate non-stationary time series. {\em J. Amer. Statist. Assoc.} {\bf 100}, 519--531.

\item Ombao, H.C. and Van Bellegem, S. (2008). Evolutionary Coherence of Nonstationary Signals. {\em IEEE Transactions on Signal Processing} {\bf 56}, 2259--2266.

\item Orbe, S., Ferreira, E., Rodriguez-Poo, R.M. (2000). A nonparametric method to estimate time varying coefficients. {\em J. Nonparam. Statist.} {\bf 12}, 779--806.

\item Orbe, S., Ferreira, E., Rodriguez-Poo, R.M. (2005). Nonparametric estimation of time varying parameters under shape restrictions. {\em Journal of Econometrics} {\bf 126}, 53--77.

\item Ozaki, T. and Tong, H. (1975). On the fitting of non-stationary autoregressive models in time series analysis. {\em Proceedings of the 8-th Hawaii International Conference on System Sciences}. Western Periodical Company, North Hollywood, California.

\item Palma, W. and Olea, R. (2010). An efficient estimator for locally stationary Gaussian long-memory processes. {\em Ann. Statist..} {\bf 38}, 2958--2997.


\item Paparoditis, E. (2009). Testing temporal constancy of the spectral structure of a time series. {\em Bernoulli} {\bf 15}, 1190--1221.

\item Paparoditis, E. (2010). Validating stationarity assumptions in time series analysis by rolling local periodograms. {\em J. Amer. Statist. Assoc.} {\bf 105}, 839--851.

\item Paparoditis, E. and Politis, D.N. (2002). Local block bootstrap. {\em C. R. Acad. Sci. Paris, Ser. I} {\bf 335}, 959--962.

\item Parzen, E. (1983). Autoregressive spectral estimation. In: {\em Handbook of Statistics} (D.R. Brillinger and P.R. Krishnaiah, eds.), {\bf 3}, 221--247, North-Holland, Amsterdam.

\item Picard, D. (1985). Testing and estimating change-points in time series. {\em Advances in Applied Probability} {\bf 17}, 841--867.

\item Preu{\ss}, P., Vetter, M., and Dette, H. (2011). Testing semiparametric hypotheses in locally stationary processes. Discussion paper 13/11. SFB 823, TU Dortmund.

\item Priestley, M.B. (1965). Evolutionary spectra and non-stationary processes. {\em J. Roy. Statist. Soc. Ser. B} {\bf 27}, 204--237.

\item Priestley, M. B. and Subba Rao, T. (1969). A test for non-stationarity of time series. {\em Journal of the Royal Statistical Society B} {\bf 31}, 140--149.

\item Priestley, M.B. (1981). {\em Spectral Analysis and Time Series}, Academic Press, London.

\item Priestley, M.B. (1988). {\em Nonlinear and Nonstationary Time Series Analysis}, Academic Press, London.

\item Robinson, P.M., (1989). Nonparametric estimation of time varying parameters. In: Hackl, P. (Ed.), {\em Statistics Analysis and Forecasting of Economic Structural Change}. Springer, Berlin, 253--264.

\item Robinson, P.M. (1995). Gaussian semiparametric estimation of long range dependence. {\em Ann. Statist.} {\bf 23}, 1630--1661.

\item Rosen, O., Stoffer, D.S. and Wood, S. (2009). Local Spectral Analysis via a Bayesian Mixture of Smoothing Splines. {\em Journal of the American Statistical Association} {\bf 104}, 249--262.

\item Roueff, F. and von Sachs, R. (2011). Locally stationary long memory estimation. {\em Stoch.  Proc. Appl.} {\bf 121}, 813--844.

\item von Sachs, R. and Neumann, M. (2000). A wavelet-based test for stationarity. {\em J. Time Ser. Anal.} {\bf 21}, 597--613.

\item Sanderson, J., Fryzlewicz, P. and Jones, M. (2010). Estimating linear dependence between nonstationary time series using the locally stationary wavelet model. {\em Biometrika} {\bf 97}, 435-446.

\item Sakiyama, K. and Taniguchi, M. (2003). Testing composite hypotheses for locally stationary processes. {\em J. Time Ser. Anal.} {\bf 24}, 483--504.

\item Sakiyama, K. and Taniguchi, M. (2004). Discriminant analysis for locally stationary processes. {\em J. Multiv. Anal.} {\bf 90}, 282--300.

\item Schack, B. and Grieszbach, G. (1994). Adaptive methods of trend detection and their application in analyzing biosignals. {\em Biom. J.} {\bf 36}, 429--452.

\item Sergides, M. and Paparoditis, E. (2008). Bootstrapping the Local Periodogram of Locally Stationary Processes. {\em J. Time Ser. Anal.} {\bf 29}, 264--299. Corrigendum: {\em J. Time Ser. Anal.} {\bf 30}, 260--261.

\item Sergides, M. and Paparoditis, E. (2009). Frequency domain tests of semiparametric hypotheses for locally stationary processes. {\em Scandin. J. Statist.} {\bf 36}, 800--821.

\item Shiraishi, H. and Taniguchi, M. (2007). Statistical estimation of optimal portfolios for locally stationary returns of assets. {\em Int. J. Theor. Appl. Finance.} {\bf 10}, 129--154.

\item Spokoiny, V. (2010). {\em Local parametric methods in nonparametric estimation.} Springer-Verlag, Berlin Heidelberg New York.

\item St$\breve{a}$ric$\breve{a}$, C. and Granger, C. (2005). Nonstationarities in stock returns.  {\em The Review of Economics and Statistics} {\bf 87}, 503--522.

\item Subba Rao, S. (2006). On some nonstationary, nonlinear random processes and their stationary approximations. {\em Advances in Applied Probability} {\bf 38}, 1155--1172.

\item Subba Rao, T. (1970). The fitting of non-stationary time series models with time-dependent parameters. {\em J. Roy Stat Soc B} {\bf 32}, 312--322.

\item Tamaki, K. (2009). Second order properties of locally stationary processes. {\em J. Time Ser. Anal.} {\bf 30}, 145--166.

\item Taniguchi, M., and Kakizawa, Y. (2000). {\em Asymptotic Theory of Statistical Inference for Time Series.} Springer Verlag, New York.

\item Taniguchi, M., Hirukawa, J. and Tamaki, K. (2008). {\em Optimal Statistical Inference in Financial Engineering.} Chapman and Hall/CRC. Boca Raton, Florida.

\item Tj{\o}stheim, D. (1976). Spectral generating operators for non-stationary processes. {\em Adv. Appl. Probab.} {\bf 8}, 831--846.

\item Tunyavetchakit, S. (2010). On the optimal segment length for tapered Yule-Walker estimates for time-varying autoregressive processes. Diploma Thesis, Heidelberg.

\item Van Bellegem, S. and Dahlhaus, R. (2006). Semiparametric estimation by model selection for locally stationary processes. {\em J. Roy. Statist. Soc. B} {\bf 68}, 721--764.

\item Van Bellegem, S. and von Sachs, R. (2004). Forecasting economic time series with unconditional time varying variance. {\em International Journal of Forecasting} {\bf 20}, 611--627.

\item Van Bellegem, S. and von Sachs, R. (2008). Locally adaptive estimation of evolutionary wavelet spectra. {\em Ann. Statist.} {\bf 36}, 1879--1924.

\item Van Bellegem, S. (2011). Locally stationary volatility models, in L. Bauwens, C. Hafner and S. Laurent (eds), Wiley Handbook in Financial Engineering and Econometrics: Volatility Models and Their Applications, Wiley, New York.

\item Vogt, M. (2011). Nonparametric regression for locally stationary time series. Preprint, University of Mannheim.

\item White, H. (1992). {\em Artificial Neural Networks}. Blackwell, Oxford.

\item Whittle, P. (1953). Estimation and information in stationary time series. {\em Ark. Mat.} {\bf 2}, 423--434.

\item Whittle, P. (1954). Some recent contributions to the theory of stationary processes. Appendix to {\em A study in the analysis of stationary time series}, by H. Wold, 2nd ed. 196--228. Almqvist and Wiksell, Uppsala.

\item Wu, W.B. and Zhou, Z. (2011). Gaussian approximations for non-stationary multiple time series. {\em Statistica Sinica} {\bf 21}, 1397--1413.

\item Zani, M. (2002). Large deviations for quadratic forms of locally stationary processes. {\em J. Multivar. Anal.} {\bf 81}, 205--228.

\item Zhou, Z. (2010). Nonparametric inference of quantile curves for nonstationary time series. {\em Ann. Statist.} {\bf 38}, 2187--2217.

\item Zhou, Z. and Wu, W.B. (2009). Local linear quantile estimation for non-stationary time series. {\em Ann. Statist.} {\bf 37}, 2696--2729.

\item Zhou, Z. and Wu, W.B. (2010). Simultaneous inference of linear models with time varying coefficients. {\em J. R. Statist. Soc. B} {\bf 72}, 513--531.
\end{description}
\end{document}